# Differentiation in Topological Vector Spaces

Jinlu Li

Department of Mathematics
Shawnee State University
Portsmouth, Ohio 45662 USA
jli@shawnee.edu

**Abstract.** Differentiation in mathematical analysis is commonly built by using $\varepsilon$-$\delta$-language. This approach also works similarly for defining continuity, Gâteaux (directional) derivative and Fréchet derivative in normed vector spaces, in particular, in Banach spaces, where Fréchet derivatives are defined as limits of ratios with respect to the norms in the considered normed vector spaces.

For general topological vector spaces, if the space is not equipped with a norm, then Fréchet derivatives cannot be similarly defined as in normed vector spaces. However, the topology of every topological vector space can be induced by a family of *F*-seminorms. The cornerstone of this paper is the development of an extended $\varepsilon$-$\delta$-language in general topological vector spaces with respect to these *F*-seminorms. By using this language, we first define the continuity of single-valued mappings. Then we define Gâteaux and Fréchet derivatives as a certain type of limits of ratios with respect to the *F*-seminorms, which generalize Gâteaux and Fréchet derivatives in normed vector spaces.

We will prove the most basic analytic properties of the generalized versions of Gâteaux and Fréchet derivatives, which came up similar to the analytic properties in normed vector spaces. Then we apply them to some specific topological vector spaces that are not normed, such as the Schwartz space (It is locally convex, but it is not normed) and other two spaces that are not even locally convex. For some single-valued mappings defined on these three spaces, we will find the explicit forms of their Gâteaux and Fréchet derivatives.

Finally, we apply the generalized Gâteaux and Fréchet derivatives to solve some vector optimization problems and investigate the order monotonic property of single-valued mappings in general topological vector spaces.

## 1. Introduction

The theory of differentiation lays the firm foundation of modern mathematics, in particular in analysis. In the analysis of Banach spaces, the differentiability of mappings between normed vector spaces (in particular in Banach spaces) has been developed rapidly, which is a natural extension of the classical differentiability of functions in calculus. It is known that the most popular and most useful concepts of differentiation in normed vector spaces are Gâteaux (directional) differentiability and Fréchet differentiability. We review these concepts below. In the following definitions of differentiability, let $(X, \|\cdot\|_X)$ and $(Y, \|\cdot\|_Y)$ be normed vector spaces with origins $\theta_X$ and $\theta_Y$, respectively. Let $T: X \to Y$ be a single-valued mapping. Let $\bar{x} \in X$.

(Gâteaux directional differentiability of *T*). Let $v \in X$ with $v \neq \theta_X$. If there is a vector in *Y* denoted by $T'(\bar{x}, v)$ such that

$$T'(\bar{x}, v) = \lim_{t \to 0} \frac{T(\bar{x}+tv)-T(\bar{x})}{t}, \tag{1.1}$$

then, $T$ is said to be Gâteaux directionally differentiable at point $\bar{x}$ along direction $v$. The point $T'(\bar{x}, v) \in Y$ is called the Gâteaux directional derivative of $T$ at point $\bar{x}$ along direction $v$. Furthermore, if $T$ is Gâteaux directionally differentiable at point $\bar{x}$ along every direction $v \in X \backslash \{\theta_X\}$, then $T$ is said to be Gâteaux differentiable at $\bar{x}$ and the Gâteaux derivative of $T$ at point $\bar{x}$ is denoted by

$$T'(\bar{x})(v) = \lim_{t \to 0} \frac{T(\bar{x}+tv)-T(\bar{x})}{t}, \text{ for any } v \in X \backslash \{\theta_X\}.$$

(Fréchet differentiability of $T$). If there is a continuous and linear mapping $\nabla T(\bar{x}): X \to Y$ such that,

$$\lim_{u \overset{X}{\to} \theta_X} \frac{T(\bar{x}+u)-T(\bar{x})-\nabla T(\bar{x})(u)}{\|u\|_X} = \theta_Y, \tag{1.2}$$

then $T$ is said to be Fréchet differentiable at $\bar{x}$ and $\nabla T(\bar{x})$ is called the Fréchet derivative of $T$ at $\bar{x}$. It is well-known that in normed vector spaces, for single-valued mappings, Fréchet differentiability is stronger than Gâteaux differentiability. That is,

$$T \text{ is Fréchet differentiable at } \bar{x} \implies T \text{ is Gâteaux differentiable at } \bar{x}. \tag{1.3}$$

More precisely speaking, if $T$ is Fréchet differentiable at $\bar{x}$, then $T$ is also Gâteaux differentiable at $\bar{x}$ and

$$T'(\bar{x})(v) = \nabla T(\bar{x})(v), \text{ for any } v \in X \backslash \{\theta_X\}. \tag{1.4}$$

This is the connection between Gâteaux differentiability and Fréchet differentiability in normed vector spaces. Both Gâteaux (directional) differentiability and Fréchet differentiability in normed vector spaces have been applied to many fields in both pure and applied mathematics such as, optimization theory, variational inequalities, control theory, and so forth (see [1−5, 8, 12−13, 15, 18, 20−21, 23−28, 30−33, 37, 39, 41, 47]).

Since the normed vector spaces are equipped with norms, the definition (1.1) and (1.2) can be studied by the $\varepsilon$-$\delta$ language or related techniques used in ordinary limit theory in calculus. And therefore, the above two definitions of differentiability of single-valued mappings between normed vector spaces are natural extension of the classical differentiability of functions in calculus. Notice again, the main reason for these extensions to exist is that the considered normed vector spaces are equipped with norms, which can define the ratios of changes in (1.1) and (1.2) of the considered mappings.

If we study the differentiability of mappings in general topological vector spaces not equipped with norms, we will have the problems that (1.2) cannot be similarly used (see [22, 34−36, 42−46]). And therefore, to define the differentiability of single-valued mappings in general topological vector spaces, it is impossible to similarly use the $\varepsilon$-$\delta$ language or related techniques used in ordinary limit theory in calculus. Without having the limits (1.1) and (1.2), the concept of Fréchet differentiability is difficult to define; even though some differentiability has been implicitly and abstractly defined in general topological vector spaces. For example, a general differentiability has been defined in Lang [22] and Yamamuro [45], which are reviewed below. A real valued function $f$ of a real variable, defined on some neighborhood of 0 is said to be $o(t)$, written $f(t) = o(t)$, if

$$\lim_{t \to 0} \frac{f(t)}{t} = 0.$$

Let $(X, \tau_X)$ and $(Y, \tau_Y)$ be topological vector spaces with origins $\theta_X$ and $\theta_Y$, respectively. Let $\varphi: X \to Y$ be a single-valued mapping. We say that φ is tangent to $\theta_X$, if given a neighborhood $W$ of $\theta_Y$ in $Y$, there is a neighborhood $V$ of $\theta_X$ in $X$ such that

$$\varphi(tV) \subset o(t)W, \text{ for all } t \text{ in a neighborhood of } 0. \tag{1.5}$$

Let $\bar{x} \in X$. Let $T$: $X \to Y$ be a single-valued mapping. $T$ is said to be differentiable at $\bar{x}$ if there is a continuous and linear mapping $\lambda$: $X \to Y$ such that,

$$T(\bar{x} + x) = T(\bar{x}) + \lambda(x) + \varphi(y), \tag{1.6}$$

for small $x$, then φ is tangent to $\theta_X$. In this case, $\lambda$ is unique and it is called the derivative of $T$ at $\bar{x}$.

In this paper, however, to overcome the difficulty of defining the Fréchet differentiability of single-valued mappings in general topological vector spaces not equipped with norms, we will apply the following fundamental theorem of topological vector spaces (see page 35 of Swartz [41] or Jarchow [19] or Zelazko [46]) to solve this problem.

**Theorem 1.1.** *The topology of every topological vector space can be induced by a family of F-seminorms.*

This theorem indeed plays the fundamental roles to defined the continuity and differentiability in topological vector spaces in this paper. Then, for the considered underlying topological vector space, by Theorem 1.1, the topology is considered to be induced by a family of $F$-seminorms. The cornerstone of this paper is by using the family of $F$-seminorms to develop an extended $\varepsilon$-$\delta$ language to define the continuity and differentiability of single-valued mappings in the considered topological vector spaces. We will explain that the generalized differentiability in topological vector spaces is an extension of the differentiability of normed vector spaces. Meanwhile, we also develop related techniques to prove some properties of continuity and differentiability of single-valued mappings. The generalized Gâteaux differentiability (Definition 4.1) and Fréchet differentiability (Definition 4.8) in general topological vector spaces in this paper are different from (1.6). Moreover, we will show that this approach is useful and practical in calculations.

This paper is organized as follows. In section 2, we review the concepts and properties of seminorms and $F$-seminorms on topological vector spaces. We also review the facts that in every topological vector space, the topology can be induced (generated) by a family of $F$-seminorms; and in every locally convex topological vector space, the topology can be induced by a countable family of seminorms. In section 3, we use the extended $\varepsilon$-$\delta$ language and related techniques to define the continuity of mappings between topological vector spaces. Some properties of continuity will be proved in section 3.

In section 4, we also use the extended $\varepsilon$-$\delta$ language and the related techniques to define the Gâteaux (directional) differentiability (Definition 4.1) and Fréchet differentiability (Definition 4.8) of single-valued mappings between topological vector spaces. Some properties of Gâteaux (directional) derivatives and Fréchet derivatives of single-valued mappings will be proved, which are similar to the properties of ordinary derivatives of functions in calculus. For example, under our definition, we prove that both the Gâteaux differentiability and Fréchet differentiability in topological vector spaces are generalizations of Gâteaux differentiability and Fréchet differentiability in normed vector spaces, respectively. However, regarding to the connection between Gâteaux differentiability and Fréchet differentiability, in contrast to the case of normed vector spaces, we cannot similarly prove (1.3) for general topological vector spaces. However, in Theorem 4.18, under some conditions (see condition (4.38)), we will prove that in general topological vector spaces, Fréchet differentiability is stronger than Gâteaux differentiability.

For the purposes of applications of the generalized Gâteaux differentiability and Fréchet differentiability defined in section 4, in sections 5, we concentrate to find the explicit forms of Fréchet derivatives and Gâteaux derivatives of some single-valued mappings such as polynomial type operators on Schwartz space, which is a Fréchet space and it is not a normed vector space. In sections 6 and 7, we concentrate to find the explicit forms of Fréchet derivatives and Gâteaux derivatives of some single-valued mappings respectively in 2 topological vector spaces that are neither locally convex, nor normed.

In section 8, we give some applications of differentiation to vector optimization problems in general topological vector spaces. More precisely speaking, we will study the connection between generalized critical points and ordered extrema; we also investigate the connection between Gâteaux and Fréchet derivatives and order monotonic properties of single-valued mappings in general topological vector spaces.

## 2. Review *F*-seminorm and Seminorm and Their Properties on topological Vector Spaces

Throughout this paper, we use $\mathbb{N}$ to denote the set of all nonnegative integers and let $\mathbb{R}$ and $\mathbb{R}_+$denote the set of real numbers and the set of nonnegative real numbers with the standard topology, respectively. In this section, we review some properties of topological vector spaces. In particularly, we review concepts of *F*-seminorms (the *F* stands for Fréchet) and seminorms on topological vector spaces. These properties play crucial roles in the definitions of continuity and differentiability of mappings between topological vector spaces, in general and in particular, in locally convex topological vector spaces. For more properties of topological vector spaces, for example, see Carothers [6], Jarchow [18], Narici and Beckenstein [35], Swartz [41], Trèves [42], Wilansky [43], Willard [44], and Zelazko [46].

Throughout this section, let $(X, \tau_X)$ be a real topological vector space with origin $\theta_X$ satisfying that both

$$\text{the vector addition: } X \times X \to X \quad \text{and} \quad \text{scalar multiplication operator: } \mathbb{R} \times X \to X \tag{2.1}$$

are continuous mappings (where the domains of these functions are endowed with product topologies of $X$ and $\mathbb{R}$).

Let $p$: $X \to \mathbb{R}_+$ be a real-valued mapping. If $p$ satisfies the following four properties:

(i) $p(x) \geq 0$, for every $x \in X$;
(ii) $p(x + y) \leq p(x) + p(y)$, for all $x, y \in X$;
(iii) $p(ax) \leq p(x)$, for $x \in X$ and every scalar $a$ satisfying $|a| \leq 1$;
(iv) Let $\{a_n\}_{n=1}^{\infty}$ be an arbitrary sequence of positive numbers converging to zero. For every $x \in X, p(a_n x) \to 0$, as $n \to \infty$, then, $p$ is called an *F*-seminorm on *X*. By these properties (i−iv) of *F*-seminorms, we can show the following lemma.

**Lemma 2.1**. *Let p be a F-seminorm on X. Then p also satisfies the following properties*:

(v) $p(\theta_X) = 0$;
(vi) $p(-x) = p(x)$, *for* $x \in X$;
(vii) $p(nx) \leq np(x)$ *and* $p(x) \geq p(\frac{1}{n}x) \geq \frac{1}{n}p(x)$, *for every* $x \in X$ *and any* $n \in \mathbb{N}$;
(viii) *Let a, b be scalars. If* $|a| \leq |b|$, *then* $p(ax) \leq p(bx)$, *for each* $x \in X$;
(ix) *For any* $x \in X$, *if* $p(x) = 0$, *then* $p(ax) = 0$, *for every scalar a.*

*Proof*. We only prove (ix). Let $x \in X$ with $p(x) = 0$. For every given scalar $a$, there is $n \in \mathbb{N}$ such that $|a| \leq n$. Then, by (viii), (vii) and by $p(x) = 0$, we have $p(ax) \leq p(nx) \leq np(x) = 0$. □

An $F$-seminorm $p$ is called an $F$-norm if in addition it satisfies: $p(x) = 0$ implies $x = \theta_X$. An $F$-seminorm $p$ on $X$ is said to be positive (not identical to 0) if there is $x_0 \in X$ such that $p(x_0) > 0$. Since zero identical $F$-seminorm is useless for us, in this paper, we assume that all considered $F$-seminorms are positive. If $X$ is equipped with a $F$-norm $p$ that induces the topology $\tau_X$, then $(X, \tau_X)$ can be rewritten as $(X, p)$ and it is called a $F$-topological vector space, which is simply called a $F$-space. Here, the $F$ stands for Fréchet. However, $F$-space is different from Fréchet space, because Fréchet spaces are defined to be complete locally convex topological vector spaces.

In particular, seminorms on topological vector spaces are considered to be special cases of $F$-seminorms, which is defined as follows. Let $p: X \to \mathbb{R}_+$ be a real-valued mapping. If $p$ satisfies the following four properties:

(i) $p(\theta) = 0$;
(ii) $p(x) \geq 0$, for all $x \in X$;
(iii) $p(sx) = |s|p(x)$, for all $x \in X$ and all scalars $s$;
(iv) $p(x + y) \leq p(x) + p(y)$, for all $x, y \in X$,

then $p$ is called a seminorm on $X$. Note that the above properties (i) and (ii) of seminorms are implied by properties (iii) and (iv).

Let $(X, \tau_X)$ be a real topological vector space. Let $\mathbb{F}_X$ be a family of $F$-seminorms on $X$. Let $\mathcal{F}_X$ be the collection of nonempty finite subsets of $\mathbb{F}_X$, that depends on $X$ and $\mathbb{F}_X$. For any $I \in \mathcal{F}_X$ and any $\lambda > 0$, let

$$U_{I,\lambda} = \cap_{p\in I}\{x \in X: p(x) < \lambda\} = \{x \in X: \max\{p(x): p \in I\} < \lambda\}. \tag{2.2}$$

Note that $U_{I,\lambda}$ depends on $\lambda$, $I$, $X$ and $\mathbb{F}_X$. More generally, for any $x_0 \in X$, $I \in \mathcal{F}_X$ and $\lambda > 0$, we write

$$U_{I,\lambda}(x_0) = \cap_{p\in I}\{x \in X: p(x - x_0) < \lambda\} = \{x \in X: \max\{p(x - x_0): p \in I\} < \lambda\}. \tag{2.3}$$

With respect to the family $\mathbb{F}_X$ of $F$-seminorms on $X$, by (2.2), we write

$$U(\mathbb{F}_X) = \{U_{I,\lambda}: I \in \mathcal{F}_X, \lambda > 0\}. \tag{2.4}$$

Similar to (2.4), by (2.3), for any $x_0 \in X$, we write

$$U(\mathbb{F}_X)(x_0) = \{U_{I,\lambda}(x_0): I \in \mathcal{F}_X, \lambda > 0\}. \tag{2.5}$$

For the above considered topological vector space $(X, \tau_X)$, that the topology $\tau_X$ is induced by the family $\mathbb{F}_X$ of $F$-seminorms on $X$ means that the following two conditions are satisfied:

(I) $U(\mathbb{F}_X)$ defined by (2.4) forms an $\tau_X$-open neighborhood basis of $X$ around $\theta_X$;
(II) For any $x_0 \in X$, $U(\mathbb{F}_X)(x_0)$ defined in (2.5) forms an $\tau_X$-open neighborhood basis of $X$ around $x_0$.

$U(\mathbb{F}_X)$ is called the $F$-seminorm basis of $X$ corresponding to the family $\mathbb{F}_X$ of $F$-seminorms on $X$.

**Definition 2.2**. Let $(X, \tau_X)$ be a topological vector space. If the topology $\tau_X$ of $X$ is induced by a family $\mathbb{F}_X$ of seminorms, then we say that $(X, \tau_X)$ has seminorm construction (or seminorm structure). Or, $(X, \tau_X)$ is said to be seminorm constructed (or seminorm structured). In particular, every locally convex topological vector space is seminorm constructed, in which its topology is induced by a countable family of seminorms.

Next, we rewrite Theorem 1.1 in detail.

**Theorem 2.3.** *Let* $(X, \tau_X)$ *be an arbitrarily given topological vector space. Then, there is a family* $\mathbb{F}_X$ *of F-seminorms on X such that the topology* $\tau_X$ *is induced by* $\mathbb{F}_X$ *satisfying* (I) *and* (II).

We list some properties of $U(\mathbb{F}_X)$ as a lemma (see Theorem 2.9.2 in Jarchow [18], or Swartz [41], or pages 1 and 2 in Zelazko [46]).

**Lemma 2.4**. *Let* $(X, \tau_X)$ *be a topological vector space. Suppose that the topology* $\tau_X$ *is induced by a family* $\mathbb{F}_X$ *of F-seminorms on X. Then,* $U(\mathbb{F}_X)$ *has the following properties.*

(i) *For each* $\lambda > 0, I \in \mathcal{F}_X$, $U_{I,\lambda}$ *is an* $\tau_X$*-open subset of X and it is called an* $\tau_X$*-open neighborhood of* $\theta_X$. *For a given* $\lambda > 0$, *there may be many* $\tau_X$*-open neighborhoods of* $\theta_X$. $U_{I,\lambda}(x_0)$ *is called an* $\tau_X$*-open neighborhood of* $x_0$, *for* $x_0 \in X$;

(ii) *The set* $\{U_{I,\lambda}(x_0): \lambda > 0, I \in \mathcal{F}_X\}$ *forms a filter base on X around point* $x_0$, *for* $x_0 \in X$;

(iii) $U_{I,\lambda_1} \subseteq U_{I,\lambda_2}$, *for* $I \in \mathcal{F}_X$ *and any* $0 < \lambda_1 < \lambda_2$;

(iv) $U_{I,\lambda} \subseteq U_{K,\lambda}$, *for any* $I, K \in \mathcal{F}_X$ *with* $I \supseteq K$ *and for any* $\lambda > 0$;

(v) $U_{I\cup K,\lambda} = U_{I,\lambda} \cap U_{K,\lambda}$, *for any* $I, K \in \mathcal{F}_X$ *and any* $\lambda > 0$;

(vi) $U_{I,\lambda_1} + U_{I,\lambda_2} \subseteq U_{I,\lambda_1+\lambda_2}$, *for* $I \in \mathcal{F}_X$ *and any* $\lambda_1, \lambda_2 > 0$. *In particular, we have*

$$U_{I,\frac{\lambda}{2}} + U_{I,\frac{\lambda}{2}} \subseteq U_{I,\lambda}, \text{ for any } I \in \mathcal{F}_X \text{ and any } \lambda > 0;$$

(vii) $U_{I,\lambda}$ *is a balanced absorbing* $\tau_X$*-open neighborhood of* $\theta_X$.

(viii) $\tau_X$ *is Hausdorff if and only if for every non-zero* $x \in X$, *there exists some* $p \in \mathbb{F}_X$ *such that* $p(x) > 0$, *that is,* $\mathbb{F}_X$ *is separating.*

(ix) *Let* $x \in X$ *and let* $\{x_\alpha\}_{\alpha\in M}$ *be a* *<u>net</u>* *in X. Then,*

$$x_\alpha \to x, \text{ with respect to } \tau_X \quad \Leftrightarrow \quad p(x_\alpha - x) \to 0, \text{ for all } p \in \mathbb{F}_X.$$

*Moreover, if* $\{x_\alpha\}_{\alpha\in M}$ *is* $\tau_X$*-Cauchy in X, then so is the net* $\{p(x_\alpha)\}_{\alpha\in M}$, *for every* $p \in \mathbb{F}_X$.

(x) *We can always take the family* $\mathbb{F}_X$ *of all F-seminorms, which are continuous in the topology* $\tau_X$. *Hence, we always suppose that, for every* $p \in \mathbb{F}_X$,

$$p \text{ is continuous in the topology } \tau_X.$$

*Proof*. Some of the properties have been proved in [18, 41, 46]. Rest of the proof is straight forward and it is omitted here. □

The topological vector space $(X, \tau_X)$, in which the topology $\tau_X$ is induced by a family $\mathbb{F}_X$ *of F*-seminorms on *X* has the following analytic properties.

**Lemma 2.5**. *Let* $(X, \tau_X)$ *be a topological vector space, in which the topology* $\tau_X$ *is induced by a family* $\mathbb{F}_X$ *of F-seminorms on X. Then,*

(i) *For any* $x_0 \in X$, $U(\mathbb{F}_X)$ and $U(\mathbb{F}_X)(x_0)$ *satisfy the following property*

$$U(\mathbb{F}_X)(x_0) = U(\mathbb{F}_X) + x_0.$$

(ii) *For any* $x_0, u \in X$ *and for every* $p \in \mathbb{F}_X$, *we have*

$$\lim_{v \xrightarrow{\tau_X} u} p(x_0 + v) = p(x_0 + u)$$

*In particular*, *we have*

$$\lim_{t \to 0} p(x_0 + tu) = p(x_0),$$

*and* $$\lim_{t \to 0} p(tu) = 0, \text{ for any } u \in X.$$

(iii) *For any* $u \in X \backslash \{\theta_X\}$ *and for every* $p \in \mathbb{F}_X$, $p(tu)$ *is a continuous and increasing function of* $t$ *on* $\mathbb{R}$ *satisfying* $p(-tu) = p(tu)$, *for any* $t \in \mathbb{R}$.

*In particular*, *if X is seminotm constructed*, *then for any* $p \in \mathbb{F}_X$, $p$ *is a convex functional such that for any* $u, v \in X$, *we have*

$$p(t_1 u + t_2 v) \leq |t_1| p(u) + |t_2| p(v), \text{ for any } t_1, t_2 \in \mathbb{R}.$$

*Proof.* The proof of this lemma follows from (2.1), part (x) in Lemma 2.4 and part (v) in Lemma 2.1 and it is omitted here. □

Next, we consider a special case of topological vector spaces, in which the topologies are induce by countable (including finite) families of *F*-seminorms.

**Lemma 2.6**. *Let* $(X, \tau_X)$ *be a Hausdorff topological vector space*, *in which the topology* $\tau_X$ *is induced by a countable family* $\mathbb{F}_X$ *of F-seminorms on X with* $\mathbb{F}_X = \{p_i : i \in \mathbb{M}_X\}$, *where* $\mathbb{M}_X \subseteq \mathbb{N}$. *Let* $\{a_i : i \in \mathbb{M}_X\}$ *be a set of positive numbers with* $\sum_{i \in \mathbb{M}_X} a_i < \infty$. *Define* $p_X : X \to \mathbb{R}_+$ *by*

$$p_X(x) = \sum_{i \in \mathbb{M}_X} a_i \frac{p_i(x)}{1 + p_i(x)}, \text{ for any } x \in X. \tag{2.6}$$

*Then*, *we have*

(i) $p_X : X \to \mathbb{R}_+$ *defined in* (2.6) *is a F-norm on X*;
(ii) $(X, \tau_X)$ *is a F-space*, *in which the topology* $\tau_X$ *is induced by this F-norm* $p_X$ *on X*;

*Proof.* The proof of this lemma is straight forward and it is omitted here. □

## 3. Continuity of Mappings between Topological Vector Spaces

### 3.1. Continuity with respect to *F*-seminorm Bases

Let $(Y, \tau_Y)$ be a topological vector space with the topology $\tau_Y$ induced by a family $\mathbb{F}_Y$ of *F*-seminorms on *Y*. Let $\mathcal{F}_Y$ be the collection of nonempty finite subsets of $\mathbb{F}_Y$. Similarly to (2.1) −(2.3), we have notations for *Y*. For $J \in \mathcal{F}_Y$ and $\lambda > 0$, let

$$V_{J,\lambda} = \cap_{q \in J} \{y \in Y : q(y) < \lambda\} = \{y \in Y : \max\{q(y) : q \in J\} < \lambda\}. \tag{3.1}$$

For any $y_0 \in Y$, $J \in \mathcal{F}_Y$ and $\lambda > 0$, we write

$$V_{J,\lambda}(y_0) = \cap_{q\in J}\{y \in Y: q(y - y_0) < \lambda\}. \tag{3.2}$$

Let $V(\mathbb{F}_Y) = \{V_{J,\lambda}: \lambda > 0, J \in \mathcal{F}_Y\}$ be the $F$-seminorm basis of $\tau_Y$ corresponding to the family $\mathbb{F}_Y$ of $F$-seminorms. By following the ordinary definition of continuity in calculus, we give the definition of the continuity of mappings between topological vector spaces with respect to their $F$-seminorm bases.

In this section, if not stated otherwise, we always let $(X, \tau_X)$ and $(Y, \tau_Y)$ be topological vector spaces, in which the topologies $\tau_X$ and $\tau_Y$ are induced by families $\mathbb{F}_X$ and $\mathbb{F}_Y$ of $F$-seminorms, respectively. Let $D$ be a nonempty $\tau_X$-open subset of $X$. Let $T: D \to Y$ be a single-valued mapping. By an extended $\varepsilon$-$\delta$ language in topological vector spaces equipped with families of $F$-seminorms, we define the continuity of single-valued mappings.

**Definition 3.1**. Let $x_0 \in D$ and $y_0 \in Y$ with $T(x_0) = y_0$. $T$ is said to be continuous at point $x_0$, if and only if, for any $J \in \mathcal{F}_Y$ and $\varepsilon > 0$, there are $I \in \mathcal{F}_X$ and $\delta > 0$ such that, for $x \in D$, we have

$$\max\{p(x - x_0): p \in I\} < \delta \quad \Longrightarrow \quad \max\{q(T(x) - T(x_0)): q \in J\} < \varepsilon. \tag{3.3}$$

In terms of the $F$-seminorm bases in $X$ and $Y$, Definition 3.1 can be rewritten as: $T$ is continuous at point $x_0$, if and only if, for any $\tau_Y$-open neighborhood $V_{J,\varepsilon}(y_0)$ of $y_0$ in $Y$, there is an $\tau_X$-open neighborhood $U_{I,\delta}(x_0)$ of $x_0$ in $X$, for some $\delta > 0$ and $I \in \mathcal{F}_X$ with $U_{I,\delta}(x_0) \subseteq D$ such that

$$T(U_{I,\delta}(x_0)) \subseteq V_{J,\varepsilon}(y_0). \tag{3.4}$$

In this case, we write

$$\lim_{x\to x_0} T(x) = T(x_0) = y_0.$$

Based on the extended $\varepsilon$-$\delta$ language (3.3) and the techniques for calculating limits in calculus, the related techniques for studying the continuity of mappings between topological vector spaces can be developed. We prove some properties of continuity defined in (3.3) and (3.4).

**Lemma 3.2**. *Let $x_0 \in D$ and let $a$ be a scalar. Suppose that $T$ is continuous at $x_0$. Then $aT: D \to Y$ is also continuous at $x_0$.*

*Proof.* For $x_0 \in D$, let $y_0 \in Y$ with $T(x_0) = y_0$. It is clear that this lemma holds for $a = 0$. Next, we prove this lemma for $a \neq 0$. For this given $a$, there is a positive integer $n$ such that $|a| \leq n$. Let $V_{J,\varepsilon}(ay_0)$ be any given $\tau_Y$-open neighborhood of $ay_0$ in $Y$. Then, $V_{J,\frac{\varepsilon}{n}}(y_0)$ is an $\tau_Y$-open neighborhood of $y_0$ in $Y$. By the assumption that $T$ is continuous at $x_0$, there are $\delta > 0$ and $I \in \mathcal{F}_X$ such that

$$\max\{p(x - x_0): p \in I\} < \delta \quad \Longrightarrow \quad \max\{q(T(x) - y_0): q \in J\} < \frac{\varepsilon}{n}, \text{ for } x \in D. \tag{3.5}$$

By properties (vii) and (iii) of $F$-seminorms and by $|a| \leq n$, we have

$$p(ax) \leq p(nx) \leq np(x), \text{ for every } x \in D.$$

For $x \in D$, by (3.5) and the above inequalities, we have

$$\max\{p(x - x_0): p \in I\} < \delta$$

$$\Longrightarrow \quad \max\{q(aT(x) - ay_0): q \in J\}$$

$$= \max\{q(a(T(x) - y_0)): q \in J\}$$

$$\leq \max\{q(n(T(x) - y_0)): q \in J\}$$

$$\leq n\max\{q(T(x) - y_0): q \in J\}$$

$$< n\frac{\varepsilon}{n} = \varepsilon.$$

This proves this lemma. □

**Lemma 3.3**. *Let $T_1, T_2: D \to Y$ be single-valued mappings. Let $x_0 \in D$. Suppose that both $T_1$ and $T_2$ are continuous at $x_0$. Then $T_1 + T_2$ is continuous at $x_0$.*

*Proof.* Suppose that $T_1(x_0) = y_1$ and $T_2(x_0) = y_2$, for $y_1, y_2 \in Y$. Let $V_{J,\varepsilon}(y_1 + y_2)$ be an arbitrarily given $\tau_Y$-open neighborhood of $y_1 + y_2$ in *Y*. For these arbitrarily given $J \in \mathcal{F}$ and $\varepsilon > 0$, $V_{J,\frac{\varepsilon}{2}}(y_1)$ is an $\tau_Y$-open neighborhood of $y_1$ and $V_{J,\frac{\varepsilon}{2}}(y_2)$ is an $\tau_Y$-open neighborhood of $y_2$ in *Y*. By the assumption that, for $k = 1, 2$, $T_k$ is continuous at $x_0$, for this $\tau_Y$-open neighborhood $V_{J,\frac{\varepsilon}{2}}(y_k)$ of $y_k$ in *Y*, there are an $\delta_k > 0$ and $I_k \in \mathcal{F}_X$ such that, for $x \in D$, by (3.3), we have

$$\max\{p(x - x_0): p \in I_k\} < \delta_k \quad \Longrightarrow \quad \max\{q(T_k(x) - y_k): q \in J\} < \frac{\varepsilon}{2}, \text{ for } k = 1, 2. \quad (3.6)$$

Let $\delta = \min\{\delta_1, \delta_2\}$ and let $I = I_1 \cup I_2$. It is clear that $I \in \mathcal{F}_X$ Then, for $x \in X$, we have

$$\max\{p(x - x_0): p \in I\} < \delta \Longrightarrow \max\{p(x - x_0): p \in I_k\} < \delta \leq \delta_k, \text{ for } k = 1, 2. \quad (3.7)$$

By (3.6) and (3.7), we have

$$\max\{p(x - x_0): p \in I\} < \delta$$

$$\Longrightarrow \quad \max\{q\big((T_1 + T_2)(x) - (y_1 + y_2)\big): q \in J\}$$

$$= \max\{q\big((T_1(x) - y_1) + ((T_2(x) - y_2))\big): q \in J\}$$

$$\leq \max\{q(T_1(x) - y_1): q \in J\} + \max\{q(T_2(x) - y_2): q \in J\}$$

$$< \frac{\varepsilon}{2} + \frac{\varepsilon}{2} = \varepsilon. \quad (3.8)$$

By (3.8), this lemma is proved. □

By Lemmas 3.4 and 3.5, we immediately obtain the following linearity of continuity of mappings on topological vector spaces.

**Theorem 3.4**. *Let $T_1, T_2: D \to Y$ be single-valued mappings. Let $x_0 \in D$ and let $a_1, a_2$ be scalars. Suppose that both $T_1$ and $T_2$ are continuous at $x_0$. Then $a_1T_1 + a_2T_2: X \to Y$ is also continuous at $x_0$.*

Let $(Z, \tau_Z)$ be another topological vector space, in which the topology $\tau_Z$ is induced by a family $\mathbb{F}_Z$ of *F*-seminorms $r$'s on *Z*. We write $\mathbb{F}_Z = \{r \in \mathbb{F}_Z\}$. The *F*-seminorm basis of *Z* corresponding to the family $\mathbb{F}_Z$ is $\{W_{K,\lambda}: \lambda > 0, K \in \mathcal{F}_X\}$, in which, for any finite subset $K \in \mathcal{F}_X$ and any $\lambda > 0$, we have

$$W_{K,\lambda} = \cap_{k\in K}\{z \in Z: r_k(z) < \lambda\}.$$

**Lemma 3.5**. *Let $T: D \to Y$ and $S: Y \to Z$ be single-valued mappings. Let $x_0 \in X$, $y_0 \in Y$ and $z_0 \in Z$ with $T(x_0) = y_0$ and $S(y_0) = z_0$. Suppose that T is continuous at $x_0$ and S is continuous at $y_0$. Then, the composite mapping $S \circ T: D \to Z$ is continuous at $x_0$.*

*Proof.* Let $W_{K,\varepsilon}(z_0)$ be an arbitrarily given $\varepsilon$-open neighborhood of $z_0$ in Z satisfying $(S \circ T)(x_0) = z_0$. By the assumption that $S: Y \to Z$ is continuous at $y_0 = T(x_0)$ in Y, there is an $\tau_Y$-open neighborhood $V_{J,\gamma}(y_0)$ of $y_0$ in Y, for some $\gamma > 0$ and $J \in \mathcal{F}_Y$ such that, for any $y \in V_{J,\gamma}(y_0)$, we have

$$\max\{q(y - y_0): q \in J\} < \gamma \quad \Longrightarrow \quad \max\{r(S(y) - z_0): r \in K\} < \varepsilon. \tag{3.9}$$

For this fixed $\tau_Y$-open neighborhood $V_{J,\gamma}(y_0)$ of $y_0$ in Y, by the assumption that $T: D \to Y$ is continuous at $x_0$ in X, there is an $\tau_X$-open neighborhood $U_{I,\delta}(x_0)$ of $x_0$ in X, for some $\delta > 0$ and $I \in \mathcal{F}_X$ such that, for $x \in U_{I,\delta}(x_0)$, by (3.3), we have

$$\max\{p(x - x_0): p \in I\} < \delta \quad \Longrightarrow \quad \max\{q(T(x) - y_0): q \in J\} < \gamma, \text{ for } x \in X. \tag{3.10}$$

By combining (3.10) and (3.9), for the arbitrarily given $\tau_Z$-open neighborhood $W_{K,\varepsilon}(z_0)$ of $z_0$ in Z, for $x \in U_{I,\delta}(x_0)$, we have

$$\max\{p(x - x_0): p \in I\} < \delta \quad \Longrightarrow \quad \max\{q(T(x) - y_0): q \in J\} < \gamma$$

$$\Longrightarrow \quad \max\{r(S(T(x)) - z_0): r \in K\} < \varepsilon. \qquad \square$$

Next, we consider a special case of topological vector spaces, in which the topologies are induced by countable (including finite) families of *F*-seminorms. Then, in this special case, the continuity of mappings between topological vector spaces defined in Definition (3.1) can be defined by the ordinary $\varepsilon$-$\delta$ formulations with respect to *F*-norms.

**Lemma 3.6**. *Let* $(X, \tau_X)$ *and* $(Y, \tau_Y)$ *be Hausdorff topological vector spaces, in which the topologies $\tau_X$ and $\tau_Y$ are induced by countable families $\mathbb{F}_X = \{p_i: i \in \mathbb{M}_X\}$ and $\mathbb{F}_Y = \{q_j: j \in \mathbb{M}_Y\}$ of F-seminorms, respectively, in which, $\mathbb{M}_X, \mathbb{M}_Y \subseteq \mathbb{N}$. Let* $\{a_i: i \in \mathbb{M}_X\}$ and $\{b_j: j \in \mathbb{M}_Y\}$ *be sets of positive numbers with $\sum_{i\in\mathbb{M}_X} a_i < \infty$ and $\sum_{j\in\mathbb{M}_Y} b_j < \infty$. By* (2.6), *we define the F-norms $p_X: X \to \mathbb{R}_+$ and $p_Y: Y \to \mathbb{R}_+$ by*

$$p_X(x) = \sum_{i\in\mathbb{M}_X} a_i \frac{p_i(x)}{1+p_i(x)}, \text{for any } x \in X,$$

*and*

$$q_Y(y) = \sum_{j\in\mathbb{M}_Y} b_j \frac{q_j(y)}{1+q_j(y)}, \text{for any } y \in Y.$$

*Let $T: D \to Y$ be a single-valued mapping. Let $x_0 \in D$ and $y_0 \in Y$ with $T(x_0) = y_0$. Then, T is continuous at $x_0$ under Definition* 3.1 *if and only if for any given $\varepsilon > 0$, there is a $\delta > 0$, such that*

$$p_X(x - x_0) < \delta \quad \Longrightarrow \quad q_Y(T(x) - y_0) < \varepsilon, \text{for } x \in D. \tag{3.11}$$

*Proof.* This lemma can be proved by using the results of Lemma 2.6, which implies that the topologies $\tau_X$ and $\tau_X$ are induced by the *F*-norm $p_X$ and $p_Y$, respectively. Since we omitted the proof of Lemma 2.6, so we give a directly proof of this lemma.

" $\Longrightarrow$ " Suppose that *T* is continuous at $x_0$ under Definition 3.1. We want to show that (3.11) holds. For any given $\varepsilon > 0$, we take an arbitrary $J \in \mathcal{F}_Y$ with $J = \{q_{j_1}, q_{j_2}, \ldots, q_{j_m}\}$, for some positive integer *m*.

For this given $\varepsilon > 0$, there is a positive integer $M$ such that

$$\sum_{j\in\mathbb{M}_Y, j>M} b_j < \frac{\varepsilon}{2}. \tag{3.12}$$

Let $B = \sum_{j\in\mathbb{M}_Y} b_j$ and let $\varepsilon_1 = \frac{\varepsilon}{2B}$. Then, we have an $\tau_Y$-open neighborhood $V_{J,\varepsilon_1}(y_0)$ of $y_0$ in $Y$. By the assumption, there is an $\tau_X$-open neighborhood $U_{I,\delta_1}(x_0)$ of $x_0$ in $X$, for some $\delta_1 > 0$ and $I \in \mathcal{F}_X$ with $J = \{p_{i_1}, p_{i_2}, \ldots, p_{i_n}\}$, for some positive integer $n$, such that, for $x \in U_{I,\delta_1}(x_0)$, by (3.3), we have

$$\max\{p(x - x_0): p \in I\} < \delta_1 \quad \Longrightarrow \quad \max\{q(T(x) - y_0): q \in J\} < \varepsilon_1, \text{ for } x \in D. \tag{3.13}$$

Let $a = \min\{a_{i_1}, a_{i_2}, \ldots, a_{i_n}\}$. Let $\delta = \frac{a\delta_1}{1+\delta_1}$. Then, for any $x \in D$, we have

$$p_X(x - x_0) = \sum_{i\in\mathbb{M}_X} a_i \frac{p_i(x-x_0)}{1+p_i(x-x_0)} < \delta$$

$$\Longrightarrow \quad \frac{p_i(x-x_0)}{1+p_i(x-x_0)} < \frac{\delta}{a_i} \leq \frac{\delta}{a} = \frac{\delta_1}{1+\delta_1}, \text{ for any } i \in \mathbb{M}_X. \tag{3.14}$$

By the strictly increasing property of the function $\frac{t}{1+t}$, for $t \in (0, \infty)$, (3.14) and (3.13) imply that

$$\sum_{i\in\mathbb{M}_X} a_i \frac{p_i(x-x_0)}{1+p_i(x-x_0)} < \delta \quad \Longrightarrow \quad p_i(x - x_0) < \delta_1, \text{ for each } i \in \mathbb{M}_X$$

$$\Longrightarrow \quad \max\{q(T(x) - y_0): q \in J\} < \varepsilon_1. \tag{3.15}$$

This implies that if

$$p_X(x - x_0) = \sum_{i\in\mathbb{M}_X} a_i \frac{p_i(x-x_0)}{1+p_i(x-x_0)} < \delta,$$

then, by $B = \sum_{j\in\mathbb{M}_Y} b_j$, $\varepsilon_1 = \frac{\varepsilon}{2B}$, (3.12) and (3.15), we have

$$\begin{aligned}
&q_Y(T(x) - y_0) \\
&= \sum_{j\in\mathbb{M}_Y} b_j \frac{q_j(T(x)-y_0)}{1+q_j(T(x)-y_0)} \\
&= \sum_{j\in\mathbb{M}_Y, j\leq M} b_j \frac{q_j(T(x)-y_0)}{1+q_j(T(x)-y_0)} + \sum_{j\in\mathbb{M}_Y, j>M} b_j \frac{q_j(T(x)-y_0)}{1+q_j(T(x)-y_0)} \\
&< \sum_{j\in\mathbb{M}_Y, j\leq M} b_j \frac{\varepsilon_1}{1+\varepsilon_1} + \sum_{j\in\mathbb{M}_Y, j>M} b_j \\
&< \sum_{j\in\mathbb{M}_Y, j\leq M} b_j \frac{\varepsilon}{2B} + \frac{\varepsilon}{2} \\
&= \varepsilon.
\end{aligned}$$

This proves the part " $\Longrightarrow$ ". Next, we prove the part " $\Longleftarrow$ ". Suppose that (3.11) holds. We want to prove that $T$ is continuous at $x_0$ under Definition 3.1. Take arbitrarily an $\tau_Y$-open neighborhood $V_{J,\varepsilon}(y_0)$ of $y_0$ in $Y$ with $J = \{q_{j_1}, q_{j_2}, \ldots, q_{j_m}\} \in \mathcal{F}_Y$, for some positive integer $m$ and $\varepsilon > 0$. Let $b = \min\{b_{j_1}, b_{j_2}, \ldots, b_{j_m}\}$. Let $\varepsilon_1 = \frac{b\varepsilon}{1+\varepsilon}$. By (3.11), there is a $\delta_1 > 0$, such that, for $x \in D$, we have

$$p_X(x - x_0) < \delta_1 \qquad \Longrightarrow \qquad q_Y(T(x) - y_0) < \varepsilon_1. \tag{3.16}$$

For this given $\delta_1 > 0$, there is a positive integer $N$ such that

$$\sum_{i \in \mathbb{M}_X, i \le N} a_i > 0 \quad \text{and} \quad \sum_{i \in \mathbb{M}_X, i > N} a_i < \frac{\delta_1}{2}. \tag{3.17}$$

Let $A = \sum_{i \in \mathbb{M}_X} a_i$ and let $\delta = \frac{\delta_1}{2A}$. Let $I = \{p_i \colon i \in \mathbb{M}_X \text{ and } i \le N\}$. By the first inequality in (3.17), we see that $I$ is finite and $I \neq \emptyset$. Hence, $I \in \mathcal{F}_X$. So, we have an $\tau_X$-open neighborhood $U_{I,\delta}(x_0)$ of $x_0$ in $X$, for some $\delta = \frac{\delta_1}{2B}$ and $I \in \mathcal{F}_X$ with $I = \{p_i \colon i \in \mathbb{M}_X \text{ and } i \le N\}$. Then, for any $x \in D$, we have that

$$x \in U_{I,\delta}(x_0) \ \Longrightarrow p_i(x - x_0) < \delta, \text{ for any } i \in \mathbb{M}_X \text{ and } i \le N. \tag{3.18}$$

This implies that if $x \in U_{I,\delta}(x_0)$, then, by $A = \sum_{i \in \mathbb{M}_X} a_i$ and let $\delta = \frac{\delta_1}{2A}$, by (3.17) and (3.18), we have

$$\begin{aligned}
& p_X(x - x_0) \\
&= \sum_{i \in \mathbb{M}_X} a_i \frac{p_i(x-x_0)}{1+p_i(x-x_0)} \\
&= \sum_{i \in \mathbb{M}_X, i \le N} a_i \frac{p_i(x-x_0)}{1+p_i(x-x_0)} + \sum_{i \in \mathbb{M}_X, i > N} a_i \frac{p_i(x-x_0)}{1+p_i(x-x_0)} \\
&< \sum_{i \in \mathbb{M}_X, i \le N} a_i \, p_i(x - x_0) + \sum_{i \in \mathbb{M}_X, i > N} a_i \\
&< \sum_{i \in \mathbb{M}_X, i \le N} a_i \, \delta + \frac{\delta_1}{2} \\
&= \sum_{i \in \mathbb{M}_X, i \le N} a_i \frac{\delta_1}{2A} + \frac{\delta_1}{2} \\
&< \delta_1.
\end{aligned} \tag{3.19}$$

By (3.19) and (3.16), we have that

$$x \in U_{I,\delta}(x_0) \ \Longrightarrow \ \sum_{j \in \mathbb{M}_Y} b_j \frac{q_j(T(x)-y_0)}{1+q_j(T(x)-y_0)} < \varepsilon_1.$$

This implies that

$$x \in U_{I,\delta}(x_0) \ \Longrightarrow \ b_j \frac{q_j(T(x)-y_0)}{1+q_j(T(x)-y_0)} < \varepsilon_1, \text{ for any } q_j \in J. \tag{3.20}$$

Since $\varepsilon_1 = \frac{b\varepsilon}{1+\varepsilon}$, by (3.20), we have

$$x \in U_{I,\delta}(x_0) \ \Longrightarrow \ b_j \frac{q_j(T(x)-y_0)}{1+q_j(T(x)-y_0)} < \frac{b\varepsilon}{1+\varepsilon}, \text{ for any } q_j \in J.$$

By $b_j \ge b$, for any $q_j \in J$, this implies that

$$x \in U_{I,\delta}(x_0) \ \Longrightarrow \ \frac{q_j(T(x)-y_0)}{1+q_j(T(x)-y_0)} < \frac{\varepsilon}{1+\varepsilon}, \text{ for any } q_j \in J. \tag{3.21}$$

By the strictly increasing property of the function $\frac{t}{1+t}$, for $t \in (0, \infty)$, (3.21) implies that

$$x \in U_{I,\delta}(x_0) \implies q_j(T(x) - y_0) < \varepsilon, \text{ for any } q_j \in J.$$

This proves that $T$ is continuous at $x_0$ under Definition 3.1. □

### 3.2. Boundness and Continuity of Linear Operators in Topological Vector Spaces

To start this subsection, we first review the concepts of bounded subsets in topological vector spaces and bounded mappings between topological vector spaces. A subset $A$ of a topological vector space $X$ is said to be bounded if for every open neighborhood $U$ of the origin $\theta_X$, there exists $t > 0$ such that $A \subseteq tU$.

Let $(X, \tau_X)$ and $(Y, \tau_Y)$ be topological vector spaces. Let $T: X \to Y$ be a single-valued mapping. If $T$ maps every bounded subset $A$ in $X$ to a bounded subset $T(A)$ in $Y$, then $T$ is said to be bounded. If there is an open neighborhood $U$ of the origin $\theta_X$ such that $T$ maps $U$ to a bounded subset $T(U)$ in $Y$, then $T$ is said to be bounded on a neighborhood. The following results in Theorem 3.7 below are well-known and they are important in analysis in topological vector spaces.

Let $(X, \tau_X)$ be a topological vector space, in which $\tau_X$ is induced by a family $\mathbb{F}_X$ of $F$-seminorms on $X$. We can construct some examples to show that, for $\delta > 0$ and $I \in \mathcal{F}_X$, the $\tau_X$-open neighborhood $U_{I,\delta}$ of $\theta_X$ is not necessarily bounded with respect to the topology $\tau_X$ induced by the $F$-seminorm basis $U(\mathbb{F}_X)$.

**Theorem 3.7 (see** Narici and Beckenstein [35] and Wilansky [43]**).**

(i) *Every continuous linear operator between topological vector spaces is always a bounded linear operator.*
(ii) *For any linear operator, if it is bounded on a neighborhood, then it is continuous.*
(iii) *If X is a bornological space (for example, a pseudometrizable topological vector space, a Fréchet space, a normed space) and Y is a locally convex topological vector space, then a linear operator from X to Y is bounded if and only if it is continuous.*
(iv) *If both X and Y are locally convex and $T: X \to Y$ is a continuous linear operator, then, for every continuous seminorm q on Y, there exists a continuous seminorm p on X such that*

$$(q \circ T)(x) \leq p(x), \textit{ for any } x \in X.$$

(v) *Let X and Y be locally convex topological vector spaces whose topologies are defined by countable families of seminorms $\{p_i: i \in \mathbb{N}_0\}$ and $\{q_j: j \in \mathbb{N}_0\}$, respectively. Then, a linear mapping $T: X \to Y$ is continuous when for every $j \in \mathbb{N}_0$ there exist a finite family $p_{i_1}, \ldots, p_{i_n}$ and a positive number $C_j$ such that*

$$q_j(Tx) \leq C_j[p_{i_1}(x) + \cdot \cdot \cdot + p_{i_n}(x)], \textit{ for any } x \in X.$$

## 4. Differentiation in Topological Vector Spaces Equipped with *F*-seminorm Bases

### 4.1. Gâteaux Directional Differentiability and Gâteaux Differentiability

In this subsection, we define the Gâteaux directional differentiability and Gâteaux differentiability of single-valued mappings in general topological vector spaces. We prove some analytic properties of Gâteaux differentiability.

To define Gâteaux directional differentiability and Gâteaux differentiability of single-valued mappings,

the $F$-seminorm construction of the domain space is not necessarily mentioned. Throughout this subsection, if not otherwise stated, we let $(X, \tau_X)$ and $(Y, \tau_Y)$ be Hausdorff topological vector spaces, in which the topology $\tau_Y$ is induced by a family $\mathbb{F}_Y$ of positive $F$-seminorms on $Y$. The family $\mathbb{F}_Y$ generates the corresponding $F$-seminorm basis $\{V_{J,\lambda}: \lambda > 0, J \in \mathcal{F}_Y\}$ of $\tau_Y$.

In this section, if not otherwise stated, we always let $D$ be an $\tau_X$-open nonempty subset of $X$ and let $T$: $D \to Y$ be a single-valued mapping.

**Definition 4.1**. Let $\bar{x} \in D$ and $v \in X \backslash \{\theta_X\}$. Define $T'(\bar{x}, v)$ to be a value from $Y$ as follows. If, for any $J \in \mathcal{F}_Y$ and for any given $\varepsilon > 0$, there is $\delta > 0$ such that, for any real number $t$, we have

$$0 < |t| < \delta \quad \Longrightarrow \quad \max\left\{q\left(\frac{T(\bar{x}+tv)-T(\bar{x})-tT'(\bar{x},v)}{t}\right): q \in J\right\} < \varepsilon, \tag{4.1}$$

then, $T$ is said to be Gâteaux directionally differentiable at point $\bar{x}$ along direction $v$ and $T'(\bar{x}, v)$ is called the Gâteaux directional derivative of $T$ at $\bar{x}$ along direction $v$. Here, we always assume that $\delta$ is taken to be small enough such that if $|t| < \delta$, then $\bar{x} + tv \in D$. This ensures that $T(\bar{x} + tv)$ is defined.

Then, Definition 4.1 is equivalently rewritten using the following limit: $T'(\bar{x}, v)$ is the Gâteaux directional derivative of $T$ at $\bar{x}$ along direction $v$, if and only if, for any $J \in \mathcal{F}_Y$, we have

$$\lim_{t \to 0} \max\left\{q\left(\frac{T(\bar{x}+tv)-T(\bar{x})-tT'(\bar{x},v)}{t}\right): q \in J\right\} = 0. \tag{4.2}$$

If $T$ is Gâteaux directionally differentiable at point $\bar{x}$ along any direction $v \in X \backslash \{\theta_X\}$, then $T$ is said to be Gâteaux differentiable at point $\bar{x}$. In this case, the Gâteaux derivative of $T$ at $\bar{x}$ is denoted by $T'(\bar{x})$ and defined as

$$T'(\bar{x})(v) = T'(\bar{x}, v), \text{ for any } v \in X \backslash \{\theta_X\}.$$

**Theorem 4.2**. *Let $T$: $D \to Y$ be a single-valued mapping. Let $\bar{x} \in D$ and $v \in X \backslash \{\theta_X\}$. Suppose that $T$ is Gâteaux directionally differentiable at $\bar{x}$ along direction $v$. Then the Gâteaux directional derivative $T'(\bar{x}, v)$ of $T$ at $\bar{x}$ along direction $v$ is unique*.

*Proof*. Suppose that there are two points $T_1'(\bar{x}, v)$ and $T_2'(\bar{x}, v)$ such that both satisfy (4.1). Let $J \in \mathcal{F}_Y$ be arbitrarily given. Take an arbitrary $\varepsilon > 0$. For the given $J \in \mathcal{F}_Y$ and $\frac{\varepsilon}{2} > 0$, there is $\delta > 0$ such that, for any real number $t$,

$$0 < |t| < \delta \quad \Longrightarrow \quad \max\left\{q\left(\frac{T(\bar{x}+tv)-T(\bar{x})-tT_k'(\bar{x},v)}{t}\right): q \in J\right\} < \frac{\varepsilon}{2}, \text{ for } k = 1, 2.$$

This implies that

$$\max\left\{q\left(T_1'(\bar{x}, v) - T_2'(\bar{x}, v)\right): q \in J\right\}$$

$$\leq \max\left\{q\left(\frac{T(\bar{x}+tv)-T(\bar{x})-tT_1'(\bar{x},v)}{t}\right): q \in J\right\} + \max\left\{q\left(\frac{T(\bar{x}+tv)-T(\bar{x})-tT_2'(\bar{x},v)}{t}\right): q \in J\right\} < \varepsilon.$$

Then we obtain

$$q\left(T_1'(\bar{x}, v) - T_2'(\bar{x}, v)\right) = 0, \text{ for any } q \in J.$$

Recall that $J \in \mathcal{F}_Y$ is arbitrarily given, this yields that

$$q\big(T_1'(\bar{x},v) - T_2'(\bar{x},v)\big) = 0, \text{ for any } q \in \mathbb{F}_Y. \tag{4.3}$$

Since $Y$ is a Hausdorff topological vector space, (4.3) implies that $T_1'(\bar{x},v) - T_2'(\bar{x},v) = \theta_Y$. □

**Corollary 4.3**. *Let $\bar{x} \in D$. Suppose that $T$ is Gâteaux differentiable at $\bar{x}$. Then the Gâteaux derivative $T'(\bar{x})$ of $T$ at $\bar{x}$ is unique.*

*Proof*. This corollary follows from Theorem 4.2 immediately. □

Now we prove some properties of Gâteaux derivatives.

**Lemma 4.4**. *Let $\bar{x} \in D$ and $v \in X\backslash\{\theta_X\}$ and let $a$ be a scalar. Suppose that $T$ is Gâteaux directionally differentiable at $\bar{x}$ along direction $v$. Then $aT$: $D \rightarrow Y$ is also Gâteaux directionally differentiable at $\bar{x}$ along direction $v$ and*

$$(aT)'(\bar{x},v) = aT'(\bar{x},v). \tag{4.4}$$

*Proof.* It is clear that this lemma holds for $a = 0$. So, we only prove (4.4) for $a \neq 0$. For the given $a \neq 0$, there is a positive integer $m$ such that $|a| \leq m$. Let $J \in \mathcal{F}_Y$ and $\varepsilon > 0$ be arbitrarily given. For this given $J \in \mathcal{F}_Y$ and $\frac{\varepsilon}{2} > 0$, by the assumption that $T$ is Gâteaux differentiable at $\bar{x}$ along direction $v$, there is $\delta > 0$ such that, for any real number $t$, we have

$$0 < |t| < \delta \quad \Longrightarrow \quad \max\Big\{q\left(\tfrac{T(\bar{x}+tv)-T(\bar{x})-tT'(\bar{x},v)}{t}\right): q \in J\Big\} < \tfrac{\varepsilon}{m}. \tag{4.5}$$

By property (vii) of $F$-seminorms in Lemma 2.1 and by (4.5), we have that if $0 < |t| < \delta$, then

$$\max\Big\{q\left(\tfrac{(aT)(\bar{x}+tv)-(aT)(\bar{x})-t(aT')(\bar{x},v)}{t}\right): q \in J\Big\}$$

$$= \max\Big\{q\left(\tfrac{a(T(\bar{x}+tv)-T(\bar{x})-tT'(\bar{x},v))}{t}\right): q \in J\Big\}$$

$$\leq \max\Big\{q\left(\tfrac{m(T(\bar{x}+tv)-T(\bar{x})-tT'(\bar{x},v))}{t}\right): q \in J\Big\}$$

$$\leq m\max\Big\{q\left(\tfrac{T(\bar{x}+tv)-T(\bar{x})-tT'(\bar{x},v)}{t}\right): q \in J\Big\}$$

$$< m\frac{\varepsilon}{m} = \varepsilon.$$

□

**Lemma 4.5**. *Let $T_1, T_2$: $D \rightarrow Y$ be single-valued mappings. Let $\bar{x} \in D$ and $v \in X\backslash\{\theta_X\}$. Suppose that both $T_1$ and $T_2$ are Gâteaux directionally differentiable at $\bar{x}$ along direction $v$. Then $T_1 + T_2$ is also Gâteaux directionally differentiable at $\bar{x}$ along direction $v$ and*

$$(T_1 + T_2)'(\bar{x},v) = T_1'(\bar{x},v) + T_2'(\bar{x},v). \tag{4.6}$$

*Proof.* Let $J \in \mathcal{F}_Y$ and $\varepsilon > 0$ be arbitrarily given. Then we consider $J \in \mathcal{F}_Y$ and $\frac{\varepsilon}{2} > 0$. By the assumption that, for $k = 1, 2$, $T_k$ is Gâteaux differentiable at $\bar{x}$ along direction $v$, for these given $J \in \mathcal{F}_Y$ and $\frac{\varepsilon}{2} > 0$, there is $\delta_k > 0$ such that, for any real number $t$, we have

$$0 < |t| < \delta_k \quad \Longrightarrow \quad \max\Big\{q\left(\tfrac{T_k(\bar{x}+tv)-T_k(\bar{x})-tT_k'(\bar{x},v)}{t}\right): q \in J\Big\} < \tfrac{\varepsilon}{2}, \text{ for } k = 1, 2. \tag{4.7}$$

Let $\delta = \min\{\delta_1, \delta_2\}$. By (4.7), for any real number $t$, if $0 < |t| < \delta$, we have

$$\max\left\{q\left(\frac{(T_1+T_2)(\bar{x}+tv)-(T_1+T_2)(\bar{x})-t(T_1'(\bar{x},v)+T_2'(\bar{x},v))}{t}\right): q \in J\right\}$$

$$\leq \max\left\{q\left(\frac{T_1(\bar{x}+tv)-T_1(\bar{x})-tT_1'(\bar{x},v))}{t}\right): q \in J\right\} + \max\left\{q\left(\frac{T_2(\bar{x}+tv)-T_2(\bar{x})-tT_2'(\bar{x},v)}{t}\right): q \in J\right\}$$

$$< \frac{\varepsilon}{2} + \frac{\varepsilon}{2} = \varepsilon.$$

(4.6) is proved. □

By Lemmas 4.23 and 4.24, we obtain the following linearity property of Gâteaux differentiation in Hausdorff topological vector spaces.

**Theorem 4.6**. *Let $T_1, T_2: D \to Y$ be single-valued mappings. Let $\bar{x} \in D$ and $v \in X\backslash\{\theta_X\}$ and let $a_1, a_2$ be scalars. Suppose that both $T_1$ and $T_2$ are Gâteaux directionally differentiable at $\bar{x}$ along direction $v$. Then $a_1T_1 + a_2T_2$ is also Gâteaux directionally differentiable at $\bar{x}$ along direction $v$ and*

$$(a_1T_1 + a_2T_2)'(\bar{x}, v) = a_1T_1'(\bar{x}, v) + a_2T_2'(\bar{x}, v). \tag{4.8}$$

*Furthermore, if both $T_1$ and $T_2$ are Gâteaux differentiable at $\bar{x}$. Then $a_1T_1 + a_2T_2$ is also Gâteaux differentiable at $\bar{x}$ and*

$$(a_1T_1 + a_2T_2)'(\bar{x})(v) = a_1T_1'(\bar{x})(v) + a_2T_2'(\bar{x})(v), \textit{for any } v \in X\backslash\{\theta_X\}.$$

*Proof*. The proof of this theorem is just by Lemmas 4.23 and 4.24. □

### 4.2. Fréchet Differentiability

In this section, in order to define Fréchet differentiation for single-valued mappings between general topological vector spaces, both the domain space and the range space are required to be equipped with $F$-seminorm structures. In this subsection, we use the notations defined in the previous sections 2 and 3. So, throughout this subsection and following subsections, if not otherwise stated, we always let $(X, \tau_X)$ and $(Y, \tau_Y)$ be Hausdorff topological vector spaces, in which the topologies $\tau_X$ and $\tau_Y$ are induced by the families $\mathbb{F}_X$ and $\mathbb{F}_Y$ of positive $F$-seminorms, respectively. Let $\mathcal{F}_X$ and $\mathcal{F}_Y$ respectively denote the collections of nonempty finite subsets of $\mathbb{F}_X$ and $\mathbb{F}_Y$.

In this subsection, we use the developed extended $\varepsilon$-$\delta$ language to define the Fréchet differentiation for single-valued mappings between general topological vector spaces. Before we pursue the definition of the Fréchet derivative, we consider some properties of $U_{I,\lambda}$, which will be used in the definitions and properties of the Fréchet derivative.

For each $p \in \mathbb{F}_X$, the kernel of $p$ is denoted and defined by $\text{Ker}(p) = \{x \in X: p(x) = 0\}$. Since $p(\theta_X) = 0$ and $p$ is positive, this implies that $\text{Ker}(p)$ is a nonempty proper subset of $X$. One more step further, for any $I \in \mathcal{F}_X$, we write

$$\text{Ker}(I) = \cap_{p\in I} \ker(p) = \cap_{p\in I}\{x \in X: p(x) = 0\}.$$

$\text{Ker}(I)$ is associated with the subset $I$ of the family $\mathbb{F}_X$ of positive $F$-seminorms on $X$. Since, for any $I \in \mathcal{F}_X$, $\theta_X \in \text{Ker}(I)$, it yields that $\text{Ker}(I)$ is a nonempty proper subset of $X$.

**Lemma 4.7.** *Let* $(X, \tau_X)$ *be a Hausdorff topological vector space, in which the topology* $\tau_X$ *is induced by a family* $\mathbb{F}_X$ *of positive F-seminorms. Then* $\text{Ker}(I)$ *has the following properties*

(i) $\text{Ker}(I)$ *is a proper subspace of X* (*possibly the trivial subspace* $\{\theta_X\}$);

(ii) $U_{I,\delta}\backslash\text{Ker}(I) \neq \emptyset$, *for any* $I \in \mathcal{F}_X$ *and any* $\delta > 0$.

*This implies that, for any* $I \in \mathcal{F}_X$ *and for any* $\delta > 0$, *we have*

$$\{u \in X: 0 < \max\{p(u): p \in I\} < \delta\} \neq \emptyset. \tag{4.9}$$

*Proof*. Part (i) can be straightforwardly proved by the properties of positive $F$-seminorms. We only prove (ii). By definition, for $u \in X$, we have

$$u \in U_{I,\delta}\backslash\text{Ker}(I) \quad \Leftrightarrow \quad 0 < \max\{p(u): p \in I\} < \delta. \tag{4.10}$$

Assume, by the way of contradiction, that there are $\delta > 0$ and $I \in \mathcal{F}_X$ such that $U_{I,\delta}\backslash\text{Ker}(I) = \emptyset$. This implies that $U_{I,\delta} \subseteq \text{Ker}(I)$, which implies that, for each $p \in I$, we have

$$p(x) = 0, \text{ for any } x \in U_{I,\delta} \subseteq \text{Ker}(I). \tag{4.11}$$

By part (e) in Lemma 2.4 , $U_{I,\delta}$ is an $\tau_X$-open balanced and absorbing neighborhood of $\theta_X$. Then, by the properties of $F$-seminorms and the definition of absorbing set, (4.11) implies that, for each $p \in I$, we have $p(x) = 0$, for any $x \in X$. This is a contradiction to the assumption that $p$ is a positive $F$-seminorm on $X$, which proves (4.9). □

Similarly to (3.3) in Definition 3.1, we define the Fréchet derivatives of mappings in topological vector spaces by using the extended $\varepsilon$-$\delta$ language. Recall that in this subsection, $D$ denotes an $\tau_X$-open nonempty subset of $X$ and $T$: $D \to Y$ denotes a single-valued mapping and let $\bar{x} \in D$.

**Definition 4.8.** (Fréchet differentiability) If there is a continuous linear mapping $\nabla T(\bar{x})$: $X \to Y$ such that, for any $J \in \mathcal{F}_Y$ and for any given $\varepsilon > 0$, there are $I \in \mathcal{F}_X$ and $\delta > 0$ such that, for $u \in X$,

(DZ) $\max\{p(u): p \in I\} = 0 \implies \max\{q\big(T(\bar{x} + u) - T(\bar{x}) - \nabla T(\bar{x})(u)\big): q \in J\} = 0;$

(DR) $$0 < \max\{p(u): p \in I\} < \delta \implies \max\left\{q\left(\frac{T(\bar{x}+u)-T(\bar{x})-\nabla T(\bar{x})(u)}{\max\{p(u):p\in I\}}\right): q \in J\right\} < \varepsilon, \tag{4.12}$$

then $T$ is said to be Fréchet differentiable at $\bar{x}$ and $\nabla T(\bar{x})$ is called the Fréchet derivative of $T$ at $\bar{x}$. Notice that as mentioned after Definition 4.1, by the continuity of the $F$-seminorms in $\mathbb{F}_X$, we always assume that $\delta$ is chosen to be small enough such that, for $u \in X$,

$$\max\{p(u): p \in I\} < \delta \implies \bar{x} + u \in D.$$

The above condition is to ensure that $T(\bar{x} + u)$ is defined. Hence, in following contents, we will always assume that condition $\bar{x} + u \in D$ is satisfied.

**Observations 4.9.** In particular, suppose that both $(X, \|\cdot\|_X)$ and $(Y, \|\cdot\|_Y)$ are normed vector spaces (in particular Banach spaces or Hilbert spaces) with both $\mathbb{F}_X = \{\|\cdot\|_X\}$ and $\mathbb{F}_Y = \{\|\cdot\|_Y\}$ singletons of norms on $X$ and $Y$, respectively. It is clear that the following conditions are satisfied

$$\{u \in X: \|u\|_X = 0\} = \{\theta_X\} \quad \text{and} \quad \{y \in Y: \|y\|_Y = 0\} = \{\theta_Y\}.$$

This implies that for the existence of $\nabla T(\bar{x})$: $X \to Y$, the condition (DZ) is automatically satisfied and the condition (DR) becomes the following ordinary limit

$$\lim_{u\to\theta_X} \frac{\|T(\bar{x}+u)-T(\bar{x})-\nabla T(\bar{x})(u)\|_Y}{\|u\|_X} = 0 \quad \Leftrightarrow \quad \lim_{u\to\theta_X} \frac{T(\bar{x}+u)-T(\bar{x})-\nabla T(\bar{x})(u)}{\|u\|_X} = \theta_Y.$$

This implies that the Fréchet derivative $\nabla T(\bar{x})$ of single-valued mapping defined in Definition 4.8 in general Hausdorff topological vector spaces is a natural extension of Fréchet derivatives in normed vector spaces, in particular in $F$-spaces, Banach spaces or Hilbert spaces.

Suppose that $T$ is Fréchet differentiable at a point $\bar{x} \in X$ with its Fréchet derivative $\nabla T(\bar{x})$. Does $\nabla T(\bar{x})$ uniquely exist? It is clear that this is a very important question. Here, we prove the uniqueness of Fréchet derivative under conditions (SC) and (SB) listed in Theorem 4.10 below.

**Theorem 4.10**. *Let $T$: $D \to Y$ be a single-valued mapping. Let $\bar{x} \in D$. Suppose $T$ is Fréchet differentiable at $\bar{x}$. If $X$ satisfies the following conditions*

- (SC) *$X$ is seminorm constructed*;
- (SB) *$\mathbb{F}_X$ satisfies the following boundedness condition*

$$\sup\{p(u): p \in \mathbb{F}_X\} < \infty, \text{ for every } u \in X, \tag{4.13}$$

*then the Fréchet derivative $\nabla T(\bar{x})$ of $T$ at $\bar{x}$ is unique.*

*Proof*. Suppose that there are two continuous and linear mappings $\nabla_1 T(\bar{x})$, $\nabla_2(\bar{x})$: $X \to Y$ both are continuous linear and both satisfy the conditions (DZ) and (DR) in Definition 4.8.

Let $w \in X\backslash\{\theta_X\}$ be arbitrarily given fixed. Let $J \in \mathcal{F}_Y$ be arbitrarily given. By the assumption (4.13) in this theorem, there is a positive integer $m$ (depending on $w$) such that

$$\sup\{p(w): p \in \mathbb{F}_X\} \leq m. \tag{4.14}$$

For the arbitrarily given (fixed) $J \in \mathcal{F}_Y$, let $\varepsilon > 0$ be arbitrarily given. We consider $J \in \mathcal{F}_Y$ and $\frac{\varepsilon}{2m} > 0$. By assumption that $\nabla_k T(\bar{x})$ is Fréchet derivative of $T$ at $\bar{x}$, for these $J \in \mathcal{F}_Y$ and $\frac{\varepsilon}{2m} > 0$, there are $\delta_k \in (0, 1)$ and $I_k \in \mathcal{F}_X$ such that, for $u \in X$ and for $k = 1, 2$, we have

$$\max\{p(u): p \in I_k\} = 0 \quad \Longrightarrow \quad \max\{q\big(T(\bar{x}+u) - T(\bar{x}) - \nabla_k T(\bar{x})(u)\big): q \in J\} = 0. \tag{4.15}$$

and

$$0 < \max\{p(u): p \in I_k\} < \delta_k \quad \Longrightarrow \quad \max\left\{q\left(\frac{T(\bar{x}+u)-T(\bar{x})-\nabla_k T(\bar{x})(u)}{\max\{p(u): p\in I_k\}}\right): q \in J\right\} < \frac{\varepsilon}{2m}. \tag{4.16}$$

By properties (ii) and (iii) of $F$-seminorms in Lemma 2.5, for the arbitrarily given $w \in X\backslash\{\theta_X\}$, there is a positive number $t_0 \in (0, 1)$ such that

$$\max\{p(t_0 w): p \in I_k\} < \delta_k, \text{ for } k = 1, 2.$$

By (4.15) and (4.16), we consider the following three cases.

Case 1. Suppose that $\max\{p(t_0 w): p \in I_k\} = 0$, for each $k \in \{1, 2\}$. By (4.15), this implies that

$$\max\{q\big(T(\bar{x}+t_0 w) - T(\bar{x}) - \nabla_k T(\bar{x})(t_0 w)\big): q \in J\} = 0, \text{ for } k = 1, 2. \tag{4.17}$$

Then, by (4.17), we estimate the following

$$\max\{q(t_0(\nabla_1 T(\bar{x})(w) - \nabla_2 T(\bar{x})(w))): q \in J\}$$

$$= \max\{q\big(\nabla_1 T(\bar{x})(t_0 w) - \nabla_2 T(\bar{x})(t_0 w)\big): q \in J\}$$

$$\leq \max\{q\big(T(\bar{x} + t_0 w) - T(\bar{x}) - \nabla_1 T(\bar{x})(t_0 w)\big): q \in J\} +$$

$$+ \max\{q\big(T(\bar{x} + t_0 w) - T(\bar{x}) - \nabla_2 T(\bar{x})(t_0 w)\big): q \in J\}$$

$$= 0. \tag{4.18}$$

Since $t_0 \in (0, 1)$, by Property (ix) of $F$-seminorms in Lemma 2.1, (4.18) implies that

$$\max\{q\big(\nabla_1 T(\bar{x})(w) - \nabla_2 T(\bar{x})(w)\big): q \in J\} = 0. \tag{4.19}$$

Case 2. Suppose that there is $j \in \{1, 2\}$ such that

$$\max\{p(t_0 w): p \in I_j\} = 0 \quad \text{and} \quad 0 < \max\{p(t_0 w): p \in I_{3-j}\} < \delta_{3-j}. \tag{4.20}$$

Then, in this case, we estimate the following

$$\max\{q\big(\nabla_1 T(\bar{x})(w) - \nabla_2 T(\bar{x})(w)\big): q \in J\}$$

$$= \max\{q\left(\frac{t_0(\nabla_1 T(\bar{x})(w) - \nabla_2 T(\bar{x})(w))}{t_0}\right): q \in J\}$$

$$= \max\{q\left(\frac{\nabla_1 T(\bar{x})(t_0 w) - \nabla_2 T(\bar{x})(t_0 w)}{t_0}\right): q \in J\}$$

$$\leq \max\{q\left(\frac{T(\bar{x}+t_0 w) - T(\bar{x}) - \nabla_j T(\bar{x})(t_0 w)}{t_0}\right): q \in J\} + \max\{q\left(\frac{T(\bar{x}+t_0 w) - T(\bar{x}) - \nabla_{3-j} T(\bar{x})(t_0 w)}{t_0}\right): q \in J\}. \tag{4.21}$$

By (4.15), the equation in (4.20) implies that

$$\max\{q\left(T(\bar{x} + t_0 w) - T(\bar{x}) - \nabla_j T(\bar{x})(t_0 w)\right): q \in J\} = 0. \tag{4.22}$$

By Property (ix) of $F$-seminorms in Lemma 2.1, (4.22) implies that

$$\max\{q\left(\frac{T(\bar{x}+t_0 w) - T(\bar{x}) - \nabla_j T(\bar{x})(t_0 w)}{t_0}\right): q \in J\} = 0. \tag{4.23}$$

By the inequalities in (4.20), by Property (vii) of $F$-seminorms in Lemma 2.1, by the condition (SC) in this theorem and by (4.14) and (4.16), we have

$$\max\{q\left(\frac{T(\bar{x}+t_0 w) - T(\bar{x}) - \nabla_{3-j} T(\bar{x})(t_0 w)}{t_0}\right): q \in J\}$$

$$= \max\{q\left(\frac{\max\{p(w): p \in I_{3-j}\}\,(T(\bar{x}+t_0 w) - T(\bar{x}) - \nabla_{3-j} T(\bar{x})(t_0 w))}{t_0 \max\{p(w): p \in I_{3-j}\}}\right): q \in J\}$$

$$= \max\{q\left(\frac{\max\{p(w): p \in I_{3-j}\}\,(T(\bar{x}+t_0 w) - T(\bar{x}) - \nabla_{3-j} T(\bar{x})(t_0 w))}{\max\{p(t_0 w): p \in I_{3-j}\}}\right): q \in J\}$$

$$\leq \max\{q\left(\frac{m\,(T(\bar{x}+t_0w)-T(\bar{x})-\nabla_{3-j}T(\bar{x})(t_0w))}{\max\{p(t_0w):p\in I_{3-j}\}}\right):q\in J\}$$

$$\leq m\max\{q\left(\frac{(T(\bar{x}+t_0w)-T(\bar{x})-\nabla_{3-j}T(\bar{x})(t_0w))}{\max\{p(t_0w):p\in I_{3-j}\}}\right):q\in J\}$$

$$< m\frac{\varepsilon}{2m}=\frac{\varepsilon}{2}. \tag{4.24}$$

Substituting (4.23) and (4.24) into (4.21), in case 2, we obtain that, for $j\in\{1,2\}$, if

$$\max\{p(t_0w):p\in I_j\}=0 \quad \text{and} \quad 0<\max\{p(t_0w):p\in I_{3-j}\}<\delta_{3-j},$$

then
$$\max\{q\big(\nabla_1T(\bar{x})(w)-\nabla_2T(\bar{x})(w)\big):q\in J\}<\frac{\varepsilon}{2}. \tag{4.25}$$

Case 3. Suppose that $0<\max\{p(t_0w):p\in I_k\}<\delta_k$, for each $k\in\{1,2\}$. Similar to the proof of case 2, we estimate

$$\max\{q\big(\nabla_1T(\bar{x})(w)-\nabla_2T(\bar{x})(w)\big):q\in J\}$$

$$\leq \max\{q\left(\frac{T(\bar{x}+t_0w)-T(\bar{x})-\nabla_1T(\bar{x})(t_0w)}{t_0}\right):q\in J\}+\max\{q\left(\frac{T(\bar{x}+t_0w)-T(\bar{x})-\nabla_2T(\bar{x})(t_0w)}{t_0}\right):q\in J\}.$$

$$\leq \max\{q\left(\frac{m\,(T(\bar{x}+t_0w)-T(\bar{x})-\nabla_1T(\bar{x})(t_0w))}{\max\{p(t_0w):p\in I_1\}}\right):q\in J\}+\max\{q\left(\frac{m\,(T(\bar{x}+t_0w)-T(\bar{x})-\nabla_2T(\bar{x})(t_0w))}{\max\{p(t_0w):p\in I_2\}}\right):q\in J\}$$

$$\leq m\max\{q\left(\frac{m\,(T(\bar{x}+t_0w)-T(\bar{x})-\nabla_1T(\bar{x})(t_0w))}{\max\{p(t_0w):p\in I_1\}}\right):q\in J\}+m\max\{q\left(\frac{m\,(T(\bar{x}+t_0w)-T(\bar{x})-\nabla_1T(\bar{x})(t_0w))}{\max\{p(t_0w):p\in I_1\}}\right):q\in J\}$$

$$< m\frac{\varepsilon}{2m}+m\frac{\varepsilon}{2m}=\varepsilon. \tag{4.26}$$

(4.26) implies that, in case 3, we obtain that if $0<\max\{p(t_0w):p\in I_k\}<\delta_k$, for each $k\in\{1,2\}$, then

$$\max\{q\big(\nabla_1T(\bar{x})(w)-\nabla_2T(\bar{x})(w)\big):q\in J\}<\varepsilon. \tag{4.27}$$

Combining (4.19) for case 1, (4.25) for case 2, and (4.27) for case 3, for the given $J\in\mathcal{F}_Y$ we have that,

$$\max\{q\big(\nabla_1T(\bar{x})(w)-\nabla_2T(\bar{x})(w)\big):q\in J\}<\varepsilon, \text{ for any } \varepsilon>0.$$

This implies that

$$\max\{q\big(\nabla_1T(\bar{x})(w)-\nabla_2T(\bar{x})(w)\big):q\in J\}=0. \tag{4.28}$$

Since $J\in\mathcal{F}_Y$ is arbitrarily given, (4.28) implies that

$$q\big(\nabla_1T(\bar{x})(w)-\nabla_2T(\bar{x})(w)\big)=0, \text{ for any } q\in\mathbb{F}_Y. \tag{4.29}$$

Notice that $w\in X\backslash\{\theta_X\}$ is arbitrarily given. By the condition that $Y$ is a Hausdorff topological vector spaces, (4.29) implies that

$$\nabla_1T(\bar{x})(w)-\nabla_2T(\bar{x})(w)=\theta_Y, \text{ for any } w\in X\backslash\{\theta_X\}.$$

This proves that $\nabla_1T(\bar{x})=\nabla_2T(\bar{x})$. □

**Corollary 4.11**. *Let* $(X, \|\cdot\|_X)$ *be a normed vector space equipped with a norm* $\|\cdot\|_X$. *Let Y be a topological vector space. Let* $T: D \to Y$ *be a single-valued mapping. Let* $\bar{x} \in D$. *Suppose T is Fréchet differentiable at* $\bar{x}$. *Then the Fréchet derivative* $\nabla T(\bar{x})$ *of T at* $\bar{x}$ *is unique.*

*Proof*. This normed vector space $(X, \|\cdot\|_X)$ is a special topological vector space, in which the topology is induced by a singleton family $\mathbb{F}_X = \{\|\cdot\|_X\}$. It is clear that the conditions (SC) and (SB) in Theorem 10 are satisfied. □

Next, we prove some properties of Fréchet differentiability in Hausdorff topological vector spaces.

**Lemma 4.12**. *Let* $T: D \to Y$ *be a single-valued mapping. Let* $\bar{x} \in D$ *and let* $a$ *be a scalar. Suppose that T is Fréchet differentiable at* $\bar{x}$. *Then* $aT: X \to Y$ *is also Fréchet differentiable at* $\bar{x}$ *and*

$$\nabla(aT)(\bar{x}) = a\nabla T(\bar{x}). \tag{4.30}$$

*Proof.* The proof of this lemma is similar to the proof of Lemma 3.2. It is clear that (4.30) holds for $a = 0$. So, we only prove (4.30) for $a \neq 0$.

For the given real number $a \neq 0$, there is a positive integer $n$ such that $|a| \leq n$. Let $J \in \mathcal{F}_Y$ and $\varepsilon > 0$ be arbitrarily given. We consider $J \in \mathcal{F}_Y$ and $\frac{\varepsilon}{n} > 0$. By the assumption that $T$ is Fréchet differentiable at $\bar{x}$, there are $I \in \mathcal{F}_X$ and $\delta > 0$ such that,

$$\max\{p(u): p \in I\} = 0 \quad \Longrightarrow \quad \max\{q\big(T(\bar{x}+u) - T(\bar{x}) - \nabla T(\bar{x})(u)\big): q \in J\} = 0. \tag{4.31}$$

and

$$0 < \max\{p(u): p \in I\} < \delta \quad \Longrightarrow \quad \max\left\{q\left(\frac{T(\bar{x}+u)-T(\bar{x})-\nabla T(\bar{x})(u)}{\max\{p(u):p\in I\}}\right): q \in J\right\} < \frac{\varepsilon}{n}. \tag{4.32}$$

By $|a| \leq n$ and by properties (vii) and (iii) of $F$-seminorms, if $\max\{p(u): p \in I\} = 0$, then (4.31) implies

$$\begin{aligned}
&\max\{q\big(aT(\bar{x}+u) - aT(\bar{x}) - a\nabla T(\bar{x})(u)\big): q \in J\}\\
&= \max\{q(a(T(\bar{x}+u) - T(\bar{x}) - \nabla T(\bar{x})(u))): q \in J\}\\
&\leq \max\{q(n(T(\bar{x}+u) - T(\bar{x}) - \nabla T(\bar{x})(u))): q \in J\}\\
&\leq n\max\{q\big(T(\bar{x}+u) - T(\bar{x}) - \nabla T(\bar{x})(u)\big): q \in J\}\\
&= 0.
\end{aligned}$$

This implies that

$$\max\{p(u): p \in I\} = 0 \quad \Longrightarrow \quad \max\{q\big(aT(\bar{x}+u) - aT(\bar{x}) - a\nabla T(\bar{x})(u)\big): q \in J\} = 0.$$

Hence, the condition (DZ) in Definition 4.8 is satisfied. Next, we show that the condition (DR) in Definition 4.8 is also satisfied. To this end, by properties (vii) and (iii) of $F$-seminorms and by $|a| \leq n$, we have

$$q(au) \leq q(nu) \leq nq(u), \text{ for every } u \in X. \tag{4.33}$$

For $u \in X$, by (4.32) and (4.33), we have

$$0 < \max\{p(u): p \in I\} < \delta$$

$$\Longrightarrow \quad \max\left\{q\left(\frac{(aT)(\bar{x}+u)-(aT)(\bar{x})-(a\nabla T(\bar{x}))(u)}{\max\{p(u):p\in I\}}\right): q \in J\right\}$$

$$= \max\left\{q\left(\frac{a(T(\bar{x}+u)-T(\bar{x})-\nabla T(\bar{x})(u))}{\max\{p(u):p\in I\}}\right): q\in J\right\}$$

$$\leq \max\left\{q\left(\frac{n(T(\bar{x}+u)-T(\bar{x})-\nabla T(\bar{x})(u))}{\max\{p(u):p\in I\}}\right): q\in J\right\}$$

$$\leq n\max\left\{q\left(\frac{T(\bar{x}+u)-T(\bar{x})-\nabla T(\bar{x})(u)}{\max\{p(u):p\in I\}}\right): q\in J\right\}$$

$$< n\frac{\varepsilon}{n} = \varepsilon. \tag{4.34}$$

Hence, the condition (DR) in Definition 4.8 is also satisfied. Then, this lemma is proved. □

**Lemma 4.13**. *Let $T_1, T_2: D \longrightarrow Y$ be single-valued mappings. Let $\bar{x} \in D$. Suppose that both $T_1$ and $T_2$ are Fréchet differentiable at $\bar{x}$. Then $T_1 + T_2: D \longrightarrow Y$ is also Fréchet differentiable at $\bar{x}$ and*

$$\nabla(T_1+T_2)(\bar{x}) = \nabla T_1(\bar{x}) + \nabla T_2(\bar{x}).$$

*Proof.* Let $J \in \mathcal{F}_Y$ and $\varepsilon > 0$ be arbitrarily given. We consider $J \in \mathcal{F}_Y$ and $\frac{\varepsilon}{2} > 0$. By the assumption that, for $k$ = 1, 2, $T_k$ is Fréchet differentiable at $\bar{x}$ ($\nabla T_k(\bar{x})$ exists), there are $\delta_k > 0$ and $I_k \in \mathcal{F}_X$ such that, for $u \in X$, we have

(i) $\max\{p(u): p \in I_k\} = 0 \quad \Longrightarrow \max\{q\big(T_k(\bar{x}+u) - T_k(\bar{x}) - \nabla T_k(\bar{x})(u)\big): q \in J\} = 0.$

(ii) $0 < \max\{p(u): p \in I_k\} < \delta_k \quad \Longrightarrow \quad \max\left\{q\left(\frac{T_k(\bar{x}+u)-T_k(\bar{x})-\nabla T_k(\bar{x})(u)}{\max\{p(u):p\in I_k\}}\right): q\in J\right\} < \frac{\varepsilon}{2}.$ (4.35)

Let $\delta = \min\{\delta_1, \delta_2\}$ and let $I = I_1 \cup I_2$. We have $\delta > 0$ and $I \in \mathcal{F}_X$, by which and by (4.35), we prove that the operator $T_1 + T_2$ satisfies both conditions ((DZ)) and (DR) in Definition 4.8.

We first check condition (DZ) in Definition 4.8 for $T_1 + T_2$. For $u \in X$, suppose that $\max\{p(u): p \in I\} = 0$. By $I = I_1 \cup I_2$, we have

$$\max\{p(u): p \in I\} = 0 \quad \Longrightarrow \quad \max\{p(u): p \in I_k\} = 0, \text{ for } k = 1, 2.$$

Then, by part (i) in (4.35), we obtain that

$$\max\{p(u): p \in I\} = 0$$

$$\Longrightarrow \quad \max\{q[(T_1+T_2)(\bar{x}+u) - (T_1+T_2)(\bar{x}) - (\nabla T_1(\bar{x}) + \nabla T_2(\bar{x}))(u))]: q \in J\}$$

$$= \max\{q[\big(T_1(\bar{x}+u) - T_1(\bar{x}) - \nabla T_1(\bar{x})(u)\big) + \big(T_2(\bar{x}+u) - T_2(\bar{x}) - \nabla T_2(\bar{x})(u)\big)]: q \in J\}$$

$$\leq \max\{q\big(T_1(\bar{x}+u) - T_1(\bar{x}) - \nabla T_1(\bar{x})(u)\big) + q(T_2(\bar{x}+u) - T_2(\bar{x}) - \nabla T_2(\bar{x})(u): q \in J\}$$

$$\leq \max\{q(T_1(\bar{x}+u) - T_1(\bar{x}) - \nabla T_1(\bar{x})(u)): q \in J\} + \max\{q(T_2(\bar{x}+u) - T_2(\bar{x}) - \nabla T_2(\bar{x})(u): q \in J\}$$

$$= 0.$$

This yields that

$$\max\{p(u): p \in I\} = 0$$

$$\Longrightarrow \max\{q[(T_1+T_2)(\bar{x}+u)-(T_1+T_2)(\bar{x})-(\nabla T_1(\bar{x})+\nabla T_2(\bar{x}))(u))]: q \in J\} = 0.$$

This proves that condition (DZ) in Definition 4.8 is satisfied for $T_1+T_2$ with respect to this $\tau_X$-open neighborhood $U_{I,\delta}$ of $\theta_X$ in $X$. Next, we show that condition (DR) in Definition 4.8 is satisfied for $T_1+T_2$ with respect to $\delta > 0$ and $I \in \mathcal{F}_X$. For $u \in X$, suppose that

$$0 < \max\{p(u): p \in I\} < \delta.$$

By the definitions that $\delta = \min\{\delta_1, \delta_2\}$ and $I = I_1 \cup I_2$, there are two cases for consideration.

(a) There is $k \in \{1,2\}$ such that

$$\max\{p(u): p \in I_k\} = 0 \quad \text{and} \quad 0 < \max\{p(u): p \in I_{3-k}\} < \delta_{3-k}.$$

For such a $k$, by part (i) in (4.35), this implies

$$\max\{q\big(T_k(\bar{x}+u)-T_k(\bar{x})-\nabla T_k(\bar{x})(u)\big): q \in J\} = 0.$$

By property (ix) of $F$-seminorms in Lemma 2.3, this implies that

$$\max\left\{q\left(\frac{T_k(\bar{x}+u)-T_k(\bar{x})-\nabla T_k(\bar{x})(u)}{\max\{p(u):p\in I\}}\right): q \in J\right\} = 0.$$

Then, by part (ii) in (4.35), we have

$$\max\left\{q\left(\frac{(T_1+T_2)(\bar{x}+u)-(T_1+T_2)(\bar{x})-(\nabla T_1(\bar{x})+\nabla T_2(\bar{x}))(u)}{\max\{p(u):p\in I\}}\right): q \in J\right\}$$

$$= \max\left\{q\left(\frac{\big(T_1(\bar{x}+u)-T_1(\bar{x})-\nabla T_1(\bar{x})(u)\big)+\big(T_2(\bar{x}+u)-T_2(\bar{x})-\nabla T_2(\bar{x})(u)\big)}{\max\{p(u):p\in I\}}\right): q \in J\right\}$$

$$\leq \max\left\{q\left(\frac{T_k(\bar{x}+u)-T_k(\bar{x})-\nabla T_k(\bar{x})(u)}{\max\{p(u):p\in I\}}\right): q \in J\right\} + \max\left\{q\left(\frac{T_{3-k}(\bar{x}+u)-T_{3-k}(\bar{x})-\nabla T_{3-k}(\bar{x})(u)}{\max\{p(u):p\in I\}}\right): q \in J\right\}$$

$$\leq \max\left\{q\left(\frac{T_k(\bar{x}+u)-T_k(\bar{x})-\nabla T_k(\bar{x})(u)}{\max\{p(u):p\in I\}}\right): q \in J\right\}$$

$$+ \max\left\{q\left(\frac{\max\{p(u):p\in I_{3-k}\}}{\max\{p(u):p\in I\}}\frac{T_{3-k}(\bar{x}+u)-T_{3-k}(\bar{x})-\nabla T_{3-k}(\bar{x})(u)}{\max\{p(u):p\in I_{3-k}\}}\right): q \in J\right\}$$

$$\leq \max\left\{q\left(\frac{T_k(\bar{x}+u)-T_k(\bar{x})-\nabla T_k(\bar{x})(u)}{\max\{p(u):p\in I\}}\right): q \in J\right\} + \max\left\{q\left(\frac{T_{3-k}(\bar{x}+u)-T_{3-k}(\bar{x})-\nabla T_{3-k}(\bar{x})(u)}{\max\{p(u):p\in I_{3-k}\}}\right): q \in J\right\}$$

$$\leq \max\left\{q\left(\frac{T_k(\bar{x}+u)-T_k(\bar{x})-\nabla T_k(\bar{x})(u)}{\max\{p(u):p\in I\}}\right): q \in J\right\} + \frac{\varepsilon}{2}$$

$$= 0 + \frac{\varepsilon}{2}.$$

$$= \varepsilon.$$

The equality comes from part (ix) in Lemma 2.1 and the assumption in case (a).

(b) Suppose that, for every $k \in \{1,2\}$, we have $0 < \max\{p(u): p \in I_k\} < \delta_k$. By $\delta = \min\{\delta_1, \delta_2\}$ and $I = I_1 \cup I_2$, by (ii) in (4.35), for every $k = 1, 2$, we have that

$$0 < \max\{p(u): p \in I\} < \delta$$

$$\Longrightarrow \quad 0 < \max\{p(u): p \in I_k\} \leq \max\{p(u): p \in I\} < \delta \leq \delta_k$$

$$\Longrightarrow \quad \max\left\{q\left(\frac{T_k(\bar{x}+u)-T_k(\bar{x})-\nabla T_k(\bar{x})(u)}{\max\{p(u):p\in I\}}\right): q \in J\right\} \leq \max\left\{q\left(\frac{T_k(\bar{x}+u)-T_k(\bar{x})-\nabla T_k(\bar{x})(u)}{\max\{p(u):p\in I_k\}}\right): q \in J\right\} < \frac{\varepsilon}{2}.$$

This yields that

$$0 < \max\{p(u): p \in I\} < \delta$$

$$\Longrightarrow \quad \max\left\{q\left(\frac{(T_1+T_2)(\bar{x}+u)-(T_1+T_2)(\bar{x})-(\nabla T_1(\bar{x})+\nabla T_2(\bar{x}))(u)}{\max\{p(u):p\in I\}}\right): q \in J\right\}$$

$$= \max\left\{q\left(\frac{\left(T_1(\bar{x}+u)-T_1(\bar{x})-\nabla T_1(\bar{x})(u)\right)+\left(T_2(\bar{x}+u)-T_2(\bar{x})-\nabla T_2(\bar{x})(u)\right)}{\max\{p(u):p\in I\}}\right): q \in J\right\}$$

$$\leq \max\left\{q\left(\frac{T_1(\bar{x}+u)-T_1(\bar{x})-\nabla T_1(\bar{x})(u)}{\max\{p(u):p\in I\}}\right): q \in J\right\} + \max\left\{q\left(\frac{T_2(\bar{x}+u)-T_2(\bar{x})-\nabla T_2(\bar{x})(u)}{\max\{p(u):p\in I\}}\right): q \in J\right\}$$

$$\leq \max\left\{q\left(\frac{T_1(\bar{x}+u)-T_1(\bar{x})-\nabla T_1(\bar{x})(u)}{\max\{p(u):p\in I_1\}}\right): q \in J\right\} + \max\left\{q\left(\frac{T_2(\bar{x}+u)-T_2(\bar{x})-\nabla T_2(\bar{x})(u)}{\max\{p(u):p\in I_2\}}\right): q \in J\right\}$$

$$< \frac{\varepsilon}{2} + \frac{\varepsilon}{2}$$

$$= \varepsilon.$$

This proves that the mapping $T_1 + T_2$ satisfies condition (DR) in Definition 4.8. This completes the proof of this lemma. □

By Lemmas 4.6 and 4.7, we obtain the following linearity property of Fréchet differentiation in Hausdorff topological vector spaces.

**Theorem 4.14**. *Let $T_1, T_2: D \to Y$ be single-valued mappings. Let $\bar{x} \in D$ and let $a_1$, $a_2$ be scalars. Suppose that both $T_1$ and $T_2$ are Fréchet differentiable at $\bar{x}$. Then $a_1T_1 + a_2T_2: D \to Y$ is also Fréchet differentiable at $\bar{x}$ and*

$$\nabla(a_1T_1 + a_2T_2)(\bar{x}) = a_1\nabla T_1(\bar{x}) + a_2\nabla T_2(\bar{x}). \tag{4.36}$$

*Proof*. The proof of this theorem is just by Lemmas 4.6 and 4.7. □

In the following proposition, we consider the Fréchet differentiability of continuous linear mappings on Hausdorff topological vector spaces, which are in particularly important in advanced analysis.

**Proposition 4.15**. *Let $T: D \to Y$ be a continuous linear mapping. Then, $T$ is Fréchet differentiable on $D$ and, for any $\bar{x} \in D$, the Fréchet derivative $\nabla T(\bar{x})$ of $T$ at point $\bar{x}$ satisfies that*

$$\nabla T(\bar{x})(u) = T(u), \text{ for every } u \in X.$$

*Proof*. For any given $\bar{x} \in X$, we have

$$T(\bar{x} + u) - T(\bar{x}) - T(u) = \theta_Y, \text{ for any } u \in X.$$

Then, this implies that conditions (DZ) and (DR) in Definition 4.8 are clearly satisfied, which proves this

proposition. □

Next, we study the connection between continuity and Fréchet differentiability of mappings between general topological vector spaces.

**Theorem 4.16**. *Let $T: D \to Y$ be a single-valued mapping. Let $\bar{x} \in D$. If $T$ is Fréchet differentiable at $\bar{x}$, then $T$ is continuous at $\bar{x}$.*

*Proof*. Let $J \in \mathcal{F}_Y$ and $\varepsilon > 0$ be arbitrarily given. We consider $J \in \mathcal{F}_Y$ and $\frac{\varepsilon}{2} > 0$. By the assumption in this proposition that $\nabla T(\bar{x})$ exists, there are $\delta_1 > 0$ satisfying $\delta_1 < 1$ and $I_1 \in \mathcal{F}_X$ such that, for $u \in X$, $u$ satisfies the conditions (DZ) and (DR) in Definition 4.8.

(i) $\quad \max\{p(u): p \in I_1\} = 0 \implies \max\{q\big(T(\bar{x}+u) - T(\bar{x}) - \nabla T(\bar{x})(u)\big): q \in J\} = 0;$

(ii) $\quad 0 < \max\{p(u): p \in I_1\} < \delta_1 \implies \max\left\{\frac{q[T(\bar{x}+u)-T(\bar{x})-\nabla T(\bar{x})(u)]}{\max\{p(u): p\in I_1\}}: q \in J\right\} < \frac{\varepsilon}{2}.$

Since $\nabla T(\bar{x}): X \to Y$ is a continuous and linear mapping, for this $\tau_Y$-open neighborhood $V_{J,\frac{\varepsilon}{2}}$ of $\theta_Y$ in $Y$, there are $\delta_2 > 0$ with $\delta_2 < 1$ and $I_2 \in \mathcal{F}_X$ such that, for $u \in U_{I_2,\delta_2}$, we have

$$\max\{p(u): p \in I_2\} < \delta_2 \implies \max\{q(\nabla T(\bar{x})(u)): q \in J\} < \frac{\varepsilon}{2}. \tag{4.37}$$

Let $\delta = \min\{\delta_1, \delta_2\}$ and let $I = I_1 \cup I_2$. Then, we have an $\tau_X$-open neighborhood $U_{I,\delta}$ of $\theta_X$ in $X$. For any $u \in U_{I,\delta}$, we consider two cases with the above two conditions (i) and (ii) for $T$.

Case 1. $\max\{p(u): p \in I\} = 0$. Since $I = I_1 \cup I_2$, this implies

$$\max\{p(u): p \in I_1\} = 0 \quad \text{and} \quad \max\{p(u): p \in I_2\} = 0 < \delta_2.$$

Hence, by the above property (i) and (4.37) of $T$, we have

$$\begin{aligned}
&\max\{p(u): p \in I\} = 0 \\
\implies\ &\max\{p(u): p \in I_1\} = 0 \\
\implies\ &\max\{q\big(T(\bar{x}+u) - T(\bar{x})\big): q \in J\} \\
&\leq \max\{q\big(T(\bar{x}+u) - T(\bar{x}) - \nabla T(\bar{x})(u)\big): q \in J\} + \max\{q\big(\nabla T(\bar{x})(u)\big): q \in J\} \\
&< 0 + \frac{\varepsilon}{2}.
\end{aligned}$$

Case 2. $0 < \max\{p(u): p \in I\} < \delta$. This implies that

$$\max\{p(u): p \in I_1\} \leq \max\{p(u): p \in I\} < \delta \leq \delta_1$$

and

$$\max\{p(u): p \in I_2\} \leq \max\{p(u): p \in I\} < \delta \leq \delta_2.$$

Subcase 2.1. $\max\{p(u): p \in I_1\} = 0$. Then, by property (i) and (4.37), we have

$$\max\{q\big(T(\bar{x}+u) - T(\bar{x})\big): q \in J\}$$

$$\leq \max\{q(T(\bar{x}+u)-T(\bar{x})-\nabla T(\bar{x})(u)): q \in J\} + \max\{q(\nabla T(\bar{x})(u)): q \in J\}$$

$$< 0 + \frac{\varepsilon}{2}.$$

Subcase 2.1. $0 < \max\{p(u): p \in I_1\} < \delta \leq \delta_1$. By the above property (ii) of $U_{I_1,\delta_1}$, this implies

$$\max\left\{q\left(\frac{T(\bar{x}+u)-T(\bar{x})-\nabla T(\bar{x})(u)}{\max\{p(u):p\in I_1\}}\right): q \in J\right\} < \frac{\varepsilon}{2}.$$

Notice that $\delta_1 < 1$, by property (viii) of $F$-seminorms in Lemma 2.2, we have

$$\max\{q(T(\bar{x}+u)-T(\bar{x})): q \in J\}$$

$$\leq \max\{q(T(\bar{x}+u)-T(\bar{x})-\nabla T(\bar{x})(u)): q \in J\} + \max\{q(\nabla T(\bar{x})(u)): q \in J\}$$

$$= \max\left\{q\left(\max\{p(u): p \in I_1\}\frac{T(\bar{x}+u)-T(\bar{x})-\nabla T(\bar{x})(u)}{\max\{p(u):p\in I_1\}}\right): q \in J\right\} + \max\{q(\nabla T(\bar{x})(u)): q \in J\}$$

$$\leq \max\left\{q\left(\delta_1\frac{T(\bar{x}+u)-T(\bar{x})-\nabla T(\bar{x})(u)}{\max\{p(u):p\in I_1\}}\right): q \in J\right\} + \max\{q(\nabla T(\bar{x})(u)): q \in J\}$$

$$\leq \max\left\{q\left(\frac{T(\bar{x}+u)-T(\bar{x})-\nabla T(\bar{x})(u)}{\max\{p(u):p\in I_1\}}\right): q \in J\right\} + \max\{q(\nabla T(\bar{x})(u)): q \in J\}$$

$$< \frac{\varepsilon}{2} + \frac{\varepsilon}{2}.$$

Combining Cases 1 and 2, we obtain that

$$\max\{p(u): p \in I\} < \delta \implies \max\{q(T(\bar{x}+u)-T(\bar{x})): q \in J\} < \varepsilon. \qquad \square$$

**Question 4.17**. Recall that the Fréchet derivative of mappings in Banach spaces satisfies the chain-rule (see Theorem 2.1 in [1]). Under the Definition 4.8, does the Fréchet differentiability of mappings in general topological vector spaces also satisfy a certain type of chain-rule? We present this question in detail in the conclusion section.

### 4.3. Connection Between Gâteaux Differentiability and Fréchet Differentiability

In this subsection, we consider the connection between Gâteaux differentiability and Fréchet differentiability of single-valued mappings in Hausdorff topological vector spaces. Here, in general topological vector spaces, we cannot prove (1.3) that Fréchet differentiability is stronger than Gâteaux differentiability, however, under some conditions given in the following theorem, we prove connection between them, which is similar to (1.3).

**Theorem 4.18**. *Let $T: D \to Y$ be a single-valued mapping. Let $\bar{x} \in X$. Suppose $T$ is Fréchet differentiable at $\bar{x}$. Let $v \in X \backslash \{\theta_X\}$. If the following conditions are satisfied,*

(SC) $(X, \tau_X)$ *is seminorm constructed*;

(sb) $\max\{p(v): p \in \mathbb{F}_X\} < \infty$, (4.38)

*then $T$ is Gâteaux directionally differentiable at $\bar{x}$ along the direction $v$ and*

$$T'(\bar{x}, v) = \nabla T(\bar{x})(v).$$

*Proof*. Let $v \in X\backslash\{\theta_X\}$. Suppose that $\max\{p(v): p \in \mathbb{F}_X\} < \infty$. Since $(X, \tau_X)$ is Hausdorff, it implies that $0 < \max\{p(v): p \in \mathbb{F}_X\}$. Then, we can assume that there is a positive integer $m$ such that

$$0 < \max\{p(v): p \in \mathbb{F}_X\} < m.$$

Let $J \in \mathcal{F}_Y$ and $\varepsilon > 0$ be arbitrarily given. We consider $J \in \mathcal{F}_Y$ and $\frac{\varepsilon}{m} > 0$. By the assumption that $T$ is Fréchet differentiable at $\bar{x}$, there are $I \in \mathcal{F}_X$ and $\delta_1 \in (0, 1)$ such that, for any $u \in X$, we have

(i) $\max\{p(u): p \in I\} = 0 \implies \max\{q(T(\bar{x}+u) - T(\bar{x}) - \nabla T(\bar{x})(u)): q \in J\} = 0;$

(ii) $0 < \max\{p(u): p \in I\} < \delta_1 \implies \max\left\{q\left(\frac{T(\bar{x}+u)-T(\bar{x})-\nabla T(\bar{x})(u)}{\max\{p(u):p\in I\}}\right): q \in J\right\} < \frac{\varepsilon}{m}.$

For the given $v \in X\backslash\{\theta_X\}$, we consider the following two cases regarding to the above conditions (i, ii).

Case 1. If $\max\{p(v): p \in I\} = 0$, then by property of seminorms, we have

$$\max\{p(tv): p \in I\} = 0, \text{ for any real } t.$$

By the above condition (i), this implies that

$$\max\{q(T(\bar{x} + tv) - T(\bar{x}) - \nabla T(\bar{x})(tv)): q \in J\} = 0, \text{ for any real } t. \tag{4.39}$$

Then by the property (ix) of $F$-seminorms in Lemma 2.2, (4.39) implies that,

$$\max\left\{q\left(\frac{T(\bar{x}+tv)-T(\bar{x})-\nabla T(\bar{x})(tv)}{t}\right): q \in J\right\} = 0 < \varepsilon, \text{ for any real } t \text{ with } t \neq 0.$$

Hence, we obtain that

$$\max\{p(v): p \in I\} = 0 \implies \max\left\{q\left(\frac{T(\bar{x}+tv)-T(\bar{x})-\nabla T(\bar{x})(tv)}{t}\right): q \in J\right\} < \varepsilon, \text{ for any } t \text{ with } t \neq 0. \tag{4.40}$$

Case 2. Suppose that $\max\{p(v): p \in I\} > 0$. By Property (iii) in Lemma 2.5 and (ix) in Lemma 2.2, there is $\delta \in (0, \delta_1)$ such that, for any real number $t$,

$$0 < |t| < \delta \implies 0 < \max\{p(tv): p \in I\} < \delta_1. \tag{4.41}$$

Since $(X, \tau_X)$ is seminorm constructed, by (4.41) and the condition (ii) in this theorem, we have

$$\begin{aligned}
&0 < |t| < \delta \\
\implies\ &\max\left\{q\left(\frac{T(\bar{x}+tv)-T(\bar{x})-t\nabla T(\bar{x})(v)}{t}\right): q \in J\right\} \\
&= \max\left\{q\left(\frac{|t|\max\{p(v):p\in I\}}{t}\frac{T(\bar{x}+tv)-T(\bar{x})-\nabla T(\bar{x})(tv)}{|t|\max\{p(v):p\in I\}}\right): q \in J\right\} \\
&= \max\left\{q\left(\frac{|t|\max\{p(v):p\in I\}}{t}\frac{T(\bar{x}+tv)-T(\bar{x})-\nabla T(\bar{x})(tv)}{\max\{p(tv):p\in I\}}\right): q \in J\right\} \\
&\leq \max\left\{q\left(\frac{\max\{p(v):p\in I\}}{1}\frac{T(\bar{x}+tv)-T(\bar{x})-\nabla T(\bar{x})(tv)}{\max\{p(tv):p\in I\}}\right): q \in J\right\}
\end{aligned}$$

$$\leq \max\left\{q\left(m\frac{T(\bar{x}+tv)-T(\bar{x})-\nabla T(\bar{x})(tv)}{\max\{p(tv):p\in I\}}\right): q\in J\right\}$$

$$\leq m\max\left\{q\left(\frac{T(\bar{x}+tv)-T(\bar{x})-\nabla T(\bar{x})(tv)}{\max\{p(tv):p\in I\}}\right): q\in J\right\}$$

$$< m\frac{\varepsilon}{m} = \varepsilon.$$

This implies that

$$0<|t|<\delta \quad \Longrightarrow \quad \max\left\{q\left(\frac{T(\bar{x}+tv)-T(\bar{x})-t\nabla T(\bar{x})(v)}{t}\right): q\in J\right\} < \varepsilon. \tag{4.42}$$

By (4.40) and (4.42), this theorem is proved. □

**Corollary 4.19**. *Let $T: D \to Y$ be a single-valued mapping. Let $\bar{x} \in D$. Suppose $T$ is Fréchet differentiable at $\bar{x}$. If the following conditions are satisfied,*

(SC) $(X, \tau_X)$ *is seminorm constructed*;
(SB) $\max\{p(v): p \in \mathbb{F}_X\} < \infty$, *for every $v \in X$*,

*then T is Gâteaux differentiable at $\bar{x}$ and*

$$T'(\bar{x})(v) = \nabla T(\bar{x})(v), \text{ for each } v \in X\backslash\{\theta_X\}.$$

*Proof*. This corollary follows from Theorem 4.18 immediately. □

**Corollary 4.20**. *Let $(X, \|\cdot\|_X)$ be a normed vector space equipped with a norm $\|\cdot\|_X$. Let Y be a topological vector space. Let $T: D \to Y$ be a single-valued mapping. Let $\bar{x} \in D$. Suppose T is Fréchet differentiable at $\bar{x}$. Then T is Gâteaux differentiable at $\bar{x}$ and*

$$T'(\bar{x})(v) = \nabla T(\bar{x})(v), \text{ for each } v \in X\backslash\{\theta_X\}.$$

*Proof*. The proof of this corollary is similar to the proofs of Corollaries 11 and 19 that is omitted here. □

### 4.4. Differentiability with Different *F*-seminorm Bases

We first consider Gâteaux differentiability with different *F*-seminorm bases. Notice that in Definition 4.1, the Gâteaux (directional) differentiability of a mapping only depends on the topology of the range space *Y*. So, in next theorem, we do not need to mention precisely the topological structure of the domain space *X*.

**Theorem 4.21**. *Let $(X, \tau_X)$ be a Hausdorff topological vector space and let $(Y, \tau_Y)$ be a Hausdorff topological vector space with two families $\mathbb{F}_Y^1$ and $\mathbb{F}_Y^2$ of positive F-seminorms on Y such that $\tau_Y$ is induced by each one of $\mathbb{F}_Y^1$ and $\mathbb{F}_Y^2$ independently. Let $T: D \to Y$ be a single-valued mapping. Let $\bar{x} \in D$ and $v \in X\backslash\{\theta_X\}$. Let $k \in \{1, 2\}$ be arbitrarily given. If T is Gâteaux directionally differentiable at $\bar{x}$ along the direction v with respect to the topology $\tau_Y$ being induced by $\mathbb{F}_Y^k$, then T is also Gâteaux directionally differentiable at $\bar{x}$ along the direction v with respect to the topology $\tau_Y$ being induced by $\mathbb{F}_Y^{3-k}$ and in both cases the Gâteaux directional derivative $T'(\bar{x}, v)$ at $\bar{x}$ along direction v is the same.*

*Proof.* Suppose that as the topology $\tau_Y$ is induced by $\mathbb{F}_Y^k$, $T$ is Gâteaux directionally differentiable at $\bar{x}$ along the direction $v$ with the Gâteaux directional derivative $T'(\bar{x}, v)$. Let $V_{J,\varepsilon}^{3-k}$ be an arbitrarily given $\tau_Y$-open neighborhood of $\theta_Y$ in *Y*, for some $\varepsilon > 0$ and $J \in \mathcal{F}_Y^{3-k}$. By the topology basis in *Y* generated by $\mathbb{F}_Y^k$, for this $\tau_Y$-open neighborhood $V_{J,\varepsilon}^{3-k}$ of $\theta_Y$ in *Y*, there is $\tau_Y$-open neighborhood $V_{K,\lambda}^k$ of $\theta_Y$ in *Y*, for some $K \in \mathcal{F}_Y^k$ and $\lambda > 0$ such that

$$V_{K,\lambda}^k \subset V_{J,\varepsilon}^{3-k}. \tag{4.43}$$

By the assumption that *T* is Gâteaux directionally differentiable at $\bar{x}$ along the direction *v* with respect to the topology $\tau_Y$ being induced by $\mathbb{F}_Y^k$ and with the Gâteaux directional derivative $T'(\bar{x}, v)$, by Definition (4.1), for this given $\tau_Y$-open neighborhood $V_{K,\lambda}^k$ of $\theta_Y$ in *Y*, there is $\delta > 0$ such that, for any real number *t*, we have that

$$0 < |t| < \delta \quad \Longrightarrow \quad \max\left\{q_k\left(\frac{T(\bar{x}+tv)-T(\bar{x})-tT'(\bar{x},v)}{t}\right): q_k \in K\right\} < \lambda.$$

This is that

$$0 < |t| < \delta \quad \Longrightarrow \quad \frac{T(\bar{x}+tv)-T(\bar{x})-tT'(\bar{x},v)}{t} \in V_{K,\lambda}^k. \tag{4.44}$$

By (4.43), (4.44) implies that

$$0 < |t| < \delta \quad \Longrightarrow \quad \frac{T(\bar{x}+tv)-T(\bar{x})-tT'(\bar{x},v)}{t} \in V_{J,\varepsilon}^{3-k}. \tag{4.45}$$

By definition, (4.45) is equivalent to

$$0 < |t| < \delta \quad \Longrightarrow \quad \max\left\{q_{3-k}\left(\frac{T(\bar{x}+tv)-T(\bar{x})-tT'(\bar{x},v)}{t}\right): q_{3-k} \in J\right\} < \varepsilon.$$

By Definition 4.1, this theorem is proved. □

By Theorem 4.21, we have the following consequence immediately.

**Corollary 4.22**. *Let* $(X, \tau_X)$ *be a Hausdorff topological vector space and let* $(Y, \tau_Y)$ *be a Hausdorff topological vector space, in which there are two families* $\mathbb{F}_Y^1$ *and* $\mathbb{F}_Y^2$ *of positive F-seminorms on Y such that* $\tau_Y$ *is induced by each one of* $\mathbb{F}_Y^1$ *and* $\mathbb{F}_Y^2$ *independently. Let T*: *D* → *Y be a single-valued mapping. Let* $\bar{x} \in D$. *Let* $k \in \{1, 2\}$ *be given. If T is Gâteaux differentiable at* $\bar{x}$ *with respect to the topology* $\tau_Y$ *being induced by* $\mathbb{F}_Y^k$, *then T is also Gâteaux differentiable at* $\bar{x}$ *with respect to the topology* $\tau_Y$ *being induced by* $\mathbb{F}_Y^{3-k}$ *and in both cases the Gâteaux derivative* $T'(\bar{x})$ *at* $\bar{x}$ *is the same.*

Now, we study Fréchet differentiability with different *F*-seminorm bases, which is more complicated than Gâteaux differentiability.

Let $(X, \tau_X)$ and $(Y, \tau_Y)$ be Hausdorff topological vector spaces. Suppose that there are two families $\mathbb{F}_X^1$ and $\mathbb{F}_X^2$ of positive *F*-seminorms on *X* such that $\tau_X$ is induced by each one of $\mathbb{F}_X^1$ and $\mathbb{F}_X^2$ independently; and that there are two families $\mathbb{F}_Y^1$ and $\mathbb{F}_Y^2$ of positive *F*-seminorms on *Y* such that $\tau_Y$ is induced by each one of $\mathbb{F}_Y^1$ and $\mathbb{F}_Y^2$ independently. For *k* = 1, 2, let $\mathcal{F}_X^k$ and $\mathcal{F}_Y^k$ be the sets of nonempty finite subsets of $\mathbb{F}_X^k$ and $\mathbb{F}_Y^k$, respectively. The families $\mathbb{F}_X^k$ and $\mathbb{F}_Y^k$ respectively induce the corresponding *F*-seminorm bases $\{U_{I,\lambda}^k : \lambda > 0, I \in \mathcal{F}_X^k\}$ and $\{V_{J,\lambda}^k : \lambda > 0, J \in \mathcal{F}_Y^k\}$ in *X* and *Y*, for *k* = 1, 2. Let *T*: *D* → *Y* be a single-valued mapping. Let $\bar{x} \in D$.

**Theorem 4.23**. *Let Z* = *X*, *Y*. *Let* $(Z, \tau_Z)$ *be a Hausdorff topological vector space, in which there are two*

*families $\mathbb{F}_Z^1$ and $\mathbb{F}_Z^2$ of positive F-seminorms on Z such that $\tau_Z$ is induced by each one of $\mathbb{F}_Z^1$ and $\mathbb{F}_Z^2$, independently. Let $k \in \{1, 2\}$ be given. Suppose that the following conditions are satisfied.*

(SX) *For any $\tau_X$-open neighborhood $U_{I,\lambda}^k$ of $\theta_X$ in X with respect to $\mathbb{F}_X^k$, there is an $\tau_X$-open neighborhood $U_{I_1,\delta}^{3-k}$ of $\theta_X$ in X with respect to $\mathbb{F}_X^{3-k}$, such that $U_{I_1,\delta}^{3-k} \subset U_{I,\lambda}^k$ and, for any $u \in X$, we have*

$$\max\{p_{3-k}(u): p_{3-k} \in I_1\} = 0 \quad \Leftrightarrow \quad \max\{p_k(u): p_k \in I\} = 0.$$

*Moreover, there is a positive integer M which is only depends on $\mathbb{F}_X^1$ and $\mathbb{F}_X^2$ such that*

$$\frac{\max\{p_k(u): p_k \in I_1\}}{\max\{p_{3-k}(u): p_{3-k} \in I\}} \leq \mathrm{M}, \quad \textit{whenever} \ \max\{p_{3-k}(u): p_{3-k} \in I_1\} > 0.$$

(SY) *For any $\tau_Y$-open neighborhood $V_{J,\varepsilon}^{3-k}$ of $\theta_Y$ in Y with respect to $\mathbb{F}_Y^{3-k}$, there is an $\tau_Y$-open neighborhood $V_{J_1,\lambda}^k$ of $\theta_Y$ in Y with respect to $\mathbb{F}_Y^k$, such that $V_{J_1,\lambda}^k \subset V_{J,\varepsilon}^{3-k}$ and, for any $y \in Y$, we have*

$$\max\{q_{3-k}(y): q_{3-k} \in J_1\} = 0 \quad \Leftrightarrow \quad \max\{q_k(y): q_k \in J\} = 0.$$

*Let T: D → Y be a single-valued mapping. If T is Fréchet differentiable at $\bar{x}$ with $\tau_X$ and $\tau_Y$ being induced by $\mathbb{F}_X^k$ and $\mathbb{F}_Y^k$, respectively, then as $\tau_X$ and $\tau_Y$ being respectively induced by $\mathbb{F}_X^{3-k}$ and $\mathbb{F}_Y^{3-k}$, T is also Fréchet differentiable at $\bar{x}$ and in both cases the Fréchet derivative $\nabla T(\bar{x})$ at $\bar{x}$ is the same.*

*Let T: D → Y be a single-valued mapping. If T is Fréchet differentiable at $\bar{x}$ with $\tau_X$ and $\tau_Y$ being induced by $\mathbb{F}_X^k$ and $\mathbb{F}_Y^k$, respectively, then since $\tau_X$ and $\tau_Y$ are also respectively induced by $\mathbb{F}_X^{3-k}$ and $\mathbb{F}_Y^{3-k}$, T is Fréchet differentiable at $\bar{x}$ in that sense as well and has the same Fréchet derivative.*

*Proof.* By the assumption in this theorem, let $\nabla T(\bar{x})$ be the Fréchet derivative of *T* at $\bar{x}$ with respect to $\tau_X$ being induced by $\mathbb{F}_X^k$ and $\tau_Y$ being induced by $\mathbb{F}_Y^k$. Let $V_{J,\varepsilon}^{3-k}$ be an arbitrarily given $\tau_Y$-open neighborhood of $\theta_Y$ in *Y*, for some $1 > \varepsilon > 0$ and $J \in \mathcal{F}_Y^{3-k}$. We consider $V_{J,\frac{\varepsilon}{M}}^{3-k}$. Since $\mathbb{F}_Y^1$ and $\mathbb{F}_Y^2$ satisfy condition (SY), there is $\tau_Y$-open neighborhood $V_{K,\lambda}^k$ of $\theta_Y$ in *Y*, for some $K \in \mathcal{F}_Y^k$ and $\lambda > 0$ such that $V_{K,\lambda}^k \subset V_{J,\varepsilon}^{3-k}$ and, for any $y \in Y$, we have

$$\max\{q_k(y): p_k \in K\} = 0 \quad \Leftrightarrow \quad \max\{q_{3-k}(y): q_{3-k} \in J\} = 0; \tag{4.46}$$

$$0 < \max\{q_k(y): p_k \in K\} < \lambda \quad \Longrightarrow \quad 0 < \max\{q_{3-k}(y): q_{3-k} \in J\} < \frac{\varepsilon}{M}. \tag{4.47}$$

By the condition in this theorem that $\nabla T(\bar{x})$ is the Fréchet derivative of *T* at $\bar{x}$ with respect to $\tau_X$ being induced by $\mathbb{F}_X^k$ and $\tau_Y$ being induced by $\mathbb{F}_Y^k$, for the $\tau_Y$-open neighborhood $V_{K,\lambda}^k$ of $\theta_Y$ in *Y*, there is an $\tau_X$-open neighborhood $U_{I_1,\delta_1}^k$ of $\theta_X$ in *X*, for some $\delta_1 > 0$ and $I_1 \in \mathcal{F}_X^k$ such that, for $u \in U_{I_1,\delta_1}^k$, *u* satisfies the conditions (DZ) and (DR) in Definition 4.8 as follows

(i) $\max\{p_k(u): p_k \in I_1\} = 0 \implies \max\left\{q_k\left(T(\bar{x}+u) - T(\bar{x}) - \nabla T(\bar{x})(u)\right): q_k \in K\right\} = 0;$

(ii) $0 < \max\{p_k(u): p_k \in I_1\} < \delta_1 \implies \max\left\{q_k\left(\frac{T(\bar{x}+u)-T(\bar{x})-\nabla T(\bar{x})(u)}{\max\{p_k(u): p_k \in I_1\}}\right): q_k \in K\right\} < \lambda.$ (4.48)

Since $\mathbb{F}_X^k$ and $\mathbb{F}_X^{3-k}$ satisfy condition (SX), for the $\tau_X$-open neighborhood $U_{I_1,\delta_1}^k$ of $\theta_X$ in *X*, there is an $\tau_X$- open neighborhood $U_{I,\delta}^{3-k}$ of $\theta_X$ such that $U_{I,\delta}^{3-k} \subset U_{I_1,\delta_1}^{k}$ and, for any $u \in U_{I,\delta}^{3-k}$, we have

$$\max\{p_{3-k}(u): p_{3-k} \in I\} = 0 \quad \Leftrightarrow \quad \max\{p_k(u): p_k \in I_1\} = 0; \tag{4.49}$$

and $\max\{p_{3-k}(u): p_{3-k} \in I\} > 0$

$$\Longrightarrow \quad \max\{p_k(u): p_k \in I_1\} > 0 \quad \text{and} \quad \frac{\max\{p_k(u): p_k \in I_1\}}{\max\{p_{3-k}(u): p_{3-k} \in I\}} \leq \mathrm{M}, \tag{4.50}$$

Then, by (4.49), above (i) and (4.46), for any $u \in X$, we have

$$\max\{p_{3-k}(u): p_{3-k} \in I\} = 0$$

$$\Leftrightarrow \quad \max\{p_k(u): p_k \in I_1\} = 0$$

$$\Longrightarrow \quad \max\left\{q_k\big(T(\bar{x}+u) - T(\bar{x}) - \nabla T(\bar{x})(u)\big): q_k \in K\right\} = 0$$

$$\Longrightarrow \quad \max\left\{q_{3-k}\big(T(\bar{x}+u) - T(\bar{x}) - \nabla T(\bar{x})(u)\big): q_{3-k} \in J\right\} = 0.$$

This implies that $\nabla T(\bar{x})$ satisfies condition (DZ) in Definition 4.8. By (4.50), (4.48) and (4.47), for any $u \in U_{I,\delta}^{3-k} \subset U_{I_1,\delta_1}^{k}$, we have

$$0 < \max\{p_{3-k}(u): p_{3-k} \in I\} < \delta$$

$$\Longrightarrow \quad 0 < \max\{p_k(u): p_k \in I_1\} < \delta_1$$

$$\Longrightarrow \quad 0 < \max\left\{q_k\left(\frac{T(\bar{x}+u) - T(\bar{x}) - \nabla T(\bar{x})(u)}{\max\{p_k(u): p_k \in I_1\}}\right): q_k \in K\right\} < \lambda$$

$$\Longrightarrow \max\left\{q_{3-k}\left(\frac{T(\bar{x}+u) - T(\bar{x}) - \nabla T(\bar{x})(u)}{\max\{p_k(u): p_k \in I_1\}}\right): q_{3-k} \in J\right\} < \frac{\varepsilon}{M}. \tag{4.51}$$

By (4.50) and (4.51), for any $u \in U_{I,\delta}^{3-k} \subset U_{I_1,\delta_1}^{k}$, if

$$0 < \max\{p_{3-k}(u): p_{3-k} \in I\} < \delta,$$

then

$$\max\left\{q_{3-k}\left(\frac{T(\bar{x}+u) - T(\bar{x}) - \nabla T(\bar{x})(u)}{\max\{p_{3-k}(u): p_{3-k} \in I\}}\right): q_{3-k} \in J\right\}$$

$$= \max\left\{q_{3-k}\left(\frac{\max\{p_k(u): p_k \in I_1\}}{\max\{p_{3-k}(u): p_{3-k} \in I\}} \, \frac{T(\bar{x}+u) - T(\bar{x}) - \nabla T(\bar{x})(u)}{\max\{p_k(u): p_k \in I_1\}}\right): q_{3-k} \in J\right\}$$

$$\leq \max\left\{q_{3-k}\left(M \, \frac{T(\bar{x}+u) - T(\bar{x}) - \nabla T(\bar{x})(u)}{\max\{p_k(u): p_k \in I_1\}}\right): q_{3-k} \in J\right\}$$

$$\leq M \max\left\{q_{3-k}\left(\frac{T(\bar{x}+u) - T(\bar{x}) - \nabla T(\bar{x})(u)}{\max\{p_k(u): p_k \in I_1\}}\right): q_{3-k} \in J\right\}$$

$$< \varepsilon.$$

This proves that $\nabla T(\bar{x})$ is also the Fréchet derivative of $T$ at $\bar{x}$ with respect to $\tau_X$ being induced by $\mathbb{F}_X^{3-k}$ and $\tau_Y$ being induced by $\mathbb{F}_Y^{3-k}$. □

## 5. Differentiation in Schwartz Space

### 5.1. Review Schwartz Space Equipped with a Countable Family of Seminorms

Schwartz space is a special class of complete locally convex topological vector spaces (It is indeed a Fréchet space). It is well known that Schwartz space has been a very important space. It has been widely used in functional analysis (see Reed and Simon [37]), in partial differential equations (see Hörmander [15]), in Fourier analysis (see [18, 30, 39, 41]), and so forth. In this section, it is natural for us to use the results of the previous section to find explicit formulas of Gâteaux and Fréchet derivatives of some single-valued mappings in Schwartz space.

In this subsection, we review the concepts of Schwartz space and the properties of a countable family of seminorms the space is equipped with. Recall that in this paper, $\mathbb{N}$ denotes the set of nonnegative integers. For any $n \in \mathbb{N}\backslash\{0\}$, let $\mathbb{N}^n \coloneqq \mathbb{N} \times \ldots \times \mathbb{N}$ be the *n*-fold Cartesian product of $\mathbb{N}$. Let $\theta_n = (0, \ldots, 0)$ be the origin in $\mathbb{N}^n$. Let $C^\infty(\mathbb{R}^n, \mathbb{C})$ be the function space of smooth functions from $\mathbb{R}^n$ into complex $\mathbb{C}$. Let $\preccurlyeq$ be the pointwise partial order in $\mathbb{N}^n$, which is defined, for $\alpha = (\alpha_1, \ldots, \alpha_n)$ and $\beta = (\beta_1, \ldots, \beta_n) \in \mathbb{N}^n$, by

$$\alpha \preccurlyeq \beta \text{ if and only if } \alpha_i \leq \beta_i, \text{ for } i = 1, 2, \ldots, n.$$

We write $\alpha + \beta = (\alpha_1 + \beta_1, \ldots, \alpha_n + \beta_n)$ and $|\alpha| = \alpha_1 + \cdots + \alpha_n$. Then, for $f \in C^\infty(\mathbb{R}^n, \mathbb{C})$, we use the following multi-index notations:

$$x^\alpha = x_1^{\alpha_1} \ldots x_n^{\alpha_n}, \text{ for any } x = (x_1, \ldots, x_n) \in \mathbb{R}^n,$$

and

$$D^\beta f = \partial_1^{\beta_1} \ldots \partial_n^{\beta_n} f, \text{ for any } f \in C^\infty(\mathbb{R}^n, \mathbb{C}).$$

In particular, for $\alpha = \beta = \theta_n$, we have $x^{\theta_n} = 1$ and $D^{\theta_n} f = f$, for any $f \in C^\infty(\mathbb{R}^n, \mathbb{C})$. Then, for any $f \in C^\infty(\mathbb{R}^n, \mathbb{C})$, we write

$$\|f\|_{\alpha,\beta} = \sup_{x \in \mathbb{R}^n} \left| x^\alpha (D^\beta f)(x) \right|.$$

In particular, we have

$$\|f\|_{\theta_n,\theta_n} = \sup_{x \in \mathbb{R}^n} |f(x)|, \text{ for any } f \in C^\infty(\mathbb{R}^n, \mathbb{C}).$$

The Schwartz space (the space of rapidly decreasing functions on $\mathbb{R}^n$) is denoted by $\mathcal{S}(\mathbb{R}^n, \mathbb{C})$, or simply by $\mathcal{S}$, which is a subspace of $C^\infty(\mathbb{R}^n, \mathbb{C})$ and it is defined by

$$\mathcal{S}(\mathbb{R}^n, \mathbb{C}) = \left\{ f \in C^\infty(\mathbb{R}^n, \mathbb{C}) : \|f\|_{\alpha,\beta} < \infty, \text{ for any } \alpha, \beta \in \mathbb{N}^n \right\}.$$

The origin of $\mathcal{S}(\mathbb{R}^n, \mathbb{C})$ is denoted by $\theta_{\mathcal{S}}$. Then, we obtain a countable family $\left\{ \|\cdot\|_{\alpha,\beta} : \alpha, \beta \in \mathbb{N}^n \right\}$ of seminorms on $\mathcal{S}(\mathbb{R}^n, \mathbb{C})$, which is denoted by

$$\mathbb{F}_{\mathcal{S}} = \left\{ \|\cdot\|_{\alpha,\beta} : \alpha, \beta \in \mathbb{N}^n \right\}.$$

The collection of nonempty finite subsets $\mathcal{F}_{\mathcal{S}}$ of $\mathbb{F}_{\mathcal{S}}$ is written as

$$\mathcal{F}_{\mathcal{S}} = \{\{\|\cdot\|_{\alpha,\beta} \in \mathbb{F}_{\mathcal{S}} : (\alpha, \beta) \in \mathbb{I}\}, \mathbb{I} \text{ is nonempty finite subset of } \mathbb{N} \times \mathbb{N}\}.$$

Let $\tau_{\mathcal{S}}$ be the topology in $\mathcal{S}(\mathbb{R}^n, \mathbb{C})$ induced by the countable family $\mathbb{F}_{\mathcal{S}} = \{\|\cdot\|_{\alpha,\beta} : \alpha, \beta \in \mathbb{N}^n\}$ of seminorms on $\mathcal{S}(\mathbb{R}^n, \mathbb{C})$. Then, it is known that $(\mathcal{S}(\mathbb{R}^n, \mathbb{C}), \tau_{\mathcal{S}})$ is a complete Hausdorff locally convex topological vector space. Recall that, for any given $\lambda > 0$, $I \in \mathcal{F}_{\mathcal{S}}$, following definition (2.2), we write

$$U_{I,\lambda} = \{f \in \mathcal{S}(\mathbb{R}^n, \mathbb{C}) : \max\{\|f\|_{\alpha,\beta} : \|\cdot\|_{\alpha,\beta} \in I\} < \lambda\}$$

$$= \left\{f \in \mathcal{S}(\mathbb{R}^n, \mathbb{C}) : \max\left\{\sup_{x \in \mathbb{R}^n} \left|x^{\alpha}(D^{\beta} f)(x)\right| : \|\cdot\|_{\alpha,\beta} \in I\right\} < \lambda\right\}.$$

$U_{I,\lambda}$ is an $\tau_{\mathcal{S}}$-open neighborhood of $(\mathcal{S}(\mathbb{R}^n, \mathbb{C}), \tau_{\mathcal{S}})$ around $\theta_{\mathcal{S}}$. By (2.1) and (2.3), the collection

$$U(\mathbb{F}_{\mathcal{S}}) = \{U_{I,\lambda} : I \in \mathcal{F}_{\mathcal{S}}, \lambda > 0\},$$

forms an open neighborhood basis of $\mathcal{S}(\mathbb{R}^n, \mathbb{C})$ around $\theta_{\mathcal{S}}$. More generally, for any $g \in \mathcal{S}(\mathbb{R}^n, \mathbb{C})$, $\lambda > 0$, and $I \in \mathcal{F}_{\mathcal{S}}$, we write

$$U_{I,\lambda}(g) = \{f \in \mathcal{S}(\mathbb{R}^n, \mathbb{C}) : \max\{\|f - g\|_{\alpha,\beta} : \|\cdot\|_{\alpha,\beta} \in I\} < \lambda\}$$

$$= \left\{f \in \mathcal{S}(\mathbb{R}^n, \mathbb{C}) : \max\left\{\sup_{x \in \mathbb{R}^n} \left|x^{\alpha}\left(D^{\beta}(f - g)(x)\right| : \|\cdot\|_{\alpha,\beta} \in I\right\} < \lambda\right\}$$

$$= \left\{f \in \mathcal{S}(\mathbb{R}^n, \mathbb{C}) : \max\left\{\sup_{x \in \mathbb{R}^n} \left|x^{\alpha}\left(D^{\beta} f\right)(x) - x^{\alpha}(D^{\beta} g)(x)\right| : \|\cdot\|_{\alpha,\beta} \in I\right\} < \lambda\right\}.$$

$U_{I,\lambda}(g)$ is an $\tau_{\mathcal{S}}$-open neighborhood of $\mathcal{S}(\mathbb{R}^n, \mathbb{C})$ around $g$. By (2.3), the following collection

$$U(\mathbb{F}_{\mathcal{S}})(g) = \{U_{I,\lambda}(g) : I \in \mathcal{F}_{\mathcal{S}}, \lambda > 0\},$$

forms an open neighborhood basis of $\mathcal{S}(\mathbb{R}^n, \mathbb{C})$ around $g$. In more details, the Schwartz space $(\mathcal{S}(\mathbb{R}^n, \mathbb{C}), \tau_{\mathcal{S}})$ has the following properties (see [16, 41, 46]).

(i) For any $\alpha, \beta \in \mathbb{N}^n$, $\|\cdot\|_{\alpha,\beta}$ is a seminorm on $\mathcal{S}(\mathbb{R}^n, \mathbb{C})$;

(ii) $(\mathcal{S}(\mathbb{R}^n, \mathbb{C}), \tau_{\mathcal{S}})$ is a complete Hausdorff locally convex topological vector space and its topology is induced by the countable family $\mathbb{F}_{\mathcal{S}} = \{\|\cdot\|_{\alpha,\beta} : \alpha, \beta \in \mathbb{N}^n\}$ of seminorms. That is, $(\mathcal{S}(\mathbb{R}^n, \mathbb{C}), \tau_{\mathcal{S}})$ is a Fréchet space. So, the space $(\mathcal{S}(\mathbb{R}^n, \mathbb{C}), \tau_{\mathcal{S}})$ has a seminorm construction;

(iii) The topology $\tau_{\mathcal{S}}$ of $\mathcal{S}(\mathbb{R}^n, \mathbb{C})$ is not defined from a norm.

(iv) $\mathcal{S}(\mathbb{R}^n, \mathbb{C})$ is metrizable; for example, we define $d : \mathcal{S} \times \mathcal{S} \to \mathbb{R}_+$ by

$$d(f, g) = \sum_{\alpha,\beta \in \mathbb{N}^n} a_{\alpha,\beta} \frac{\|f - g\|_{\alpha,\beta}}{1 + \|f - g\|_{\alpha,\beta}}, \text{ for any } f, g \in \mathcal{S}(\mathbb{R}^n, \mathbb{C}).$$

Where, $a_{\alpha,\beta} > 0$ are such that $\sum_{\alpha,\beta \in \mathbb{N}^n} a_{\alpha,\beta} < \infty$. Then, $d$ is a metric on $\mathcal{S}$ which generates the topology $\tau_{\mathcal{S}}$, such that $\mathcal{S}$ is complete with respect to this metric $d$.

(v) $f, g \in \mathcal{S}(\mathbb{R}^n, \mathbb{C}) \Longrightarrow fg \in \mathcal{S}(\mathbb{R}^n, \mathbb{C})$;

(vi) Let $f \in \mathcal{S}(\mathbb{R}^n, \mathbb{C})$ and $h \in C^{\infty}(\mathbb{R}^n, \mathbb{C})$. If $h$ has bounded derivatives of all orders, then, $hf \in \mathcal{S}(\mathbb{R}^n, \mathbb{C})$.

### 5.2. Some Continuous and Linear Operators in Schwartz space

In this subsection, we study some continuous and linear operators in Schwartz space $\mathcal{S}(\mathbb{R}^n, \mathbb{C})$. By Proposition 4.15, we find explicit formulas of their Fréchet derivatives, which are listed as the following lemmas.

**Lemma 5.1**. *For any* $\alpha, \beta, \gamma \in \mathbb{N}^n$, *and for any* $f \in \mathcal{S}(\mathbb{R}^n, \mathbb{C})$, *we have*

(i) $D^{\beta}(D^{\gamma}f) = D^{\beta+\gamma}f = \partial_1^{\beta_1+\gamma_1} \dots \partial_n^{\beta_n+\gamma_n} f$;

(ii) $\|D^{\gamma}f\|_{\alpha,\beta} = \|f\|_{\alpha,\beta+\gamma}$.

*Proof*. The proof of this lemma is straight forward and it is omitted here. □

**Lemma 5.2**. *For any* $\gamma \in \mathbb{N}^n$, *let* $\mathfrak{D}^{\gamma}\colon \mathcal{S}(\mathbb{R}^n, \mathbb{C}) \to \mathcal{S}(\mathbb{R}^n, \mathbb{C})$ *be the differential operator of order* $\gamma$ *with*

$$\mathfrak{D}^{\gamma}(f) = D^{\gamma}f, \text{ for any } f \in \mathcal{S}(\mathbb{R}^n, \mathbb{C}).$$

*Then*, *we have*

(i) $\mathfrak{D}^{\gamma}$ *is a continuous and linear mapping on* $\mathcal{S}(\mathbb{R}^n, \mathbb{C})$;

(ii) $\mathfrak{D}^{\gamma}$ *is Fréchet differentiable on* $\mathcal{S}(\mathbb{R}^n, \mathbb{C})$ *such that*, *for any* $\bar{f} \in \mathcal{S}(\mathbb{R}^n, \mathbb{C})$, *the Fréchet derivative of* $\mathfrak{D}^{\gamma}$ *at point* $\bar{f}$ *satisfies that*

$$\nabla(\mathfrak{D}^{\gamma})(\bar{f})(u) = \mathfrak{D}^{\gamma}(u) = D^{\gamma}u, \text{ for every } u \in \mathcal{S}(\mathbb{R}^n, \mathbb{C}).$$

*Proof*. Proof of (i). The linearity of $\mathfrak{D}^{\gamma}\colon \mathcal{S}(\mathbb{R}^n, \mathbb{C}) \to \mathcal{S}(\mathbb{R}^n, \mathbb{C})$ is clear. To prove the continuity of $\mathfrak{D}^{\gamma}$, for any seminorm $\|\cdot\|_{\alpha,\beta}$ with arbitrarily given $\alpha, \beta \in \mathbb{N}^n$, by part (ii) of Lemma 5.1, we have

$$\|\mathfrak{D}^{\gamma}(f)\|_{\alpha,\beta} = \|D^{\gamma}f\|_{\alpha,\beta} = \|f\|_{\alpha,\beta+\gamma}, \text{ for any } f \in \mathcal{S}(\mathbb{R}^n, \mathbb{C}).$$

In part (v) of Theorem 3.9, let $C_{\alpha,\beta} = 1$. This proves the continuity of $\mathfrak{D}^{\gamma}\colon \mathcal{S}(\mathbb{R}^n, \mathbb{C}) \to \mathcal{S}(\mathbb{R}^n, \mathbb{C})$. After part (i) of this lemma is proved, part (ii) follows from part (i) and Proposition 4.15 immediately. □

**Lemma 5.3**. *Let* $g \in \mathcal{S}(\mathbb{R}^n, \mathbb{C})$. *Let* $T_g\colon \mathcal{S}(\mathbb{R}^n, \mathbb{C}) \to \mathcal{S}(\mathbb{R}^n, \mathbb{C})$ *be a multiplication operator defined by*

$$T_g(f) = gf, \text{ for any } f \in \mathcal{S}(\mathbb{R}^n, \mathbb{C}).$$

*Then*, *we have*

(i) $T_g$ *is a continuous and linear operator on* $(\mathcal{S}(\mathbb{R}^n, \mathbb{C}), \tau_{\mathcal{S}})$.

(ii) $T_g$ *is Fréchet differentiable on* $\mathcal{S}(\mathbb{R}^n, \mathbb{C})$ *such that*, *for any* $\bar{f} \in \mathcal{S}(\mathbb{R}^n, \mathbb{C})$, *the Fréchet derivative of* $T_g$ *at point* $\bar{f}$ *satisfies that*

$$\nabla T_g(\bar{f})(u) = T_g(u), \text{ for every } u \in \mathcal{S}(\mathbb{R}^n, \mathbb{C}).$$

*Proof*. Proof of (i). By property (vi) of the Schwartz space $\mathcal{S}(\mathbb{R}^n, \mathbb{C})$, the operator $T_g\colon \mathcal{S}(\mathbb{R}^n, \mathbb{C}) \to \mathcal{S}(\mathbb{R}^n, \mathbb{C})$ is well-defined and the linearity of $T_g$ is clear. Next, we prove the continuity of $T_g$. Let $\beta = (\beta_1, \dots, \beta_n) \in \mathbb{N}^n$. We can show that

$$\partial_n^{\beta_n}(gf) = \sum_{k_n=0}^{\beta_n} \binom{\beta_n}{k_n} \partial_n^{\beta_n - k_n} g \partial_n^{k_n} f, \text{ for any } f \in \mathcal{S}(\mathbb{R}^n, \mathbb{C}).$$

This implies that, for any $f \in \mathcal{S}(\mathbb{R}^n, \mathbb{C})$, we have

$$\begin{aligned}
&\partial_{n-1}^{\beta_{n-1}} \partial_n^{\beta_n}(gf) \\
&= \partial_{n-1}^{\beta_{n-1}} \left( \sum_{k=0}^{\beta_n} \binom{\beta_n}{k_n} (\partial_n^{\beta_n - k_n} g \partial_n^{k_n} f) \right) \\
&= \sum_{k=0}^{\beta_n} \binom{\beta_n}{k_n} \partial_{n-1}^{\beta_{n-1}} (\partial_n^{\beta_n - k_n} g \partial_n^{k_n} f) \\
&= \sum_{k_n=0}^{\beta_n} \binom{\beta_n}{k_n} \sum_{k_{n-1}=0}^{\beta_{n-1}} \binom{\beta_{n-1}}{k_{n-1}} \partial_{n-1}^{\beta_{n-1} - k_{n-1}} \left( \partial_n^{\beta_n - k_n} g \right) \partial_{n-1}^{k_{n-1}} (\partial_n^{k_n} f) \\
&= \sum_{k_n=0}^{\beta_n} \binom{\beta_n}{k_n} \sum_{k_{n-1}=0}^{\beta_{n-1}} \binom{\beta_{n-1}}{k_{n-1}} \partial_{n-1}^{\beta_{n-1} - k_{n-1}} \left( \partial_n^{\beta_n - k_n} g \right) \partial_{n-1}^{k_{n-1}} (\partial_n^{k_n} f) \\
&= \sum_{k_n=0}^{\beta_n} \binom{\beta_n}{k_n} \sum_{k_{n-1}=0}^{\beta_{n-1}} \binom{\beta_{n-1}}{k_{n-1}} \left( \partial_{n-1}^{\beta_{n-1} - k_{n-1}} \partial_n^{\beta_n - k_n} g \right) (\partial_{n-1}^{k_{n-1}} \partial_n^{k_n} f). \qquad (5.1)
\end{aligned}$$

By repeating the above process (5.1), for any $f \in \mathcal{S}(\mathbb{R}^n, \mathbb{C})$, we obtain that,

$$\partial_1^{\beta_1} \dots \partial_n^{\beta_n}(gf) = \sum_{k_n=0}^{\beta_n} \binom{\beta_n}{k_n} \dots \sum_{k_1=0}^{\beta_1} \binom{\beta_1}{k_1} \left( \partial_1^{\beta_1 - k_1} \dots \partial_n^{\beta_n - k_n} g \right) \left( \partial_1^{k_1} \dots \partial_n^{k_n} f \right). \qquad (5.2)$$

Then, $\partial_1^{\beta_1} \dots \partial_n^{\beta_n}(gf)$ is considered to be a "polynomial" of $\partial_1^{k_1} \dots \partial_n^{k_n} f$ with coefficients of some constants multiplied by functions $\partial_1^{\beta_1 - k_1} \dots \partial_n^{\beta_n - k_n} g$, which has maximum $\prod_{i=1}^n (\beta_i + 1)$ terms. Let $(\alpha, \beta) \in \mathbb{N}^n \times \mathbb{N}^n$ and let $k = (k_1, \dots, k_n) \in \mathbb{N}^n$ runs over all possible cases $\theta_n \preccurlyeq k \preccurlyeq \beta$. Then we have

$$\begin{aligned}
&\left\| T_g f \right\|_{\alpha,\beta} = \| gf \|_{\alpha,\beta} \\
&= \sup_{x \in \mathbb{R}^n} \left| x^\alpha (D^\beta g f)(x) \right| \\
&= \sup_{x \in \mathbb{R}^n} \left| x^\alpha \left( \sum_{k_n=0}^{\beta_n} \binom{\beta_n}{k_n} \dots \sum_{k_1=0}^{\beta_1} \binom{\beta_1}{k_1} \left( \partial_1^{\beta_1 - k_1} \dots \partial_n^{\beta_n - k_n} g \right) \left( \partial_1^{k_1} \dots \partial_n^{k_n} f \right) \right) (x) \right| \\
&\leq \sup_{x \in \mathbb{R}^n} \left( \sum_{k_n=0}^{\beta_n} \binom{\beta_n}{k_n} \dots \sum_{k_1=0}^{\beta_1} \binom{\beta_1}{k_1} \left| \left( \partial_1^{\beta_1 - k_1} \dots \partial_n^{\beta_n - k_n} g \right) (x) \right| \left| x^\alpha \left( \partial_1^{k_1} \dots \partial_n^{k_n} f \right)(x) \right| \right) \\
&\leq \sum_{k_n=0}^{\beta_n} \binom{\beta_n}{k_n} \dots \sum_{k_1=0}^{\beta_1} \binom{\beta_1}{k_1} \| g \|_{\theta_n, \beta - k} \| f \|_{\alpha, k}. \qquad (5.\ 3)
\end{aligned}$$

Notice that $g \in \mathcal{S}(\mathbb{R}^n, \mathbb{C})$ is given satisfying $\| g \|_{\theta_n, \beta - k} < \infty$, for any $k \in \mathbb{N}^n$ with $\theta_n \preccurlyeq k \preccurlyeq \beta$. By part (v) in Theorem 3.9 and (5.3), it proves that $T_g \colon \mathcal{S}(\mathbb{R}^n, \mathbb{C}) \to \mathcal{S}(\mathbb{R}^n, \mathbb{C})$ is continuous. After part (i) of this lemma is proved, part (ii) follows from part (i) and Proposition 4.15 immediately. □

Notice that in property (vi) of the Schwartz space $\mathcal{S}(\mathbb{R}^n, \mathbb{C})$, for $h \in C^\infty(\mathbb{R}^n, \mathbb{C})$, the conditions that $h$ has bounded derivatives of all orders are sufficient conditions to have $hf \in \mathcal{S}(\mathbb{R}^n, \mathbb{C})$, for any $f \in \mathcal{S}(\mathbb{R}^n, \mathbb{C})$. As a matter of fact, we will show that, for any polynomial function $P \in C^\infty(\mathbb{R}^n, \mathbb{C})$ with any given order, we have $Pf \in \mathcal{S}(\mathbb{R}^n, \mathbb{C})$, for any $f \in \mathcal{S}(\mathbb{R}^n, \mathbb{C})$. For this purpose, we start with monomials.

**Lemma 5.4**. *For any* $\lambda = (\lambda_1, \dots, \lambda_n) \in \mathbb{N}^n$, *define an operator* $T_\lambda$ *on* $\mathcal{S}(\mathbb{R}^n, \mathbb{C})$ *by*

$$T_\lambda(f) = x^\lambda f, \textit{ for any } f \in \mathcal{S}(\mathbb{R}^n, \mathbb{C}).$$

*Then,*

(i) *$T_\lambda$ is a $\tau_S$-continuous linear operator on $\mathcal{S}(\mathbb{R}^n, \mathbb{C})$ satisfying*

$$x^\lambda f \in \mathcal{S}(\mathbb{R}^n, \mathbb{C}), \textit{ for any } f \in \mathcal{S}(\mathbb{R}^n, \mathbb{C}).$$

(ii) *$T_\lambda$ is Fréchet differentiable on $\mathcal{S}(\mathbb{R}^n, \mathbb{C})$ such that, for any $\bar{f} \in \mathcal{S}(\mathbb{R}^n, \mathbb{C})$, the Fréchet derivative of $T_\lambda$ at point $\bar{f}$ satisfies that*

$$\nabla T_\lambda(\bar{f})(u) = T_\lambda(u) = x^\lambda u, \textit{ for every } u \in \mathcal{S}(\mathbb{R}^n, \mathbb{C}).$$

*Proof*. Proof of (i). Let $\beta = (\beta_1, \dots, \beta_n) \in \mathbb{N}^n$. By (5.2) in the proof of Lemma 5.3, we have

$$\partial_1^{\beta_1} \dots \partial_n^{\beta_n}\left(x^\lambda f\right) = \sum_{k_n=0}^{\beta_n} \binom{\beta_n}{k_n} \dots \sum_{k_1=0}^{\beta_1} \binom{\beta_1}{k_1} \left(\partial_1^{\beta_1-k_1} \dots \partial_n^{\beta_n-k_n} x^\lambda\right)\left(\partial_1^{k_1} \dots \partial_n^{k_n} f\right).$$

Let $\mathbb{Z}$ be the set of all integers. For $i, j \in \mathbb{Z}$, let $i \vee j = \max\{i, j\}$ and $i \wedge j = \min\{i, j\}$. More generally, for $\beta = (\beta_1, \dots, \beta_n) \in \mathbb{N}^n$ and $\lambda = (\lambda_1, \dots, \lambda_n) \in \mathbb{N}^n$, we write

$$\beta \vee \lambda = (\beta_1 \vee \lambda_1, \dots, \beta_n \vee \lambda_n) \quad \text{and} \quad \beta \wedge \lambda = (\beta_1 \wedge \lambda_1, \dots, \beta_n \wedge \lambda_n).$$

Then, for seminorm $\|\cdot\|_{\alpha,\beta}$ with arbitrary $\alpha, \beta \in \mathbb{N}^n$ and for any $f \in \mathcal{S}(\mathbb{R}^n, \mathbb{C})$, similar to (5.3), we have

$$\|T_\lambda f\|_{\alpha,\beta} = \left\|x^\lambda f\right\|_{\alpha,\beta}$$

$$= \sup_{x \in \mathbb{R}^n} \left|x^\alpha\left(D^\beta\left(x^\lambda f\right)(x)\right)\right|$$

$$\leq \sup_{x \in \mathbb{R}^n} \left(\sum_{k_n=0}^{\beta_n} \binom{\beta_n}{k_n} \dots \sum_{k_1=0}^{\beta_1} \binom{\beta_1}{k_1} \left|\partial_1^{\beta_1-k_1} \dots \partial_n^{\beta_n-k_n} x^\lambda\right| \left|x^\alpha\left(\partial_1^{k_1} \dots \partial_n^{k_n} f\right)(x)\right|\right)$$

$$= \sup_{x \in \mathbb{R}^n} \left(\sum_{k_n=0 \vee (\beta_n-\lambda_n)}^{\beta_n} \binom{\beta_n}{k_n} \dots \sum_{k_1=0 \vee (\beta_1-\lambda_1)}^{\beta_1} \binom{\beta_1}{k_1} \left|\partial_1^{\beta_1-k_1} \dots \partial_n^{\beta_n-k_n} x^\lambda\right| \left|x^\alpha\left(\partial_1^{k_1} \dots \partial_n^{k_n} f\right)(x)\right|\right)$$

$$= \sup_{x \in \mathbb{R}^n} \left(\sum_{k_n=0 \vee (\beta_n-\lambda_n)}^{\beta_n} \binom{\beta_n}{k_n} \dots \sum_{k_1=0 \vee (\beta_1-\lambda_1)}^{\beta_1} \binom{\beta_1}{k_1} \left|A_k x_1^{\lambda_1-\beta_1+k_1} \dots x_n^{\lambda_n-\beta_n+k_n}\right| \left|x^\alpha\left(\partial_1^{k_1} \dots \partial_n^{k_n} f\right)(x)\right|\right)$$

$$= \sup_{x \in \mathbb{R}^n} \left(\sum_{k_n=0 \vee (\beta_n-\lambda_n)}^{\beta_n} \binom{\beta_n}{k_n} \dots \sum_{k_1=0 \vee (\beta_1-\lambda_1)}^{\beta_1} \binom{\beta_1}{k_1} \left|A_k x_1^{\alpha_1+\lambda_1-\beta_1+k_1} \dots x_n^{\alpha_n+\lambda_n-\beta_n+k_n}\left(\partial_1^{k_1} \dots \partial_n^{k_n} f\right)(x)\right|\right).$$

Here, $k = (k_1, \dots, k_n) \in \mathbb{N}^n$, which satisfies that $0 \vee (\beta - \lambda) \preccurlyeq k \preccurlyeq \beta$, and $A_k$ is a positive number, which only depends on $k$ such that, for

$$A_k = \Pi_{i=1}^n \left(\lambda_i(\lambda_i - 1) \dots \left(\lambda_i - \left(\beta_i - k_i\right) + 1\right)\right), \text{ for any } \lambda = (\lambda_1, \dots, \lambda_n) \in \mathbb{N}^n \text{ with } 0 \vee (\beta - \lambda) \preccurlyeq k \preccurlyeq \beta.$$

Then, we have that

$$\|T_\lambda f\|_{\alpha,\beta}$$

$$\leq \sup_{x \in \mathbb{R}^n} \left(\sum_{k_n=0 \vee (\beta_n \wedge \lambda_n)}^{\beta_n} \binom{\beta_n}{k_n} \dots \sum_{k_1=0 \vee (\beta_1-\lambda_1)}^{\beta_1} \binom{\beta_1}{k_1} \left|A_k x_1^{\alpha_1+\lambda_1-\beta_1+k_1} \dots x_n^{\alpha_n+\lambda_n-\beta_n+k_n}\left(\partial_1^{k_1} \dots \partial_n^{k_n} f\right)(x)\right|\right)$$

$$= \sup_{x\in\mathbb{R}^n} \left( \sum_{k_n=0\vee(\beta_n\wedge\lambda_n)}^{\beta_n} \binom{\beta_n}{k_n} \cdots \sum_{k_1=0\vee(\beta_1-\lambda_1)}^{\beta_1} \binom{\beta_1}{k_1} \left| A_k t^{\alpha+\lambda-\beta+k} \left(D^k f\right)(x) \right| \right)$$

$$\leq \sum_{k_n=0\vee(\beta_n\wedge\lambda_n)}^{\beta_n} \binom{\beta_n}{k_n} \cdots \sum_{k_1=0\vee(\beta_1-\lambda_1)}^{\beta_1} \binom{\beta_1}{k_1} A_k \, \|f\|_{\alpha+\lambda-\beta+k,k}. \tag{5.4}$$

By (5.4), we obtain that $T_\lambda\colon \mathcal{S}(\mathbb{R}^n,\mathbb{C}) \to \mathcal{S}(\mathbb{R}^n,\mathbb{C})$ is well-defined. Then, by part (v) in Theorem 3.9 and (5.4) again, it proves that $T_\lambda\colon \mathcal{S}(\mathbb{R}^n,\mathbb{C}) \to \mathcal{S}(\mathbb{R}^n,\mathbb{C})$ is $\tau_{\mathcal{S}}$-continuous. This proves part (i) of this lemma. Then, part (ii) of this lemma follows from part (i) and Proposition 4.15 immediately. □

Next, we study the Fourier transform on the Schwartz space $\mathcal{S}$. We first review the definition of Fourier transform on $\mathcal{S}(\mathbb{R}^n,\mathbb{C})$. Let $f \in \mathcal{S}(\mathbb{R}^n,\mathbb{C})$. The Fourier transform of $f$ is a function $\hat{f}\colon \mathbb{R}^n \to \mathbb{C}$ defined by

$$\hat{f}(\xi) = \int_{\mathbb{R}^n} f(t) e^{-i2\pi t\xi} dt, \text{ for any } \xi \in \mathbb{R}^n.$$

Since in this case $f$ is (Lebesgue) integrable over the whole space $\mathbb{R}^n$, the above integral converges to a continuous function $\hat{f}$ at all $\xi \in \mathbb{R}^n$, which is decaying to zero as $\xi \to \infty$. The inverse Fourier transform of $g \in \mathcal{S}(\mathbb{R}^n,\mathbb{C})$ is denoted by $\check{g}\colon \mathbb{R}^n \to \mathbb{C}$ and is defined by

$$\check{g}(t) = \int_{\mathbb{R}^n} g(\xi) e^{i2\pi t\xi} d\xi, \text{ for any } t \in \mathbb{R}^n.$$

We write the Fourier transform by $\mathfrak{F}$ and its inverse by $\mathfrak{F}^{-1}$ such that

$$\mathfrak{F}(f) = \hat{f} \quad \text{and} \quad \mathfrak{F}^{-1}(g) = \check{g}, \text{ for any } f, g \in \mathcal{S}(\mathbb{R}^n,\mathbb{C}).$$

Then, with respect to the topology on the Fréchet space $\mathcal{S}(\mathbb{R}^n,\mathbb{C})$ generated by the countable family $\{\|\cdot\|_{\alpha,\beta}\colon \alpha,\beta \in \mathbb{N}^n\}$ of seminorms, both the Fourier transform $\mathfrak{F}\colon \mathcal{S}(\mathbb{R}^n,\mathbb{C}) \to \mathcal{S}(\mathbb{R}^n,\mathbb{C})$ and its inverse $\mathfrak{F}^{-1}\colon \mathcal{S}(\mathbb{R}^n,\mathbb{C}) \to \mathcal{S}(\mathbb{R}^n,\mathbb{C})$ are $\tau_{\mathcal{S}}$-continuous linear mappings. They have the following properties:

**Theorem 5.5** (see Stein and Weiss [19] or Hunter [7]). *Both $\mathfrak{F}$ and $\mathfrak{F}^{-1}$ are homeomorphisms on the Schwartz space $\mathcal{S}(\mathbb{R}^n,\mathbb{C})$, for the Fréchet space with topology generated by the countable family $\{\|\cdot\|_{\alpha,\beta}\colon \alpha,\beta \in \mathbb{N}^n\}$ of seminorms on $\mathcal{S}(\mathbb{R}^n,\mathbb{C})$.*

By the above properties of Fourier transform $\mathfrak{F}$ and its inverse $\mathfrak{F}^{-1}$, we obtain the following results immediately.

**Proposition 5.6**. *For the Fourier transform $\mathfrak{F}$ and its inverse $\mathfrak{F}^{-1}$ on $\mathcal{S}(\mathbb{R}^n,\mathbb{C})$, we have that*

(i) *$\mathfrak{F}$ is Fréchet differentiable on $\mathcal{S}(\mathbb{R}^n,\mathbb{C})$ such that, for any $\bar{f} \in \mathcal{S}(\mathbb{R}^n,\mathbb{C})$, the Fréchet derivative of $\mathfrak{F}$ at point $\bar{f}$ satisfies that*

$$\nabla\mathfrak{F}(\bar{f})(u) = \mathfrak{F}(u), \text{ for every } u \in \mathcal{S}(\mathbb{R}^n,\mathbb{C}).$$

(ii) *$\mathfrak{F}^{-1}$ is Fréchet differentiable on $\mathcal{S}(\mathbb{R}^n,\mathbb{C})$ such that, for any $\bar{g} \in \mathcal{S}(\mathbb{R}^n,\mathbb{C})$, the Fréchet derivative of $\mathfrak{F}^{-1}$ at point $\bar{g}$ satisfies that*

$$\nabla\mathfrak{F}^{-1}(\bar{g})(v) = \mathfrak{F}^{-1}(v), \text{ for every } v \in \mathcal{S}(\mathbb{R}^n,\mathbb{C}).$$

*Proof*. This proposition is proved by Theorem 5.5 and Proposition 4.15. □

### 5.3. The Fréchet Differentiability of Polynomial Type Operators in Schwartz Space

In this subsection, we consider the Fréchet differentiability of some polynomial type operators, which are not linear in general but considered as special class of mappings on Schwartz space $\mathcal{S}(\mathbb{R}^n, \mathbb{C})$. We will use the definition of Fréchet derivatives (4.12) and (4.7) to precisely calculate the Fréchet derivatives of polynomial functional on $\mathcal{S}(\mathbb{R}^n, \mathbb{C})$. We first prove the following lemma.

**Lemma 5.7**. *Let m be a positive integer. Let* $\alpha, \beta \in \mathbb{N}^n$. *Then, for any* $u \in \mathcal{S}(\mathbb{R}^n, \mathbb{C})$, *we have*

$$\|u^m\|_{\alpha,\beta} \leq 2^{m|\beta|}\left(\max\{\|u\|_{\alpha,\gamma} : \gamma \in \mathbb{N}^n, \theta_n \preccurlyeq \gamma \preccurlyeq \beta\}\right)\left(max\{\|u\|_{\theta_n,\gamma} : \gamma \in \mathbb{N}^n, \theta_n \preccurlyeq \gamma \preccurlyeq \beta\}\right)^{m-1}.$$

*In particular, if* $\alpha = \theta_n$, *then we have*

$$\|u^m\|_{\theta_n,\beta} \leq 2^{m|\beta|}\left(\max\{\|u\|_{\theta_n,\gamma} : \gamma \in \mathbb{N}^n, \theta_n \preccurlyeq \gamma \preccurlyeq \beta\}\right)^m. \tag{5.5}$$

*Proof*. It is clear that when $m = 1$, this lemma is true. Hence, we only prove this lemma for $m \geq 2$. Let $m \geq 2$ and $\beta = (\beta_1, \dots, \beta_n) \in \mathbb{N}^n$ be given (fixed). For any $u \in \mathcal{S}(\mathbb{R}^n, \mathbb{C})$, repeating the product rule given in (5.2), we have that

$$\partial_1^{\beta_1} \dots \partial_n^{\beta_n}(u^m) = \partial_1^{\beta_1} \dots \partial_n^{\beta_n}(uu^{m-1})$$

$$= \sum_{k_{n1}=0}^{\beta_n} \binom{\beta_n}{k_{n1}} \dots \sum_{k_{11}=0}^{\beta_1} \binom{\beta_1}{k_{11}} \left(\partial_1^{\beta_1-k_{11}} \dots\ \partial_n^{\beta_n-k_{n1}} u\right)\left(\partial_1^{k_{11}} \dots \partial_n^{k_{n1}} u^{m-1}\right)$$

$$= \sum_{k_{n1}=0}^{\beta_n} \binom{\beta_n}{k_n} \dots \sum_{k_{11}=0}^{\beta_1} \binom{\beta_1}{k_{11}} \left(\partial_1^{\beta_1-k_{11}} \dots\ \partial_n^{\beta_n-k_{n1}} u\right)\left(\partial_1^{k_{11}} \dots \partial_n^{k_{n1}} uu^{m-2}\right)$$

$$= \sum_{k_{n1}=0}^{\beta_n} \binom{\beta_n}{k_{n1}} \dots \sum_{k_{11}=0}^{\beta_1} \binom{\beta_1}{k_{11}} \left(\partial_1^{\beta_1-k_{11}} \dots\ \partial_n^{\beta_n-k_{n1}} u\right)$$

$$\left[\sum_{k_{n2}=0}^{k_{n1}} \binom{k_{n1}}{k_{n2}} \dots \sum_{k_{12}=0}^{k_{11}} \binom{k_{11}}{k_{12}} \left(\partial_1^{k_{11}-k_{12}} \dots\ \partial_n^{k_{n1}-k_{n2}} u\right)\left(\partial_1^{k_{12}} \dots \partial_n^{k_{n2}} u^{m-2}\right)\right]$$

$$= \sum_{k_{n1}=0}^{\beta_n} \binom{\beta_n}{k_{n1}} \dots \sum_{k_{11}=0}^{\beta_1} \binom{\beta_1}{k_{11}} \left(\partial_1^{\beta_1-k_{11}} \dots\ \partial_n^{\beta_n-k_{n1}} u\right)$$

$$\left\{\sum_{k_{n2}=0}^{k_{n1}} \binom{k_{n1}}{k_{n2}} \dots \sum_{k_{12}=0}^{k_{11}} \binom{k_{11}}{k_{12}} \left(\partial_1^{k_{11}-k_{12}} \dots\ \partial_n^{k_{n1}-k_{n2}} u\right)\right.$$

$$\left.\dots \sum_{k_{nm-1}=0}^{k_{nm-2}} \binom{k_{nm-2}}{k_{nm-1}} \dots \sum_{k_{1m-1}=0}^{k_{1m-2}} \binom{k_{1m-2}}{k_{1m-1}} \left(\partial_1^{k_{1m-2}-k_{1m-1}} \dots\ \partial_n^{k_{nm-2}-k_{nm-1}} u\right)\left(\partial_1^{k_{1m-1}} \dots \partial_n^{k_{nm-1}} u\right)\right\}. \tag{5.6}$$

Here, for every $j = 1, 2, \dots, m-1$, we have $(k_{1j}, k_{2j}, \dots, k_{nj}) \in \mathbb{N}^n$ such that

$$\theta_n \preccurlyeq (k_{1m-1}, k_{2m-1}, \dots, k_{nm-1}) \preccurlyeq (k_{1m-2}, k_{2m-2}, \dots, k_{nm-2}) \preccurlyeq \cdots$$

$$\dots \preccurlyeq (k_{11}, k_{21}, \dots, k_{n1}) \preccurlyeq (\beta_1, \dots, \beta_n) = \beta.$$

By (5.6) and the above order inequality, we estimate that

$$\|u^m\|_{\alpha,\beta} = \sup_{x \in \mathbb{R}^n} \left|x^\alpha \partial_1^{\beta_1} \dots\ \partial_n^{\beta_n} u^m(x)\right|$$

$$= \sup_{x \in \mathbb{R}^n} \left|x^\alpha \left(\sum_{k_{n1}=0}^{\beta_n} \binom{\beta_n}{k_{n1}} \dots \sum_{k_{11}=0}^{\beta_1} \binom{\beta_1}{k_{11}} \left(\partial_1^{\beta_1-k_{11}} \dots\ \partial_n^{\beta_n-k_{n1}} u(x)\right)\right.\right.$$

$$\left\{\sum_{k_{n2}=0}^{k_{n1}} \binom{k_{n1}}{k_{n2}} \dots \sum_{k_{12}=0}^{k_{11}} \binom{k_{11}}{k_{12}} \left(\partial_1^{k_{11}-k_{12}} \dots\ \partial_n^{k_{n1}-k_{n2}} u(x)\right) \dots \sum_{k_{nm-1}=0}^{k_{nm-2}} \binom{k_{nm-2}}{k_{nm-1}}\right.$$

$$\dots\sum_{k_{1m-1}=0}^{k_{1m-2}}\binom{k_{1m-2}}{k_{1m-1}}\left(\partial_1^{k_{1m-2}-k_{1m-1}}\dots\,\partial_n^{k_{nm-2}-k_{nm-1}}u(x)\right)\left(\partial_1^{k_{1m-1}}\dots\partial_n^{k_{nm-1}}u(x)\right)\Big\}\Big)\Big|$$

$$=\sup_{x\in\mathbb{R}^n}\Big|\Big(\sum_{k_{n1}=0}^{\beta_n}\binom{\beta_n}{k_{n1}}\dots\sum_{k_{11}=0}^{\beta_1}\binom{\beta_1}{k_{11}}\left(x^\alpha\partial_1^{\beta_1-k_{11}}\dots\,\partial_n^{\beta_n-k_{n1}}u(x)\right)$$

$$\Big\{\sum_{k_{n2}=0}^{k_{n1}}\binom{k_{n1}}{k_{n2}}\dots\sum_{k_{12}=0}^{k_{11}}\binom{k_{11}}{k_{12}}\left(\partial_1^{k_{11}-k_{12}}\dots\,\partial_n^{k_{n1}-k_{n2}}u(x)\right)\dots\sum_{k_{nm-1}=0}^{k_{nm-2}}\binom{k_{nm-2}}{k_{nm-1}}$$

$$\dots\sum_{k_{1m-1}=0}^{k_{1m-2}}\binom{k_{1m-2}}{k_{1m-1}}\left(\partial_1^{k_{1m-2}-k_{1m-1}}\dots\,\partial_n^{k_{nm-2}-k_{nm-1}}u(x)\right)\left(\partial_1^{k_{1m-1}}\dots\partial_n^{k_{nm-1}}u(x)\right)\Big\}\Big)\Big|$$

$$\leq\sum_{k_{n1}=0}^{\beta_n}\binom{\beta_n}{k_{n1}}\dots\sum_{k_{11}=0}^{\beta_1}\binom{\beta_1}{k_{11}}\sup_{x\in\mathbb{R}^n}\left|x^\alpha\partial_1^{\beta_1-k_{11}}\dots\,\partial_n^{\beta_n-k_{n1}}u(x)\right|$$

$$\Big\{\sum_{k_{n2}=0}^{k_{n1}}\binom{k_{n1}}{k_{n2}}\dots\sum_{k_{12}=0}^{k_{11}}\binom{k_{11}}{k_{12}}\sup_{x\in\mathbb{R}^n}\left|\partial_1^{k_{11}-k_{12}}\dots\,\partial_n^{k_{n1}-k_{n2}}u(x)\right|\dots\sum_{k_{nm-1}=0}^{k_{nm-2}}\binom{k_{nm-2}}{k_{nm-1}}$$

$$\dots\sum_{k_{1m-1}=0}^{k_{1m-2}}\binom{k_{1m-2}}{k_{1m-1}}\sup_{x\in\mathbb{R}^n}\left|\partial_1^{k_{1m-2}-k_{1m-1}}\dots\,\partial_n^{k_{nm-2}-k_{nm-1}}u(x)\right|\sup_{x\in\mathbb{R}^n}\left|\partial_1^{k_{1m-1}}\dots\partial_n^{k_{nm-1}}u(x)\right|\Big\}$$

$$\leq 2^{|\beta|}\max\left\{\sup_{x\in\mathbb{R}^n}\left|x^\alpha\partial_1^{\beta_1-k_{11}}\dots\,\partial_n^{\beta_n-k_{n1}}u(x)\right|:\theta_n\preccurlyeq(k_{11},k_{21},\dots,k_{n1})\preccurlyeq\beta\right\}$$

$$\Big\{2^{|\beta|}\max\left\{\sup_{x\in\mathbb{R}^n}\left|\partial_1^{k_{11}-k_{12}}\dots\,\partial_n^{k_{n1}-k_{n2}}u(x)\right|:\theta_n\preccurlyeq(k_{12},k_{22},\dots,k_{n2})\preccurlyeq\dots\preccurlyeq(k_{11},k_{21},\dots,k_{n1})\right\}$$

$$\dots2^{|\beta|}\max\left\{\sup_{x\in\mathbb{R}^n}\left|\partial_1^{k_{1m-2}-k_{1m-1}}\dots\,\partial_n^{k_{nm-2}-k_{nm-1}}u(x)\right|:\theta_n\preccurlyeq(k_{1m-1},\dots,k_{nm-1})\preccurlyeq(k_{1m-2},\dots,k_{nm-2})\right\}$$

$$\max\left\{\sup_{x\in\mathbb{R}^n}\left|\partial_1^{k_{1m-1}}\dots\partial_n^{k_{nm-1}}u(x)\right|:\theta_n\preccurlyeq(k_{1m-1},\dots,k_{nm-1})\preccurlyeq(k_{1m-2},\dots,k_{nm-2})\right\}\Big\}$$

$$\leq 2^{|\beta|}\left\{\max\{\|u\|_{\alpha,\gamma}:\gamma\in\mathbb{N}^n,\theta_n\preccurlyeq\gamma\preccurlyeq\beta\}\right\}\Big\{2^{|\beta|}\left\{\max\{\|u\|_{\theta_n,\gamma}:\gamma\in\mathbb{N}^n,\theta_n\preccurlyeq\gamma\preccurlyeq\beta\}\right\}$$

$$\dots2^{|\beta|}\left\{\max\{\|u\|_{\theta_n,\gamma}:\gamma\in\mathbb{N}^n,\theta_n\preccurlyeq\gamma\preccurlyeq\beta\}\right\}\left\{\max\{\|u\|_{\theta_n,\gamma}:\gamma\in\mathbb{N}^n,\theta_n\preccurlyeq\gamma\preccurlyeq\beta\}\right\}\Big\}$$

$$=2^{|\beta|}\left\{\max\{\|u\|_{\alpha,\gamma}:\gamma\in\mathbb{N}^n,\theta_n\preccurlyeq\gamma\preccurlyeq\beta\}\right\}2^{(m-1)|\beta|}\left\{\max\{\|u\|_{\theta_n,\gamma}:\gamma\in\mathbb{N}^n,\theta_n\preccurlyeq\gamma\preccurlyeq\beta\}\right\}^{m-1}$$

$$=2^{m|\beta|}\left\{\max\{\|u\|_{\alpha,\gamma}:\gamma\in\mathbb{N}^n,\theta_n\preccurlyeq\gamma\preccurlyeq\beta\}\right\}\left\{\max\{\|u\|_{\theta_n,\gamma}:\gamma\in\mathbb{N}^n,\theta_n\preccurlyeq\gamma\preccurlyeq\beta\}\right\}^{m-1}.$$ □

Before we investigate the Fréchet derivatives of some polynomial type operators on $\mathcal{S}(\mathbb{R}^n,\mathbb{C})$, we first study the power operators, which are the special cases. It is because any polynomial type operator on $\mathcal{S}(\mathbb{R}^n,\mathbb{C})$ is defined by some power operators. Let *m* be a positive integer. Define a power operator $P^m:\mathcal{S}(\mathbb{R}^n,\mathbb{C})\to\mathcal{S}(\mathbb{R}^n,\mathbb{C})$ by

$$P^m(f)=f^m,\text{ for every } f\in\mathcal{S}(\mathbb{R}^n,\mathbb{C}).$$

**Theorem 5.8**. *Let m be a positive integer. Then, $P^m$ is Fréchet differentiable on $\mathcal{S}(\mathbb{R}^n,\mathbb{C})$ such that for any given $\bar{f}\in\mathcal{S}(\mathbb{R}^n,\mathbb{C})$, the Fréchet derivative of $P^m$ at $\bar{f}$ satisfies that*

(i) *If $m = 1$, then*

$$\nabla P^1(\bar{f})(u) = u, \textit{for every } u \in \mathcal{S}(\mathbb{R}^n, \mathbb{C}).$$

(ii) *If* $m > 1$, *then*

(a) *for* $\bar{f} \equiv \theta_{\mathcal{S}}$, *we have*

$$\nabla P^m(\theta_{\mathcal{S}})(u) \equiv \theta_{\mathcal{S}}, \textit{for every } u \in \mathcal{S}(\mathbb{R}^n, \mathbb{C}).$$

(b) *for* $\bar{f} \not\equiv \theta_{\mathcal{S}}$, *we have*

$$\nabla P^m(\bar{f})(u) = m\bar{f}^{m-1}u, \textit{for every } u \in \mathcal{S}(\mathbb{R}^n, \mathbb{C}). \tag{5.7}$$

*More precisely, for any* $u \in \mathcal{S}(\mathbb{R}^n, \mathbb{C})$,

$$\nabla P^m(\bar{f})(u)(x) = m\bar{f}^{m-1}(x)u(x), \textit{for every } x \in \mathbb{R}^n.$$

*Proof*. Proof of (i). When $m = 1$, $P^1$ is the identity mapping that is a $\tau_{\mathcal{S}}$-continuous and linear mapping on $\mathcal{S}(\mathbb{R}^n, \mathbb{C})$. By Proposition 4.15, we have

$$\nabla P^1(\bar{f})(u) = u, \text{ for every } u \in \mathcal{S}(\mathbb{R}^n, \mathbb{C}).$$

Proof of (b) of (ii) of this theorem. Let $\bar{f} \in \mathcal{S}(\mathbb{R}^n, \mathbb{C})$ with $\bar{f} \not\equiv \theta_{\mathcal{S}}$.

Recall that $\preccurlyeq$ is the pointwise partial order on $\mathbb{N}^n$. By $\preccurlyeq$ we define the pointwise partial order $\preccurlyeq^2$ on $\mathbb{N}^n \times \mathbb{N}^n$. For any $(\alpha^{(1)}, \beta^{(1)}), (\alpha^{(2)}, \beta^{(2)}) \in \mathbb{N}^n \times \mathbb{N}^n$, we say that

$$(\alpha^{(1)}, \beta^{(1)}) \preccurlyeq^2 (\alpha^{(2)}, \beta^{(2)}) \quad \text{if and only if} \quad \alpha^{(1)} \preccurlyeq \alpha^{(2)} \ \text{ and } \ \beta^{(1)} \preccurlyeq \beta^{(2)}.$$

Let $J \in \mathcal{F}_{\mathcal{S}}$ and $\varepsilon > 0$ be arbitrarily given, where $0 < \varepsilon < 1$. Since $J$ is finite, there is a positive integer $N$ such that $J$ can be written as

$$J = \left\{ \|\cdot\|_{\alpha^{(j)}, \beta^{(j)}} : (\alpha^{(j)}, \beta^{(j)}) \in \mathbb{N}^n \times \mathbb{N}^n, \text{for } j = 1, 2, \ldots, N \right\}.$$

Here, for each $j = 1, 2, \ldots, N$, we have $\alpha^{(j)} = (\alpha_1^{(j)}, \ldots, \alpha_n^{(j)}) \in \mathbb{N}^n$, $\beta^{(j)} = (\beta_1^{(j)}, \ldots, \beta_n^{(j)}) \in \mathbb{N}^n$. For the given $\bar{f} \in \mathcal{S}(\mathbb{R}^n, \mathbb{C})$ with $\bar{f} \not\equiv \theta_{\mathcal{S}}$, with respect to the finite set $J$, we write

$$[\![\bar{f}]\!] = \max\left\{ \|\bar{f}^{m-i}\|_{\alpha^{(j)}, k^{(j)}} : i = 0, 1, 2, 3, \ldots m, j = 1, 2, \ldots, N \right\}.$$

It is clear that for given $\bar{f} \in \mathcal{S}(\mathbb{R}^n, \mathbb{C})$, $[\![\bar{f}]\!]$ depends on $J \in \mathcal{F}_{\mathcal{S}}$. Meanwhile, since $\bar{f} \in \mathcal{S}(\mathbb{R}^n, \mathbb{C})$ and $\bar{f} \not\equiv \theta_{\mathcal{S}}$, it yields that $0 < [\![\bar{f}]\!] < \infty$. We will prove this proposition by using Proposition 4.4 with respect to the extended $\varepsilon$-$\delta$ language. That is, we need to verify that $P^m$ satisfies (5.7). To this end, for the arbitrarily given $J \in \mathcal{F}_{\mathcal{S}}$ and $\varepsilon > 0$ with $0 < \varepsilon < 1$, we have to find some $\delta > 0$ and $I \in \mathcal{F}_{\mathcal{S}}$ such that $\nabla P^m(\bar{f})$ given by (5.7) will satisfy (DZ) and (DR) in Definition 4.8. For this purpose, with respect to the given $J \in \mathcal{F}_{\mathcal{S}}$ and $\varepsilon > 0$, let $\alpha = (\alpha_1, \ldots, \alpha_n), \beta = (\beta_1, \ldots, \beta_n) \in \mathbb{N}^n$ be defined by

$$(\alpha, \beta) = \vee_{1 \le j \le N} (\alpha^{(j)}, \beta^{(j)}). \tag{5.8}$$

More precisely speaking, $\alpha = (\alpha_1, \ldots, \alpha_n)$ and $\beta = (\beta_1, \ldots, \beta_n)$ satisfy the following equations

$$\alpha_k = \max\{\alpha_k^{(j)} : j = 1, 2, \ldots, N\} \quad \text{and} \quad \beta_k = \max\{\beta_k^{(j)} : j = 1, 2, \ldots, N\}, \text{ for } k = 1, 2, \ldots, n.$$

Recall that for the above defined $\beta = (\beta_1, \dots, \beta_n) \in \mathbb{N}^n$, we have $|\beta| = \beta_1 + \cdots + \beta_n$. Now, with respect to the arbitrarily given $J \in \mathcal{F}_{\mathcal{S}}$ and $\varepsilon > 0$, we take $\delta > 0$ and $I \in \mathcal{F}_{\mathcal{S}}$ as follows

$$\delta = \frac{\varepsilon}{(m-1)m!(\llbracket \bar{f} \rrbracket + 1)2^{(m+1)|\beta|}} \quad \text{and} \quad I = \{\|\cdot\|_{\lambda,\gamma} \in \mathbb{F}_{\mathcal{S}}: (\sigma,\gamma) \preccurlyeq^2 (\alpha,\beta)\}. \tag{5.9}$$

It is clearly to see that $0 < \delta < 1$ and $I \in \mathcal{F}_{\mathcal{S}}$. Notice that $\|\cdot\|_{\alpha^{(j)},\, \beta^{(j)}}$ is a seminorm on $\mathcal{S}(\mathbb{R}^n, \mathbb{C})$. Then, for $u \in \mathcal{S}(\mathbb{R}^n, \mathbb{C})$, if $\max\{\|u\|_{\sigma,\gamma}: \|\cdot\|_{\sigma,\gamma} \in I\} > 0$, then we have

$$\max\left\{\left\|\frac{P^m(\bar{f}+u) - P^m(\bar{f}) - m\bar{f}^{m-1}u}{\max\{\|u\|_{\sigma,\gamma}: \|\cdot\|_{\sigma,\gamma} \in I\}}\right\|_{\alpha^{(j)},\, \beta^{(j)}} : \|\cdot\|_{\alpha^{(j)},\, \beta^{(j)}} \in J\right\} = \frac{\max\{\|P^m(\bar{f}+u) - P^m(\bar{f}) - m\bar{f}^{m-1}u\|_{\alpha^{(j)},\, \beta^{(j)}} : \|\cdot\|_{\alpha^{(j)},\, \beta^{(j)}} \in J\}}{\max\{\|u\|_{\sigma,\gamma}: \|\cdot\|_{\sigma,\gamma} \in I\}}.$$

Then, the property (DR) in Definition 4.8 is equivalent to

$$0 < \max\{\|u\|_{\sigma,\gamma}: \|\cdot\|_{\sigma,\gamma} \in I\} < \delta \implies \frac{\max\{\|P^m(\bar{f}+u) - P^m(\bar{f}) - m\bar{f}^{m-1}u\|_{\alpha^{(j)},\, \beta^{(j)}} : \|\cdot\|_{\alpha^{(j)},\, \beta^{(j)}} \in J\}}{\max\{\|u\|_{\sigma,\gamma}: \|\cdot\|_{\sigma,\gamma} \in I\}} < \varepsilon.$$

Before we precede the proof that $\nabla P^m(\bar{f})$ given by (5.7) satisfies conditions (DZ) and (DR) in Definition 4.8 with respect to the $\tau_{\mathcal{S}}$-open neighborhood $U_{I,\delta}$ of $\theta_{\mathcal{S}}$ constructed in (5.9), we first investigate the following properties of the set $I \in \mathcal{F}_{\mathcal{S}}$.

(a) $\|\cdot\|_{\theta_n,\theta_n} \in I$ and $\|\cdot\|_{\alpha,\beta} \in I$;

(b) $J \subseteq I$;

(c) For any $(\sigma,\gamma) \in \mathbb{N}^n \times \mathbb{N}^n$, we have that

$$(\theta_n, \theta_n) \preccurlyeq^2 (\sigma,\gamma) \preccurlyeq^2 (\alpha,\beta) \quad \Leftrightarrow \quad \|\cdot\|_{\lambda,\gamma} \in I.$$

(d) By the definition (2.1), the $\tau_{\mathcal{S}}$-open neighborhood $U_{I,\delta}$ of $\theta_{\mathcal{S}}$ in $\mathcal{S}(\mathbb{R}^n, \mathbb{C})$ satisfies that

$$U_{I,\delta} = \{u \in \mathcal{S}(\mathbb{R}^n, \mathbb{C}): \max\{\|u\|_{\sigma,\gamma}: \|\cdot\|_{\sigma,\gamma} \in I\} < \delta\}$$

$$= \{u \in \mathcal{S}(\mathbb{R}^n, \mathbb{C}): \max\{\|u\|_{\sigma,\gamma}: (\theta_n, \theta_n) \preccurlyeq^2 (\sigma,\gamma) \preccurlyeq^2 (\alpha,\beta)\} < \delta\}.$$

(e) By the inclusion property (b) of the finite subsets $I$ and $J$, for any $u \in \mathcal{S}(\mathbb{R}^n, \mathbb{C})$, we have

$$\max\{\|u\|_{\alpha^{(j)},\, \beta^{(j)}}: \|\cdot\|_{\alpha^{(j)},\, \beta^{(j)}} \in J\} \le \max\{\|u\|_{\sigma,\gamma}: \|\cdot\|_{\sigma,\gamma} \in I\}. \tag{5.10}$$

Now, we prove that $\nabla P^m(\bar{f})$ given by (5.7) satisfies condition (DZ) in Definition 4.8 with respect to this $\tau_{\mathcal{S}}$-open neighborhood $U_{I,\delta}$ of $\theta_{\mathcal{S}}$ in $\mathcal{S}(\mathbb{R}^n, \mathbb{C})$ constructed in (5.9).

Let $u \in U_{I,\delta}$. Suppose that

$$\max\{\|u\|_{\sigma,\gamma}: \|\cdot\|_{\sigma,\gamma} \in I\} = \max\{\|u\|_{\sigma,\gamma}: (\theta_n, \theta_n) \preccurlyeq^2 (\sigma,\gamma) \preccurlyeq^2 (\alpha,\beta)\} = 0. \tag{5.11}$$

By the above property (c) and (5.11), we have

$$(\theta_n, \theta_n) \preccurlyeq^2 (\alpha,\beta) \text{ and } \|u\|_{\theta_n,\theta_n} \le \max\{\|u\|_{\sigma,\gamma}: (\sigma,\gamma) \preccurlyeq^2 (\alpha,\beta)\} = 0.$$

This implies that

$$\sup_{x \in \mathbb{R}^n} |u(x)| = \|u\|_{\theta_n, \theta_n} = 0.$$

That is, $u \equiv \theta_{\mathcal{S}}$, which implies that

$$P^m(\bar{f} + \theta_{\mathcal{S}}) - P^m(\bar{f}) - m\bar{f}^{m-1}\theta_{\mathcal{S}} \equiv \theta_{\mathcal{S}}.$$

Hence, we obtain that

$$\max\{\|u\|_{\sigma,\gamma} \colon \|\cdot\|_{\sigma,\gamma} \in I\} = 0 \quad \Longrightarrow \quad P^m(\bar{f} + u) - P^m(\bar{f}) - m\bar{f}^{m-1}u = \theta_{\mathcal{S}}$$

$$\Longrightarrow \max\left\{\left\|P^m(\bar{f} + u) - P^m(\bar{f}) - m\bar{f}^{m-1}u\right\|_{\alpha^{(j)},\, \beta^{(j)}} \colon \|\cdot\|_{\alpha^{(j)},\, \beta^{(j)}} \in J\right\} = 0.$$

This proves that $\nabla P^m(\bar{f})(u) = m\bar{f}^{m-1}u$ satisfies condition (DZ) in Definition 4.8 with respect to $\delta$ and $I \in \mathcal{F}_{\mathcal{S}}$ given in (5.9).

Next, we prove that $\nabla P^m(\bar{f})$ given by (5.7) satisfies condition (DR) in Definition 4.8 with respect to the $\tau_{\mathcal{S}}$-open neighborhood $U_{I,\delta}$ of $\theta_{\mathcal{S}}$ in $\mathcal{S}(\mathbb{R}^n, \mathbb{C})$ constructed in (5.9).

Let $u \in U_{I,\delta}$. Suppose $0 < \max\{\|u\|_{\sigma,\gamma} \colon \|\cdot\|_{\sigma,\gamma} \in I\} < \delta$. In the case (b) in (ii) of this theorem, since $\bar{f} \in \mathcal{S}(\mathbb{R}^n, \mathbb{C})$ and $\bar{f} \not\equiv \theta_{\mathcal{S}}$, this yields that $0 < [\![\bar{f}]\!] < \infty$. By (5.10), we have

$$\frac{\max\left\{\|P^m(\bar{f}+u) - P^m(\bar{f}) - m\bar{f}^{m-1}u\|_{\alpha^{(j)},\, \beta^{(j)}} \colon \|\cdot\|_{\alpha^{(j)},\, \beta^{(j)}} \in J\right\}}{\max\{\|u\|_{\sigma,\gamma} \colon \|\cdot\|_{\sigma,\gamma} \in I\}}$$

$$= \frac{\max\left\{\|(\bar{f}+u)^m - \bar{f}^m - m\bar{f}^{m-1}u\|_{\alpha^{(j)},\, \beta^{(j)}} \colon \|\cdot\|_{\alpha^{(j)},\, \beta^{(j)}} \in J\right\}}{\max\{\|u\|_{\sigma,\gamma} \colon \|\cdot\|_{\sigma,\gamma} \in I\}}$$

$$= \frac{\max\left\{\|\binom{m}{2}\bar{f}^{m-2}u^2 + \binom{m}{3}\bar{f}^{m-3}u^3 + \dots + u^m\|_{\alpha^{(j)},\, \beta^{(j)}} \colon \|\cdot\|_{\alpha^{(j)},\, \beta^{(j)}} \in J\right\}}{\max\{\|u\|_{\sigma,\gamma} \colon \|\cdot\|_{\sigma,\gamma} \in I\}}$$

$$\leq \frac{\max\left\{\|\binom{m}{2}\bar{f}^{m-2}u^2\|_{\alpha^{(j)},\, \beta^{(j)}} + \|\binom{m}{3}\bar{f}^{m-3}u^3\|_{\alpha^{(j)},\, \beta^{(j)}} + \dots + \binom{m}{m}\|u^m\|_{\alpha^{(j)},\, \beta^{(j)}} \colon \|\cdot\|_{\alpha^{(j)},\, \beta^{(j)}} \in J\right\}}{\max\{\|u\|_{\sigma,\gamma} \colon \|\cdot\|_{\sigma,\gamma} \in I\}}. \tag{5.12}$$

For any $j = 1, 2, \dots, N$, by (5.3), for $i = 2, 3, \dots, m$, we estimate

$$\left\|\binom{m}{i}\bar{f}^{m-i}u^i\right\|_{\alpha^{(j)},\, \beta^{(j)}} = \binom{m}{i}\left\|u^i\bar{f}^{m-i}\right\|_{\alpha^{(j)},\, \beta^{(j)}}$$

$$\leq \binom{m}{i}\sum_{k_n^{(j)}=0}^{\beta_n^{(j)}}\binom{\beta_n^{(j)}}{k_n^{(j)}} \dots \sum_{k_1^{(j)}=0}^{\beta_1^{(j)}}\binom{\beta_1^{(j)}}{k_1^{(j)}}\left\|u^i\right\|_{\theta_n,\, \beta^{(j)} - k^{(j)}}\left\|\bar{f}^{m-i}\right\|_{\alpha^{(j)},\, k^{(j)}}. \tag{5.13}$$

Here, $k^{(j)} = (k_1^{(j)}, \dots, k_n^{(j)}) \in \mathbb{N}^n$ satisfying $\theta_n \preccurlyeq k^{(j)} \preccurlyeq \beta^{(j)}$, for $j = 1, 2, \dots, N$. By (5.13) and (5.5), for $i = 2, 3, \dots m$, we have

$$\left\|\binom{m}{i}\bar{f}^{m-i}u^i\right\|_{\alpha^{(j)},\, \beta^{(j)}}$$

$$\leq \binom{m}{i}[\![\bar{f}]\!]\sum_{k_n^{(j)}=0}^{\beta_n^{(j)}}\binom{\beta_n^{(j)}}{k_n^{(j)}} \dots \sum_{k_1^{(j)}=0}^{\beta_1^{(j)}}\binom{\beta_1^{(j)}}{k_1^{(j)}}\left\|u^i\right\|_{\theta_n,\, \beta^{(j)} - k^{(j)}}$$

$$\leq \binom{m}{i}[\![\bar{f}]\!]\sum_{k_n^{(j)}=0}^{\beta_n^{(j)}}\binom{\beta_n^{(j)}}{k_n^{(j)}}\dots\sum_{k_1^{(j)}=0}^{\beta_1^{(j)}}\binom{\beta_1^{(j)}}{k_1^{(j)}}2^{i|\beta^{(j)}-k^{(j)}|}\left\{\max\left\{\|u\|_{\theta_n,\gamma}:\gamma\in\mathbb{N}^n,\theta_n\preccurlyeq\gamma\preccurlyeq\ \beta^{(j)}-\ k^{(j)}\right\}\right\}^i$$

$$\leq \binom{m}{i}[\![\bar{f}]\!]\sum_{k_n^{(j)}=0}^{\beta_n^{(j)}}\binom{\beta_n^{(j)}}{k_n^{(j)}}\dots\sum_{k_1^{(j)}=0}^{\beta_1^{(j)}}\binom{\beta_1^{(j)}}{k_1^{(j)}}2^{m|\beta^{(j)}|}\left\{\max\left\{\|u\|_{\theta_n,\gamma}:\gamma\in\mathbb{N}^n,\theta_n\preccurlyeq\gamma\preccurlyeq\ \beta^{(j)}\right\}\right\}^i$$

$$\leq \binom{m}{i}[\![\bar{f}]\!]\sum_{k_n^{(j)}=0}^{\beta_n^{(j)}}\binom{\beta_n^{(j)}}{k_n^{(j)}}\dots\sum_{k_1^{(j)}=0}^{\beta_1^{(j)}}\binom{\beta_1^{(j)}}{k_1^{(j)}}2^{m|\beta|}\left\{\max\left\{\|u\|_{\theta_n,\gamma}:\gamma\in\mathbb{N}^n,\theta_n\preccurlyeq\gamma\preccurlyeq\beta\right\}\right\}^i$$

$$\leq \binom{m}{i}[\![\bar{f}]\!]2^{\beta_n^{(j)}}\dots 2^{\beta_1^{(j)}}2^{m|\beta|}\left\{\max\left\{\|u\|_{\theta_n,\gamma}:\gamma\in\mathbb{N}^n,\theta_n\preccurlyeq\gamma\preccurlyeq\beta\right\}\right\}^i$$

$$= \binom{m}{i}[\![\bar{f}]\!]2^{|\beta^{(j)}|}\,2^{m|\beta|}\left\{\max\left\{\|u\|_{\theta_n,\gamma}:\gamma\in\mathbb{N}^n,\theta_n\preccurlyeq\gamma\preccurlyeq\beta\right\}\right\}^i$$

$$\leq m!\,[\![\bar{f}]\!]2^{(m+1)|\beta|}\left\{\max\left\{\|u\|_{\theta_n,\gamma}:\gamma\in\mathbb{N}^n,\theta_n\preccurlyeq\gamma\preccurlyeq\beta\right\}\right\}^i,\ \text{for } i=2,3,\dots m. \tag{5.14}$$

Substituting the results of (5.14) into (5.12), by the assumption $\delta = \frac{\varepsilon}{(m-1)m!([\![\bar{f}]\!]+1)2^{(m+1)|\beta|}}$ satisfying $\delta < 1$, we obtain that

$$\frac{\max\left\{\|P^m(\bar{f}+u)-P^m(\bar{f})-m\bar{f}^{m-1}u\|_{\alpha^{(j)},\beta^{(j)}}:\|\cdot\|_{\alpha^{(j)},\beta^{(j)}}\in J\right\}}{\max\left\{\|u\|_{\sigma,\gamma}:\|\cdot\|_{\sigma,\gamma}\in I\right\}}$$

$$\leq \frac{m![\![\bar{f}]\!]2^{(m+1)|\beta|}\max\left\{\left\{\max\{\|u\|_{\theta_n,\gamma}:\gamma\in\mathbb{N}^n,\theta_n\preccurlyeq\gamma\preccurlyeq\beta\}\right\}^2+\dots+\left\{\max\{\|u\|_{\theta_n,\gamma}:\gamma\in\mathbb{N}^n,\theta_n\preccurlyeq\gamma\preccurlyeq\beta\}\right\}^m\right\}}{\max\left\{\|u\|_{\sigma,\gamma}:\|\cdot\|_{\sigma,\gamma}\in I\right\}}$$

$$\leq \frac{m![\![\bar{f}]\!]2^{(m+1)|\beta|}\left\{\left\{\max\{\|u\|_{\theta_n,\gamma}:\gamma\in\mathbb{N}^n,\theta_n\preccurlyeq\gamma\preccurlyeq\beta\}\right\}^2+\dots+\left\{\max\{\|u\|_{\theta_n,\gamma}:\gamma\in\mathbb{N}^n,\theta_n\preccurlyeq\gamma\preccurlyeq\beta\}\right\}^m\right\}}{\max\left\{\|u\|_{\sigma,\gamma}:\|\cdot\|_{\sigma,\gamma}\in I\right\}}$$

$$\leq \frac{m![\![\bar{f}]\!]2^{(m+1)|\beta|}\left\{\left\{\max\{\|u\|_{\sigma,\gamma}:\|\cdot\|_{\sigma,\gamma}\in I\}\right\}^2+\dots+\left\{\max\{\|u\|_{\sigma,\gamma}:\|\cdot\|_{\sigma,\gamma}\in I\}\right\}^m\right\}}{\max\left\{\|u\|_{\sigma,\gamma}:\|\cdot\|_{\sigma,\gamma}\in I\right\}}$$

$$= m!\,[\![\bar{f}]\!]2^{(m+1)|\beta|}\left\{\max\left\{\|u\|_{\sigma,\gamma}:\|\cdot\|_{\sigma,\gamma}\in I\right\}+\dots+\max\left\{\|u\|_{\sigma,\gamma}:\|\cdot\|_{\sigma,\gamma}\in I\right\}^{m-2}\right\}$$

$$= m!\,[\![\bar{f}]\!]2^{(m+1)|\beta|}\left\{\max\left\{\|u\|_{\sigma,\gamma}:\|\cdot\|_{\sigma,\gamma}\in I\right\}+\dots+\max\left\{\|u\|_{\sigma,\gamma}:\|\cdot\|_{\sigma,\gamma}\in I\right\}^{m-2}\right\}$$

$$< m!\,[\![\bar{f}]\!]2^{(m+1)|\beta|}\{\delta+\dots+\delta^{m-2}\}$$

$$< (m-1)m!\,[\![\bar{f}]\!]2^{(m+1)|\beta|}\delta$$

$$= (m-1)m!\,[\![\bar{f}]\!]2^{(m+1)|\beta|}\frac{\varepsilon}{(m-1)m!([\![\bar{f}]\!]+1)2^{(m+1)|\beta|}}$$

$$= \varepsilon.$$

This proves that $\nabla P^m(\bar{f})(u) = m\bar{f}^{m-1}u$ satisfies the condition (DR) in Definition 4.8 with respect to the constructed $\tau_{\mathcal{S}}$-open neighborhood $U_{I,\delta}$ of $\theta_{\mathcal{S}}$ with $I \in \mathcal{F}_{\mathcal{S}}$ and $\delta > 0$ given in (5.9) and with the hypothesis that $\bar{f} \not\equiv \theta_{\mathcal{S}}$. This proved part (b) in (ii) of this theorem.

Proof of (a) in (ii). $\bar{f} \equiv \theta_{\mathcal{S}}$. In this case, for this arbitrarily given $\tau_{\mathcal{S}}$-open neighborhood $V_{J,\varepsilon}$ of $\theta_{\mathcal{S}}$, similarly to (5.8) and (5.9), let $\alpha = (\alpha_1, \dots, \alpha_n), \beta = (\beta_1, \dots, \beta_n) \in \mathbb{N}^n$ be defined by

$$(\alpha, \beta) = \vee_{1 \le j \le N} (\alpha^{(j)},\ \beta^{(j)}). \tag{5.15}$$

Then, we define

$$\lambda = \left(\frac{\varepsilon}{2^{m|\beta|}}\right)^{\frac{1}{m-1}} \quad \text{and} \quad I = \left\{\|\cdot\|_{\lambda,\gamma} \in \mathbb{F}_{\mathcal{S}} : (\sigma, \gamma) \preccurlyeq^2 (\alpha, \beta)\right\}. \tag{5.16}$$

We have an $\tau_{\mathcal{S}}$-open neighborhood $U_{I,\lambda}$ of $\theta_{\mathcal{S}}$ with $I \in \mathcal{F}_{\mathcal{S}}$ defined in (5.15) and (5.16). Now, we prove that $\nabla P^m(\theta_{\mathcal{S}})$ given by part (a) in (ii) satisfies condition (DZ) in Definition 4.8 with respect to the $\tau_{\mathcal{S}}$-open neighborhood $U_{I,\lambda}$ of $\theta_{\mathcal{S}}$ in $\mathcal{S}(\mathbb{R}^n, \mathbb{C})$ constructed in (5.16).

Let $u \in U_{I,\lambda}$. Suppose that $\max\left\{\|u\|_{\sigma,\gamma} : \|\cdot\|_{\sigma,\gamma} \in I\right\} = 0$. This is rewritten as

$$\max \{\|u\|_{\sigma,\gamma} : (\theta_n, \theta_n) \preccurlyeq^2 (\sigma, \gamma) \preccurlyeq^2 (\alpha, \beta)\} = 0.$$

By (5.15), we have that $(\theta_n, \theta_n) \preccurlyeq^2 (\alpha, \beta)$. Then, the above equation implies that

$$\|u\|_{\theta_n,\theta_n} \le \max \{\|u\|_{\sigma,\gamma} : (\sigma, \gamma) \preccurlyeq^2 (\alpha, \beta)\} = 0.$$

This implies that $\sup_{x \in \mathbb{R}^n} |u(x)| = \|u\|_{\theta_n,\theta_n} = 0$. We must have $u \equiv \theta_{\mathcal{S}}$. For this given $u \in U_{I,\lambda}$, we obtain

$$P^m(\theta_{\mathcal{S}} + u) - P^m(\theta_{\mathcal{S}}) - \theta_{\mathcal{S}} = P^m(\theta_{\mathcal{S}} + \theta_{\mathcal{S}}) - P^m(\theta_{\mathcal{S}}) - \theta_{\mathcal{S}} \equiv \theta_{\mathcal{S}}.$$

Hence, we have that

$$\max\left\{\|u\|_{\sigma,\gamma} : \|\cdot\|_{\sigma,\gamma} \in I\right\} = 0$$

$$\Longrightarrow \quad P^m\left(\bar{f} + u\right) - P^m\left(\bar{f}\right) - m\bar{f}^{m-1}u = \theta_{\mathcal{S}}$$

$$\Longrightarrow \quad \max\left\{\|P^m(\theta_{\mathcal{S}} + u) - P^m(\theta_{\mathcal{S}}) - \theta_{\mathcal{S}}\|_{\alpha^{(j)},\, \beta^{(j)}} : \|\cdot\|_{\alpha^{(j)},\, \beta^{(j)}} \in J\right\} = 0.$$

This proves that $\nabla P^m(\theta_{\mathcal{S}})(u) \equiv \theta_{\mathcal{S}}$ satisfies condition (DZ) in Definition 4.8 with respect to the $\tau_{\mathcal{S}}$-open neighborhood $U_{I,\lambda}$ of $\theta_{\mathcal{S}}$ in $\mathcal{S}(\mathbb{R}^n, \mathbb{C})$ constructed in (5.15) and (5.16). Next, we prove that $\nabla P^m(\theta_{\mathcal{S}})$ given by (5.7) satisfies condition (DR) in Definition 4.8 with respect to the $\tau_{\mathcal{S}}$-open neighborhood $U_{I,\lambda}$ of $\theta_{\mathcal{S}}$ in $\mathcal{S}(\mathbb{R}^n, \mathbb{C})$ constructed in (5.15) and (5.16).

Let $u \in U_{I,\lambda}$. Suppose $0 < \max\left\{\|u\|_{\sigma,\gamma} : \|\cdot\|_{\sigma,\gamma} \in I\right\} < \lambda$. By Lemma 5.7 and by $\lambda$ in (5.16), we have

$$\frac{\max\left\{\|P^m(\theta_{\mathcal{S}}+u) - P^m(\theta_{\mathcal{S}}) - \theta_{\mathcal{S}}\|_{\alpha^{(j)},\, \beta^{(j)}} : \|\cdot\|_{\alpha^{(j)},\, \beta^{(j)}} \in J\right\}}{\max\{\|u\|_{\sigma,\gamma} : \|\cdot\|_{\sigma,\gamma} \in I\}}$$

$$= \frac{\max\left\{\|u^m\|_{\alpha^{(j)},\, \beta^{(j)}} : \|\cdot\|_{\alpha^{(j)},\, \beta^{(j)}} \in J\right\}}{\max\{\|u\|_{\sigma,\gamma} : \|\cdot\|_{\sigma,\gamma} \in I\}}$$

$$\le \frac{\mathrm{Max}\left\{2^{m\left|\beta^{(j)}\right|}\left\{\max\left\{\|u\|_{\alpha^{(j)},\gamma} : \gamma \in \mathbb{N}^n, \theta_n \preccurlyeq \gamma \preccurlyeq \beta^{(j)}\right\}\right\}\left\{\max\{\|u\|_{\theta_n,\gamma} : \gamma \in \mathbb{N}^n, \theta_n \preccurlyeq \gamma \preccurlyeq \beta^{(j)}\}\right\}^{m-1} : \|\cdot\|_{\alpha^{(j)},\, \beta^{(j)}} \in J\right\}}{\max\{\|u\|_{\sigma,\gamma} : \|\cdot\|_{\sigma,\gamma} \in I\}}$$

$$\leq \frac{\mathrm{Max}\left\{2^{m\left|\beta^{(j)}\right|}\left\{\max\left\{\|u\|_{\alpha^{(j)},\gamma}:\gamma\in\mathbb{N}^n,\theta_n\preccurlyeq\gamma\preccurlyeq\beta\right\}\right\}\left\{\max\left\{\|u\|_{\theta_n,\gamma}:\gamma\in\mathbb{N}^n,\theta_n\preccurlyeq\gamma\preccurlyeq\beta\right\}\right\}^{m-1}:\|\cdot\|_{\alpha^{(j)},\beta^{(j)}}\in J\right\}}{\max\left\{\|u\|_{\sigma,\gamma}:\|\cdot\|_{\sigma,\gamma}\in I\right\}}$$

$$\leq \frac{\mathrm{Max}\left\{2^{m\left|\beta^{(j)}\right|}\left\{\max\left\{\|u\|_{\sigma,\gamma}:\|\cdot\|_{\sigma,\gamma}\in I\right\}\right\}\left\{\max\left\{\|u\|_{\sigma,\gamma}:\|\cdot\|_{\sigma,\gamma}\in I\right\}\right\}^{m-1}:\|\cdot\|_{\alpha^{(j)},\beta^{(j)}}\in J\right\}}{\max\left\{\|u\|_{\sigma,\gamma}:\|\cdot\|_{\sigma,\gamma}\in I\right\}}$$

$$\leq \frac{2^{m|\beta|}\left\{\max\left\{\|u\|_{\sigma,\gamma}:\|\cdot\|_{\sigma,\gamma}\in I\right\}\right\}^{m}}{\max\left\{\|u\|_{\sigma,\gamma}:\|\cdot\|_{\sigma,\gamma}\in I\right\}}$$

$$\leq 2^{m|\beta|}\left\{\max\left\{\|u\|_{\sigma,\gamma}:\|\cdot\|_{\sigma,\gamma}\in I\right\}\right\}^{m-1}$$

$$< \varepsilon.$$

This proves that for arbitrarily given $\tau_{\mathcal{S}}$-open neighborhood $V_{J,\varepsilon}$ of $\theta_{\mathcal{S}}$, $\nabla P^m(\theta_{\mathcal{S}})(u) \equiv \theta_{\mathcal{S}}$ satisfies the conditions (DZ) and (DR) in Definition 4.8 with respect to the constructed $U_{I,\lambda}$ of $\theta_{\mathcal{S}}$ with $I \in \mathcal{F}_{\mathcal{S}}$ and $\lambda > 0$ defined in (5.16), which proves (a) in (ii) of this theorem. The proof is completed. □

By the linearity of Fréchet derivatives (Theorem 4.14), we have the following results for polynomial type operators in $\mathcal{S}(\mathbb{R}^n, \mathbb{C})$. Let *m* be a positive integer. Let $a_m$, $a_{m-1}$, … $a_1$ be complex numbers. Define a polynomial type operator $P_m$: $\mathcal{S}(\mathbb{R}^n, \mathbb{C}) \to \mathcal{S}(\mathbb{R}^n, \mathbb{C})$ by

$$P_m(f) = a_m f^m + a_{m-1} f^{m-1} + \ldots + a_1 f, \text{ for any } f \in \mathcal{S}(\mathbb{R}^n, \mathbb{C}). \tag{5.17}$$

**Theorem 5.9**. *Let m be a positive integer and let $P_m$ be a polynomial type operator on $\mathcal{S}(\mathbb{R}^n, \mathbb{C})$ defined by* (5.17). *Then, $P_m$ is Fréchet differentiable on $\mathcal{S}(\mathbb{R}^n, \mathbb{C})$ such that for any given $\bar{f} \in \mathcal{S}(\mathbb{R}^n, \mathbb{C})$, the Fréchet derivative of $P_m$ at $\bar{f}$ satisfies that*

(i) *If m* = 1, *then*

$$\nabla P_m(\bar{f})(u) = a_1 u, \textit{ for every } u \in \mathcal{S}(\mathbb{R}^n, \mathbb{C}).$$

(ii) *If m* > 1, *then*

(a) *for $\bar{f} \equiv \theta_{\mathcal{S}}$, we have*

$$\nabla P_m(\theta_{\mathcal{S}})(u) \equiv \theta_{\mathcal{S}}, \textit{ for every } u \in \mathcal{S}(\mathbb{R}^n, \mathbb{C}).$$

(b) *for $\bar{f} \not\equiv \theta_{\mathcal{S}}$, we have*

$$\nabla P_m(\bar{f})(u) = \sum_{i=1}^{m} i a_i \bar{f}^{i-1} u, \textit{ for every } u \in \mathcal{S}(\mathbb{R}^n, \mathbb{C}).$$

*More precisely, for any $u \in \mathcal{S}(\mathbb{R}^n, \mathbb{C})$,*

$$\nabla P_m(\bar{f})(u)(x) = \sum_{i=1}^{m} i a_i \bar{f}^{i-1}(x) u(x), \textit{ for every } x \in \mathbb{R}^n.$$

*Proof*. This theorem follows from Theorem 4.14 and (5.7) in Theorem 5.8 immediately. □

### 5.4. The Gâteaux Differentiability of Polynomial Type Operators in Schwartz Space

Recall that Theorem 4.18 provides the connection between Gâteaux directional differentiability and Fréchet differentiability of mappings in Hausdorff topological vector spaces. It shows that with two conditions (SC) and (sb), the Fréchet differentiability implies the Gâteaux directional differentiability.

In Theorem 5.9, the Fréchet differentiability of a polynomial type operator $P_m$ on $\mathcal{S}(\mathbb{R}^n, \mathbb{C})$ is proved; and for any given $\bar{f} \in \mathcal{S}(\mathbb{R}^n, \mathbb{C})$, the Fréchet derivative $\nabla P_m(\bar{f})$ of $P_m$ at $\bar{f}$ was precisely calculated. Now, we consider the Gâteaux differentiability of this polynomial type operators $P_m$ on $\mathcal{S}(\mathbb{R}^n, \mathbb{C})$. For this purpose, notice that the space $\mathcal{S}(\mathbb{R}^n, \mathbb{C})$ is Hausdorff and seminorm constructed that satisfies the condition (SC) in Theorem 4.18. But, for a given $v \in \mathcal{S}(\mathbb{R}^n, \mathbb{C})$, the (sb) boundedness condition given in (4.38) in the space $\mathcal{S}(\mathbb{R}^n, \mathbb{C})$ is difficult to check. Hence, we cannot apply Theorem 4.18 and Theorem 5.9 to prove the Gâteaux differentiability of polynomial type operators $P_m$ on $\mathcal{S}(\mathbb{R}^n, \mathbb{C})$. In this subsection, we will directly prove the Gâteaux differentiability of polynomial type operators $P_m$ in Schwartz space without using its Fréchet differentiability. Meanwhile, we will find that the Gâteaux derivatives of polynomial type operators $P_m$ coincide with their Fréchet derivatives. We start with the power operator $P^m$.

**Proposition 5.10**. *Let m be a positive integer. Let $\bar{f} \in \mathcal{S}(\mathbb{R}^n, \mathbb{C})$. Then, $P^m$ is Gâteaux differentiable at $\bar{f}$ and the Gâteaux derivative of $P^m$ at $\bar{f}$ satisfies that, for every $u \in \mathcal{S}(\mathbb{R}^n, \mathbb{C}) \backslash \{\theta_{\mathcal{S}}\}$*

(i) *If m = 1, then we have*

$$P'(\bar{f})(u) = \nabla P(\bar{f})(u) = u;$$

(ii) *If m > 1, then we have*

$$(P^m)'(\bar{f})(u) = \nabla P^m(\bar{f})(u) = m\bar{f}^{m-1}u. \tag{5.18}$$

*Proof.* It is clear for $m = 1$. Hence, we prove (5.18) with $m > 1$. Let $\bar{f} \in \mathcal{S}(\mathbb{R}^n, \mathbb{C})$, $u \in \mathcal{S}(\mathbb{R}^n, \mathbb{C}) \backslash \{\theta_{\mathcal{S}}\}$. It is also clear that (5.18) holds for $\bar{f} \equiv \theta_{\mathcal{S}}$. So, we prove (5.18) for $\bar{f} \not\equiv \theta_{\mathcal{S}}$.

Let $J \in \mathcal{F}_{\mathcal{S}}$ and $0 < \varepsilon < 1$ be arbitrarily given. By (5.10) and (5.12) and by the proof of Theorem 5.8, for any real $t$ with $0 < |t| < 1$, we have

$$\max\left\{\left\|\frac{P^m(\bar{f}+tu)-P^m(\bar{f})-tm\bar{f}^{m-1}u}{t}\right\|_{\alpha^{(j)},\beta^{(j)}} : \|\cdot\|_{\alpha^{(j)},\beta^{(j)}} \in J\right\}$$

$$= \frac{\max\left\{\left\|P^m(\bar{f}+tu)-P^m(\bar{f})-tm\bar{f}^{m-1}u\right\|_{\alpha^{(j)},\beta^{(j)}} : \|\cdot\|_{\alpha^{(j)},\beta^{(j)}} \in J\right\}}{|t|}$$

$$= \frac{\max\left\{\left\|(\bar{f}+tu)^m-\bar{f}^m-tm\bar{f}^{m-1}u\right\|_{\alpha^{(j)},\beta^{(j)}} : \|\cdot\|_{\alpha^{(j)},\beta^{(j)}} \in J\right\}}{|t|}$$

$$= \frac{\max\left\{\left\|\binom{m}{2}\bar{f}^{m-2}u^2+\binom{m}{3}\bar{f}^{m-3}u^3+\cdots+u^m\right\|_{\alpha^{(j)},\beta^{(j)}} : \|\cdot\|_{\alpha^{(j)},\beta^{(j)}} \in J\right\}}{|t|}$$

$$\leq \frac{t^2\max\left\{\left\|\binom{m}{2}\bar{f}^{m-2}u^2\right\|_{\alpha^{(j)},\beta^{(j)}}+\left\|\binom{m}{3}\bar{f}^{m-3}u^3\right\|_{\alpha^{(j)},\beta^{(j)}}+\cdots+\binom{m}{m}\|u^m\|_{\alpha^{(j)},\beta^{(j)}} : \|\cdot\|_{\alpha^{(j)},\beta^{(j)}} \in J\right\}}{|t|}$$

$$= |t|\max\left\{\left\|\binom{m}{2}\bar{f}^{m-2}u^2\right\|_{\alpha^{(j)},\beta^{(j)}} + \left\|\binom{m}{3}\bar{f}^{m-3}u^3\right\|_{\alpha^{(j)},\beta^{(j)}} + \cdots + \binom{m}{m}\|u^m\|_{\alpha^{(j)},\beta^{(j)}} : \|\cdot\|_{\alpha^{(j)},\beta^{(j)}} \in J\right\}$$

$$\to 0, \text{ as } t \to 0.$$

Here, for given $\bar{f} \in \mathcal{S}(\mathbb{R}^n, \mathbb{C})$ and for fixed $u \in \mathcal{S}(\mathbb{R}^n, \mathbb{C}) \backslash \{\theta_{\mathcal{S}}\}$, we have that

$$\max\left\{\left\|\binom{m}{2}\bar{f}^{m-2}u^2\right\|_{\alpha^{(j)},\ \beta^{(j)}} + \left\|\binom{m}{3}\bar{f}^{m-3}u^3\right\|_{\alpha^{(j)},\ \beta^{(j)}} + \cdots + \binom{m}{m}\|u^m\|_{\alpha^{(j)},\ \beta^{(j)}} : \|\cdot\|_{\alpha^{(j)},\ \beta^{(j)}} \in J\right\} < \infty.$$

This completes the proof of this proposition. □

Next, we consider the Gâteaux differentiability of polynomial type operator $P_m$ on $\mathcal{S}(\mathbb{R}^n, \mathbb{C})$.

**Corollary 5.11**. *Let m be a positive integer and let $P_m$ be a polynomial type operator on $\mathcal{S}(\mathbb{R}^n, \mathbb{C})$ defined by* (5.17). *Then, $P_m$ is Gâteaux differentiable at each point $\bar{f} \in \mathcal{S}(\mathbb{R}^n, \mathbb{C})$ such that*

$$P_m'(\bar{f})(u) = \nabla P_m(\bar{f})(u) = \sum_{i=1}^{m} ia_i\bar{f}^{i-1}u, \textit{ for every } u \in \mathcal{S}(\mathbb{R}^n, \mathbb{C})\backslash\{\theta_{\mathcal{S}}\}.$$

*Proof*. This theorem follows from Proposition 5.10 and the linearity of Gâteaux differentiability immediately. □

## 6. An Example of Not Seminorm-Constructed Topological Vector Space $(\sigma_\rho, \tau_\rho)$

In this section, we consider a topological vector space, denoted by $\sigma_\rho$ with $0 < \rho < 1$, which is not locally convex. This space $\sigma_\rho$ has been used as an example of not seminorm-constructed topological vector spaces. In contrast to the Schwartz space $\mathcal{S}(\mathbb{R}^n, \mathbb{C})$ studied in the previous section, the space $\sigma_\rho$ does not have the property of local convexity. Not surprisingly, the continuity and Fréchet differentiability of mappings in $\sigma_\rho$ are more complicated than that in $\mathcal{S}(\mathbb{R}^n, \mathbb{C})$.

Recall that, in this paper, let $\mathbb{N}$ denotes the set of nonnegative integers. We review the definition of $\sigma_\rho$. Let $\rho$ be a positive number with $0 < \rho < 1$. The space $\sigma_\rho$ consists of sequences of real numbers with finite sum of $\rho$-th power of absolute values of its terms. More precisely, for any sequence of real numbers $\{t_n\}_{n=1}^{\infty}$, we define

$$\{t_n\}_{n=1}^{\infty} \in \sigma_\rho \quad \Leftrightarrow \quad \sum_{n=1}^{\infty}|t_n|^\rho < \infty. \tag{6.1}$$

The origin of $\sigma_\rho$ is denoted by $\theta_{\sigma_\rho}$. Notice that, in (6.1) the space is denoted by $\sigma_\rho$ instead of $l_\rho$, which is to distinguish the ordinary Banach space $l_p$, for $1 \leq p \leq \infty$. By definition (6.1), a translation invariant metric $d_\rho$ is defined on $\sigma_\rho$ by

$$d_\rho(x, y) = \sum_{n=1}^{\infty}|t_n - s_n|^\rho, \text{ for } x = \{t_n\}_{n=1}^{\infty},\ y = \{s_n\}_{n=1}^{\infty} \in \sigma_\rho. \tag{6.2}$$

In this section, we define a topology on $\sigma_\rho$ that is induced by a family of *F*-seminorms on $\sigma_\rho$ as follows. For each $k = 1, 2, \ldots$, we define a *F*-seminorm $\|\cdot\|_{\rho,k}$ on $\sigma_\rho$ by

$$\|x\|_{\rho,k} = |t_k|^\rho, \text{ for any } x = \{t_n\}_{n=1}^{\infty} \in \sigma_\rho. \tag{6.3}$$

Let $\mathbb{F}_{\sigma_\rho} = \left\{\|\cdot\|_{\rho,k} : k = 1, 2, \ldots\right\}$, which is a countable family of *F*-seminorms on $\sigma_\rho$. The collection of nonempty finite subsets $\mathcal{F}_{\sigma_\rho}$ of $\mathbb{F}_{\sigma_\rho}$ is written as

$$\mathcal{F}_{\sigma_\rho} = \left\{\{\|\cdot\|_{\rho,i} \in \mathbb{F}_{\sigma_\rho} : i \in A\} : A \text{ is a nonempty finite subset of } \mathbb{N}\right\}.$$

(6.4)

$(\sigma_\rho, \tau_\rho)$ is a topological vector space equipped with the topology $\tau_\rho$ induced by the countable family $\mathbb{F}_{\sigma_\rho} = \{\|\cdot\|_{\rho,k} : \text{for } k = 1, 2, \ldots .\}$ of *F*-seminorms on $\sigma_\rho$. $(\sigma_\rho, \tau_\rho)$ is not seminorm-constructed. Recall

that, for any given $I \in \mathcal{F}_{\sigma_\rho}$ and $\delta > 0$, we have an $\tau_\rho$-open neighborhood $U_{I,\delta}$ in $\sigma_\rho$ around the origin $\theta_{\sigma_\rho}$ given by

$$U_{I,\delta} = \left\{x \in \sigma_\rho : \max\{\|x\|_{\rho,i} : \|\cdot\|_{\rho,i} \in I\} < \delta\right\}.$$

More precisely speaking, $U_{I,\delta}$ is rewritten as

$$U_{I,\delta} = \left\{x = \{t_n\}_{n=1}^{\infty} \in \sigma_\rho : \max\{|t_i|^\rho : \|\cdot\|_{\rho,i} \in I\} < \delta\right\}.$$

By (2.1) and (2.3), the following collection

$$U(\mathbb{F}_{\sigma_\rho}) = \left\{U_{I,\delta} : I \in \mathcal{F}_{\sigma_\rho}, \delta > 0\right\}, \tag{6.5}$$

forms an $\tau_\rho$-open neighborhood basis of $\sigma_\rho$ around $\theta_{\sigma_\rho}$. More generally, for $y = \{s_n\}_{n=1}^{\infty} \in \sigma_\rho$, $I \in \mathcal{F}_{\sigma_\rho}$ and $\delta > 0$, we write

$$\begin{aligned} U_{I,\delta}(y) &= \left\{x \in \sigma_\rho : \max\{\|x - y\|_{\rho,i} : \|\cdot\|_{\rho,i} \in I\} < \delta\right\} \\ &= \left\{x = \{t_n\}_{n=1}^{\infty} \in \sigma_\rho : \max\{|t_i - s_i|^\rho : \|\cdot\|_{\rho,i} \in I\} < \delta\right\}. \end{aligned}$$

Then, $U_{I,\delta}(y)$ is an $\tau_\rho$-open neighborhood of $\sigma_\rho$ around $y$. By (2.1) and (2.3), the following collection

$$U(\mathbb{F}_{\sigma_\rho})(y) = \left\{U_{I,\delta}(y) : I \in \mathcal{F}_{\sigma_\rho}, \delta > 0\right\}, \tag{6.6}$$

forms an $\tau_\rho$-open neighborhood basis of $\sigma_\rho$ around $y$. The space $\sigma_\rho$ has the following properties.

**Lemma 6.1**. *Let $\rho$ and $\gamma$ be positive numbers. If $0 < \rho < \gamma < 1$, then $\sigma_\rho \subseteq \sigma_\gamma$ and*

*$\sigma_\rho$ is a proper subspace of $\sigma_\gamma$, for $0 < \rho < \gamma < 1$.*

*Proof*. Let $x = \{t_n\}_{n=1}^{\infty} \in \sigma_\rho$. Then, by $\sum_{n=1}^{\infty} |t_n|^\rho < \infty$, there is a positive integer $N$ such that $|t_n|^\rho < 1$, for all $n \geq N$. Since $0 < \rho < \gamma < 1$, it yields that

$$|t_n|^\gamma \leq |t_n|^\rho < 1, \text{ for all } n \geq N.$$

This implies that $x = \{t_n\}_{n=1}^{\infty} \in \sigma_\gamma$. □

The following lemma is a special case of Lemma 3.2. In order to deeply understand the construction of $\sigma_\rho$, we give a direct proof of this lemma.

**Lemma 6.2**. *Let a be a real number. Define the multiplication operator $M_a$ on $\sigma_\rho$ by*

$$M_a x = \{a t_n\}_{n=1}^{\infty}, \text{ for any } x = \{t_n\}_{n=1}^{\infty} \in \sigma_\rho.$$

*Then, we have*

(i) *$M_a : \sigma_\rho \to \sigma_\rho$ is a continuous and linear operator on $\sigma_\rho$.*

(ii) *$M_a$ is Fréchet differentiable on $\sigma_\rho$ such that for any $\bar{x} \in \sigma_\rho$, we have*

$$\nabla(M_a)(\bar{x})(x) = ax = \{at_n\}_{n=1}^{\infty}, \textit{for any } x = \{t_n\}_{n=1}^{\infty} \in \sigma_\rho.$$

*Proof*. Proof of (i). The linearity of $M_a$is clear. We only prove its continuity at point $\theta_{\sigma_\rho}$. The continuity of $M_a$ is clear for $a = 0$. So, we only need to prove the continuity of $M_a$ at point $\theta_{\sigma_\rho}$ for $a \neq 0$. Let $V_{J,\varepsilon}$ be any $\tau_\rho$-open neighborhood of $\theta_{\sigma_\rho}$ in $\sigma_\rho$. Since $J \in \mathcal{F}_{\sigma_\rho}$, ($J$ is a nonempty finite subset in $\mathbb{F}_{\sigma_\rho}$), there is a positive integer $M$ such that

$$J \subseteq \{\|\cdot\|_{\rho,k} : k = 1, 2, \dots, M\}. \tag{6.7}$$

We take an $\tau_\rho$-open neighborhood $U_{I,\delta}$ of $\theta_{\sigma_\rho}$ in $\sigma_\rho$ with $I = \{\|\cdot\|_{\rho,k} : k = 1, 2, \dots, M\}$ and $\delta = \frac{1}{|a|^\rho}\varepsilon$. This induces that, for any $x = \{t_n\}_{n=1}^{\infty} \in U_{I,\delta}$, it satisfies

$$\max\{\|x\|_{\rho,k} : k = 1, 2, \dots, M\} = \max\{\|x\|_{\rho,k} : \|\cdot\|_{\rho,k} \in I\} < \delta.$$

Then, by (6.7), we have

$$\begin{aligned}
&\max\{\|M_a x\|_{\rho,j} : \|\cdot\|_{\rho,j} \in J\} \\
&\leq \max\{\|M_a x\|_{\rho,k} : \|\cdot\|_{\rho,k} \in I\} \\
&= |a|^\rho \max\{\|x\|_{\rho,k} : \|\cdot\|_{\rho,k} \in I\} \\
&< |a|^\rho \delta = \varepsilon.
\end{aligned}$$

By the extended $\varepsilon$-$\delta$ language of the continuity given in Definition 3.3, part (i) of this lemma is proved. Then, by (i) of this lemma and (4.30) in Lemma 4.12, part (ii) of this lemma is proved immediately. □

Let $m$ be a positive integer. We define the $m$th power operator $Q^m : \sigma_\rho \to \sigma_\rho$ as follows.

$$Q^m(x) = \{t_n^m\}_{n=1}^{\infty}, \text{ for any } x = \{t_n\}_{n=1}^{\infty} \in \sigma_\rho. \tag{6.8}$$

The following lemma shows the continuity of the $m^{th}$ power mapping $Q^m$ on $\sigma_\rho$. As a matter of fact, in Theorem 6.4 below, we will prove the Fréchet differentiability of $Q^m$ on $\sigma_\rho$, which implies the continuity of $Q^m$, by part (i) of Theorem 4.16. However, in order to practice the definition of the continuity in general topological vector spaces, in the following lemma, we give a direct proof of the continuity of $Q^m$. We define a mapping $\mathcal{P} : \sigma_\rho \to \mathbb{R}_+$ as follows

$$\mathcal{P}(x) = \sup\{|t_n| : n = 1, 2, \dots\}, \text{ for any } x = \{t_n\}_{n=1}^{\infty} \in \sigma_\rho.$$

It is clear that $\mathcal{P}$ is well-defined and it satisfies that

$$0 \leq \mathcal{P}(x) < \infty, \text{ for any } x = \{t_n\}_{n=1}^{\infty} \in \sigma_\rho.$$

Furthermore, by Definition (6.1), for any $x = \{t_n\}_{n=1}^{\infty} \in \sigma_\rho$, we have that,

$$\mathcal{P}(x) = \sup\{|t_n| : n = 1, 2, \dots\} = \max\{|t_n| : n = 1, 2, \dots\}.$$

**Lemma 6.3**. *Let m be a positive integer. Then, the* $m^{th}$ *power operator* $Q^m$ *is continuous on* $\sigma_\rho$.

*Proof*. It is clear that the lemma holds for $m = 1$. So, we suppose $m > 1$. We first prove that the $m$th power operator $Q^m$ on $\sigma_\rho$ is well-defined. To this end, for $x = \{t_n\}_{n=1}^{\infty} \in \sigma_\rho$ satisfying $\sum_{n=1}^{\infty}|t_n|^\rho < \infty$. We calculate

$$\sum_{n=1}^{\infty}|t_n^m|^\rho = \sum_{n=1}^{\infty}|t_n^{m-1}|^\rho|t_n|^\rho \leq p(x)^{m-1}\sum_{n=1}^{\infty}|t_n|^\rho < \infty.$$

This implies that $Q^m(x) \in \sigma_\rho$, for any $x \in \sigma_\rho$. Next, for any given $\bar{x} = (\bar{t}_1, \bar{t}_2, \ldots) \in \sigma_\rho$, we prove that $Q^m$ is continuous at $\bar{x}$. Let $V_{J,\varepsilon}(\bar{x})$ be an arbitrarily given $\tau_\rho$-open neighborhood of $\bar{x}$ in $\sigma_\rho$. Since $J \in \mathcal{F}_{\sigma_\rho}$, there is a positive integer $M$ such that

$$J \subseteq \left\{\|\cdot\|_{\rho,k} : k = 1, 2, \ldots, M\right\}. \tag{6.9}$$

We take an $\tau_\rho$-open neighborhood $U_{I,\delta}$ of $\theta_{\sigma_\rho}$ with $I = \left\{\|\cdot\|_{\rho,k} : k = 1, \ldots, M\right\}$ and $\delta = \frac{\varepsilon}{m(p(\bar{x})+1)^{m-1}(\varepsilon+1)}$. It is clear that $\delta < 1$. For $x = \{t_n\}_{n=1}^{\infty} \in U_{I,\delta}$, it satisfies that

$$\max\left\{|t_k - \bar{t}_k|^\rho : \|\cdot\|_{\rho,k} \in I\right\} = \max\left\{\|x - \bar{x}\|_{\rho,k} : \|\cdot\|_{\rho,k} \in I\right\} < \delta. \tag{6.10}$$

Then, by $\delta < 1$, it implies that, for $x = \{t_n\}_{n=1}^{\infty} \in U_{I,\delta}$, we have

$$|t_k| < \mathcal{P}(\bar{x}) + 1, \text{ for } k = 1, 2, \ldots, M.$$

By (6.9), (6.10) and the above inequality, we have

$$\begin{aligned}
&\max\left\{\|Q^m(x) - Q^m(\bar{x})\|_{\rho,k} : \|\cdot\|_{\rho,k} \in J\right\} \\
&\leq \max\left\{\|Q^m(x) - Q^m(\bar{x})\|_{\rho,k} : \|\cdot\|_{\rho,k} \in I\right\} \\
&= \max\left\{|t_k^m - \bar{t}_k^m|^\rho : \|\cdot\|_{\rho,k} \in I\right\} \\
&= \max\left\{\left|t_k^{m-1} + t_k^{m-2}\bar{t}_k + \cdots + \bar{t}_k^{m-1}\right|^\rho |t_k - \bar{t}_k|^\rho : \|\cdot\|_{\rho,k} \in I\right\} \\
&\leq \max\left\{\left(\left|t_k^{m-1}\right|^\rho + \left|t_k^{m-2}\bar{t}_k\right|^\rho + \cdots + \left|\bar{t}_k^{m-1}\right|^\rho\right)|t_k - \bar{t}_k|^\rho : \|\cdot\|_{\rho,k} \in I\right\} \\
&\leq \max\left\{m(\mathcal{P}(\bar{x}) + 1)^{\rho(m-1)}|t_k - \bar{t}_k|^\rho : \|\cdot\|_{\rho,k} \in I\right\} \\
&\leq \max\left\{m(\mathcal{P}(\bar{x}) + 1)^{m-1}|t_k - \bar{t}_k|^\rho : \|\cdot\|_{\rho,k} \in I\right\} \\
&= m(\mathcal{P}(\bar{x}) + 1)^{m-1}\max\left\{|t_k - \bar{t}_k|^\rho : \|\cdot\|_{\rho,k} \in I\right\} \\
&< m(\mathcal{P}(\bar{x}) + 1)^{m-1}\,\delta \\
&= \varepsilon.
\end{aligned}$$

By the extended $\varepsilon$-$\delta$ language of the continuity given in Definition 3.3, this lemma is proved. □

**Theorem 6.4**. *Let m be a positive integer. Then, the m*[th] *power operator $Q^m$ is Fréchet differentiable on $\sigma_\rho$ such that for given $\bar{x} = (\bar{t}_1, \bar{t}_2, \ldots) \in \sigma_\rho$, the Fréchet derivative $\nabla(Q^m)(\bar{x})$ of $Q^m$ at $\bar{x}$ satisfies that,*

(i) *For m > 1, we have*

$$\nabla(Q^m)(\bar{x}) = \begin{pmatrix} m\bar{t}_1^{m-1} & 0 & \dots \\ 0 & m\bar{t}_2^{m-1} & 0 \\ \vdots & 0 & \ddots \end{pmatrix}. \tag{6.12}$$

*Here,* $\nabla(Q^m)(\bar{x})$ *is an* $\infty \times \infty$ *matrix that defines a pointwise multiplication operator on* $\sigma_\rho$ *satisfying that, for every* $x = \{t_n\}_{n=1}^{\infty} \in \sigma_\rho$,

$$\nabla(Q^m)(\bar{x})(x) = \{t_n\}_{n=1}^{\infty} \begin{pmatrix} m\bar{t}_1^{m-1} & 0 & \dots \\ 0 & m\bar{t}_2^{m-1} & 0 \\ \vdots & 0 & \ddots \end{pmatrix} = \{m\bar{t}_n^{m-1} t_n\}_{n=1}^{\infty} \in \sigma_\rho. \tag{6.13}$$

(ii) *When* $m = 1$, *we have*

$$\nabla(Q)(\bar{x}) = \begin{pmatrix} 1 & 0 & \dots \\ 0 & 1 & 0 \\ \vdots & 0 & \ddots \end{pmatrix}.$$

*It is defined as a constant pointwise multiplication operator on* $\sigma_\rho$ *such that*

$$\nabla(Q)(\bar{x})(x) = x = \{t_n\}_{n=1}^{\infty}, \textit{for every } x = \{t_n\}_{n=1}^{\infty} \in \sigma_\rho.$$

*Proof*. Part (ii) is clear. We only prove (i) of this theorem. For an arbitrarily given $\bar{x} = (\bar{t}_1, \bar{t}_2, \dots) \in \sigma_\rho$, by Lemma 6.3 and (6.8), we need to prove the equation (6.13). That is, $m\bar{x}^{m-1} = \{m\bar{t}_n^{m-1}\}_{n=1}^{\infty} \in \sigma_\rho$. By the pointwise multiplication (6.13), it is clear that the operator $\nabla(Q^m)(\bar{x}): \sigma_\rho \to \sigma_\rho$ is linear. Before we prove (6.12), we prove the continuity of $\nabla(Q^m)(\bar{x})$, which is defined by (6.13). Since $\nabla(Q^m)(\bar{x})$ is linear, to this end, we only prove that $\nabla(Q^m)(\bar{x})$ is continuous at the origin $\theta_{\sigma_\rho}$ in $\sigma_\rho$.

Let $V_{J,\varepsilon}(\bar{x})$ be an arbitrarily given $\tau_\rho$-open neighborhood of $\theta_{\sigma_\rho}$ in $\sigma_\rho$, in which $J \in \mathcal{F}_{\sigma_\rho}$ and $\varepsilon > 0$. This implies that there is a positive integer $M$ such that

$$J \subseteq \{\|\cdot\|_{\rho,k}: k = 1, 2, \dots, M\}.$$

We take an $\tau_\rho$-open neighborhood $U_{I,\delta}$ of $\theta_{\sigma_\rho}$ with $I = \{\|\cdot\|_{\rho,k}: k = 1, \dots, M\}$, $\delta = \frac{\varepsilon}{m^\rho (p(\bar{x})+1)^{\rho(m-1)}(\varepsilon+1)}$. It is clear that $I \in \mathcal{F}_{\sigma_\rho}$ and $0 < \delta < 1$. For any $x = \{t_n\}_{n=1}^{\infty} \in U_{I,\delta}$, it satisfies

$$\max\{|t_k|^\rho: \|\cdot\|_{\rho,k} \in I\} = \max\{\|x\|_{\rho,k}: \|\cdot\|_{\rho,k} \in I\} < \delta.$$

Then, by (6.13), we have

$$\begin{aligned} &\max\{\|\nabla(Q^m)(\bar{x})(x)\|_{\rho,k}: \|\cdot\|_{\rho,k} \in J\} \\ &\leq \max\{\|\nabla(Q^m)(\bar{x})(x)\|_{\rho,k}: \|\cdot\|_{\rho,k} \in I\} \\ &= \max\{|m\bar{t}_k^{m-1} t_k|^\rho: \|\cdot\|_{\rho,k} \in I\} \\ &= m^\rho \max\{|\bar{t}_k^{m-1}|^\rho |t_k|^\rho: \|\cdot\|_{\rho,k} \in I\} \\ &\leq m^\rho \mathcal{P}(\bar{x})^{\rho(m-1)} \max\{|t_k|^\rho: \|\cdot\|_{\rho,k} \in I\} \end{aligned}$$

$$< m^{\rho}\mathcal{P}(\bar{x})^{\rho(m-1)}\,\delta$$

$$< \varepsilon.$$

By the extended $\varepsilon$-$\delta$ language of the continuity given in Definition 3.3, this proves the continuity of $\nabla(Q^m)(\bar{x})$ at the origin $\theta_{\sigma_\rho}$ in $\sigma_\rho$. Next, we prove the Fréchet differentiability of $Q^m\colon \sigma_\rho \to \sigma_\rho$ and we prove (6.12). For an arbitrarily given $\bar{x} = (\bar{t}_1, \bar{t}_2, \dots) \in \sigma_\rho$ and for any $u = \{s_n\}_{n=1}^{\infty} \in \sigma_\rho$, we calculate

$$Q^m(\bar{x}+u) - Q^m(\bar{x}) - \nabla(Q^m)(\bar{x})(u)$$

$$= \{(\bar{t}_n + s_n)^m - \bar{t}_n^m - m\bar{t}_n^{m-1}s_n\}_{n=1}^{\infty}$$

$$= \left\{\binom{m}{2}\bar{t}_n^{m-2}s_n^2 + \binom{m}{3}\bar{t}_n^{m-3}s_n^3 + \cdots + \binom{m}{m}s_n^m\right\}_{n=1}^{\infty}.$$

This implies that, for any $k = 1, 2, \dots$, we have

$$\|Q^m(\bar{x}+u) - Q^m(\bar{x}) - \nabla(Q^m)(\bar{x})(u)\|_{\rho,k}$$

$$= \left|\binom{m}{2}\bar{t}_k^{m-2}s_k^2 + \binom{m}{3}\bar{t}_k^{m-3}s_k^3 + \cdots + \binom{m}{m}s_k^m\right|^{\rho}$$

$$= |s_k|^{2\rho}\left|\binom{m}{2}\bar{t}_k^{m-2} + \binom{m}{3}\bar{t}_k^{m-3}s_k + \cdots + \binom{m}{m}s_k^{m-2}\right|^{\rho}. \tag{6.14}$$

Let $W_{K,\varepsilon}$ be an arbitrarily given $\tau_\rho$-open neighborhood of $\theta_{\sigma_\rho}$ in $\sigma_\rho$, in which $K \in \mathcal{F}_{\sigma_\rho}$ and $0 < \varepsilon < 1$. This implies that there is a positive integer $M$ such that

$$K \subseteq \left\{\|\cdot\|_{\rho,k} \colon k = 1, 2, \dots, M\right\}.$$

Recall that, by (6.1), $\mathcal{P}\colon \sigma_\rho \to \mathbb{R}_+$ is well-defined. For the given $\bar{x} = (\bar{t}_1, \bar{t}_2, \dots) \in \sigma_\rho$, we have that

$$\mathcal{P}(\bar{x}) = \sup\{|\bar{t}_n| \colon n = 1, 2, \dots\} = \max\{|\bar{t}_n| \colon n = 1, 2, \dots\} < \infty.$$

We take an $\tau_\rho$-open neighborhood $U_{I,\delta}$ of $\theta_{\sigma_\rho}$ with $I = \left\{\|\cdot\|_{\rho,k} \colon k = 1, \dots, M\right\}$ and $\delta = \frac{\varepsilon}{(\mathcal{P}(\bar{x})+1)^{m\rho}}$. It is clear that $0 < \delta < 1$. Let $u = \{s_n\}_{n=1}^{\infty} \in U_{I,\delta}$, which satisfies

$$\max\left\{|s_i|^{\rho} \colon \|\cdot\|_{\rho,i} \in I\right\} = \max\left\{\|u\|_{\rho,i} \colon \|\cdot\|_{\rho,i} \in I\right\} < \delta. \tag{6.15}$$

By $\delta < 1$, (6.15) implies that $\mathcal{P}_M(u) < 1$, for any $u = \{s_n\}_{n=1}^{\infty} \in U_{I,\delta}$. This means that $u = \{s_n\}_{n=1}^{\infty} \in U_{I,\delta}$, then $|s_n| < 1$, for all $n = 1, 2, \dots M$. Then by (6.14), for any $k = 1, 2, \dots, M$, we estimate

$$\|Q^m(\bar{x}+u) - Q^m(\bar{x}) - \nabla(Q^m)(\bar{x})(u)\|_{\rho,k}$$

$$\leq |s_k|^{2\rho}(|\bar{t}_k| + 1)^{m\rho}$$

$$\leq |s_k|^{2\rho}(\mathcal{P}_M(\bar{x}) + 1)^{m\rho} \tag{6.16}$$

For the arbitrarily given $\tau_\rho$-open neighborhood $W_{K,\varepsilon}$ of $\theta_{\sigma_\rho}$ in $\sigma_\rho$, we first prove that $\nabla(Q^m)(\bar{x})$ satisfies condition (DZ) in Definition 4.8 with respect to the constructed $\tau_\rho$-open neighborhood $U_{I,\delta}$ of $\theta_{\sigma_\rho}$.

Let $u = \{s_n\}_{n=1}^{\infty} \in U_{I,\delta}$. Suppose that $u$ satisfies the following condition

$$\max\{|s_k|^\rho : k = 1, 2, \dots, M\} = \max\{\|u\|_{\rho,k} : k = 1, 2, \dots, M\} = \max\{\|u\|_{\rho,k} : \|\cdot\|_{\rho,k} \in I\} = 0.$$

By $K \subseteq \{\|\cdot\|_{\rho,k} : k = 1, 2, \dots, M\} = I$. This implies that $s_k = 0$, for $k = 1, 2, \dots, M$. Then, we obtain that

$$\begin{aligned}
&\max\{\|Q^m(\bar{x} + u) - Q^m(\bar{x}) - \nabla(Q^m)(\bar{x})(u)\|_{\rho,k} : \|\cdot\|_{\rho,k} \in K\} \\
&\leq \max\{\|Q^m(\bar{x} + u) - Q^m(\bar{x}) - \nabla(Q^m)(\bar{x})(u)\|_{\rho,k} : \|\cdot\|_{\rho,k} \in I\} \\
&= \max\{\|Q^m(\bar{x}) - Q^m(\bar{x})\|_{\rho,k} : \|\cdot\|_{\rho,k} \in I\} \\
&= 0.
\end{aligned}$$

Hence, we proved that for the arbitrary $\tau_\rho$-open neighborhood $W_{K,\varepsilon}$ of $\theta_{\sigma_\rho}$ in $\sigma_\rho$, $\nabla(Q^m)(\bar{x})$ satisfies condition (DZ) in Definition 4.8 with respect to the constructed $\tau_\rho$-open neighborhood $U_{I,\delta}$ of $\theta_{\sigma_\rho}$.

Next, for the arbitrary $\tau_\rho$-open neighborhood $W_{K,\varepsilon}$ of $\theta_{\sigma_\rho}$ in $\sigma_\rho$, we prove that $\nabla(Q^m)(\bar{x})$ satisfies condition (DR) in Definition 4.8 with respect to the constructed $\tau_\rho$-open neighborhood $U_{I,\delta}$ of $\theta_{\sigma_\rho}$. Let $u = \{s_n\}_{n=1}^\infty \in U_{I,\delta}$. Suppose that

$$0 < \max\{|s_k|^\rho : k = 1, 2, \dots, M\} = \max\{\|u\|_{\rho,k} : k = 1, 2, \dots, M\} = \max\{\|u\|_{\rho,k} : \|\cdot\|_{\rho,k} \in I\} < \delta.$$

Notice that, for this given fixed $0 < \rho < 1$, $t^\rho$ is a strictly increasing function for $t \in [0, \infty)$. Then by (6.15), (6.14) and (6.16), for $u = \{s_n\}_{n=1}^\infty \in U_{I,\delta}$, we have

$$\begin{aligned}
&\max\left\{\left\|\frac{Q^m(\bar{x}+u) - Q^m(\bar{x}) - \nabla(Q^m)(\bar{x})(u)}{\max\{|s_i|^\rho : \|\cdot\|_i \in I\}}\right\|_{\rho,k} : \|\cdot\|_{\rho,k} \in K\right\} \\
&= \max\left\{\left|\frac{\binom{m}{2}\bar{t}_k^{m-2} s_k^2 + \binom{m}{3}\bar{t}_k^{m-3} s_k^3 + \cdots + \binom{m}{m} s_k^m}{\max\{|s_i|^\rho : \|\cdot\|_i \in I\}}\right|^\rho : \|\cdot\|_{\rho,k} \in K\right\} \\
&= \max\left\{\frac{|s_k|^{2\rho}\left|\binom{m}{2}\bar{t}_k^{m-2} + \binom{m}{3}\bar{t}_k^{m-3} s_k + \cdots + \binom{m}{m} s_k^{m-2}\right|^\rho}{(\max\{|s_i|^\rho : \|\cdot\|_i \in I\})^\rho} : \|\cdot\|_{\rho,k} \in K\right\} \\
&\leq \max\left\{\frac{|s_k|^{2\rho}(\mathcal{P}_M(\bar{x}) + 1)^{m\rho}}{(\max\{|s_i|^\rho : \|\cdot\|_i \in I\})^\rho} : \|\cdot\|_{\rho,k} \in K\right\} \\
&= (\mathcal{P}(\bar{x}) + 1)^{m\rho} \max\left\{\frac{|s_k|^{2\rho}}{(\max\{|s_i|^\rho : \|\cdot\|_i \in I\})^\beta} : \|\cdot\|_{\rho,k} \in K\right\} \\
&= (\mathcal{P}(\bar{x}) + 1)^{m\rho} \frac{\max\{|s_k|^{2\rho} : \|\cdot\|_{\rho,k} \in K\}}{(\max\{|s_i|^\rho : \|\cdot\|_i \in I\})^\rho} \\
&= (\mathcal{P}(\bar{x}) + 1)^{m\rho} \frac{\left(\max\{|s_k|^\rho : \|\cdot\|_{\rho,k} \in K\}\right)^2}{\left(\max\{|s_i|^\rho : \|\cdot\|_{\rho,i} \in I\}\right)^\rho} \\
&\leq (\mathcal{P}(\bar{x}) + 1)^{m\rho} \frac{\left(\max\{|s_k|^\rho : \|\cdot\|_{\rho,i} \in I\}\right)^2}{\left(\max\{|s_i|^\rho : \|\cdot\|_{\rho,i} \in I\}\right)^\rho} \\
&= (\mathcal{P}(\bar{x}) + 1)^{m\rho} \left(\max\{|s_k|^\rho : \|\cdot\|_{\rho,i} \in I\}\right)^{2-\rho}
\end{aligned}$$

$$< (\mathcal{P}(\bar{x}) + 1)^{m\rho}\, \delta^{2-\rho}$$

$$< (\mathcal{P}(\bar{x}) + 1)^{m\rho}\, \delta$$

$$< \varepsilon.$$

Hence, we proved that for the arbitrarily given $\tau_\beta$-open neighborhood $W_{K,\varepsilon}$ of $\theta_{\sigma_\rho}$ in $\sigma_\rho$, $\nabla(Q^m)(\bar{x})$ satisfies condition (DR) in Definition 4.8 with respect to the constructed $\tau_\rho$-open neighborhood $U_{I,\delta}$ of $\theta_{\sigma_\rho}$. This completes the proof of this Theorem. □

**Remarks 6.5**. By the property of Fréchet differentiability given in Theorem 4.16, the continuity of the *m*th power operator $Q^m: \sigma_\rho \to \sigma_\rho$, Lemma 6.3 follows immediately from the Fréchet differentiability of $Q^m$ proved in Theorem 6.4.

Similarly to (5.17), we define a polynomial type operator in $\sigma_\rho$. Let *m* be a positive integer. Let $a_m$, $a_{m-1}$, … $a_1$ be real numbers. Define a polynomial type operator $Q_m: \sigma_\rho \to \sigma_\rho$, for $x = \{t_n\}_{n=1}^\infty \in \sigma_\rho$, by

$$Q_m(x) = a_m x^m + a_{m-1} x^{m-1} + \ldots + a_1 x = \left\{\sum_{i=1}^m a_i t_n^i\right\}_{n=1}^\infty. \tag{6.17}$$

**Theorem 6.6**. *Let m be a positive integer and let $Q_m$ be a polynomial type operator on $\sigma_\rho$ defined by* (6.17). *Then, $Q_m$ is Fréchet differentiable on $\sigma_\rho$ such that for any given $\bar{x} = \{\bar{t}_n\}_{n=1}^\infty \in \sigma_\rho$, the Fréchet derivative $\nabla Q_m(\bar{x})$ of $Q_m$ at $\bar{x}$ satisfies that*

(i) *If m* > 1, *we have*

$$\nabla Q_m(\bar{x}) = \begin{pmatrix} \sum_{i=1}^m i a_i \bar{t}_1^{i-1} & 0 & \cdots \\ 0 & \sum_{i=1}^m i a_i \bar{t}_2^{i-1} & 0 \\ \vdots & 0 & \ddots \end{pmatrix}.$$

*Here, $\nabla Q_m(\bar{x})$ is represented by an $\infty \times \infty$ matrix that defines a pointwise multiplication operator on $\sigma_\rho$ satisfying that, for every $x = \{t_n\}_{n=1}^\infty \in \sigma_\rho$,*

$$\nabla Q_m(\bar{x})(x) = \{t_n\}_{n=1}^\infty \begin{pmatrix} \sum_{i=1}^m i a_i \bar{t}_1^{i-1} & 0 & \cdots \\ 0 & \sum_{i=1}^m i a_i \bar{t}_2^{i-1} & 0 \\ \vdots & 0 & \ddots \end{pmatrix} = \left\{\sum_{i=1}^m i a_i \bar{t}_n^{i-1} t_n\right\}_{n=1}^\infty \in \sigma_\rho.$$

(ii) *When m* = 1, *we have*

$$\nabla Q_1(\bar{x}) = \begin{pmatrix} a_1 & 0 & \cdots \\ 0 & a_1 & 0 \\ \vdots & 0 & \ddots \end{pmatrix}.$$

*For every $x = \{t_n\}_{n=1}^\infty \in \sigma_\rho$,*

$$\nabla Q_1(\bar{x})x = a_1 x = \{a_1 t_n\}_{n=1}^\infty \in \sigma_\rho.$$

*Proof*. This theorem follows from Theorem 6.4 and Theorem 4.14 immediately. □

The space $(\sigma_\rho, \tau_\rho)$ is a Hausdorff topological vector space and it is not seminorm-constructed. Hence, we cannot apply Theorem 4.18 and Theorem 6.4 to prove the Gâteaux differentiability of the *m*th power operator $Q^m$ $\sigma_\rho$. To end this section, we will directly prove the Gâteaux differentiability of polynomial type operators $Q^m$ in the space $(\sigma_\rho, \tau_\rho)$ without using its Fréchet differentiability. Meanwhile, similarly to Theorem 5.10, we will find that the Gâteaux derivatives of polynomial type operators $Q^m$ coincide with their Fréchet derivatives.

**Theorem 6.7**. *Let m be a positive integer. Then, the mth power operator $Q^m$ is Gâteaux differentiable on $\sigma_\rho$ such that for any $\bar{x} = (\bar{t}_1, \bar{t}_2, \dots) \in \sigma_\rho$, the Gâteaux derivative $(Q^m)'(\bar{x})$ of $Q^m$ at $\bar{x}$ satisfies that,*

(i) *For m > 1, we have*

$$(Q^m)'(\bar{x}) = \nabla(Q^m)(\bar{x}) = \begin{pmatrix} m\bar{t}_1^{m-1} & 0 & \dots \\ 0 & m\bar{t}_2^{m-1} & 0 \\ \vdots & 0 & \ddots \end{pmatrix}.$$

*Here, for every $v = \{s_n\}_{n=1}^{\infty} \in \sigma_\rho \backslash \{\theta_{\sigma_\rho}\}$*

$$(Q^m)'(x)(v) = \nabla(Q^m)(\bar{x})(v) = \{s_n\}_{n=1}^{\infty} \begin{pmatrix} m\bar{t}_1^{m-1} & 0 & \dots \\ 0 & m\bar{t}_2^{m-1} & 0 \\ \vdots & 0 & \ddots \end{pmatrix} = \{m\bar{t}_n^{m-1} s_n\}_{n=1}^{\infty} \in \sigma_\rho.$$

(ii) *When m = 1, we have*

$$Q'(\bar{x}) = \nabla(Q)(\bar{x}) = \begin{pmatrix} 1 & 0 & \dots \\ 0 & 1 & 0 \\ \vdots & 0 & \ddots \end{pmatrix}.$$

*It is defined as a constant pointwise multiplication operator on $\sigma_\rho$ such that*

$$Q'(\bar{x})(v) = v = \{s_n\}_{n=1}^{\infty}, \textit{for every } v = \{s_n\}_{n=1}^{\infty} \in \sigma_\rho \backslash \{\theta_{\sigma_\rho}\}.$$

*Proof*. Part (ii) is clear. So, we only prove (i). Let $\bar{x} = (\bar{t}_1, \bar{t}_2, \dots) \in \sigma_\rho$. Let $v = \{s_n\}_{n=1}^{\infty} \in \sigma_\rho \backslash \{\theta_{\sigma_\rho}\}$ be arbitrarily given (fixed). We want to prove that the Gâteaux derivative $(Q^m)'(\bar{x})(v)$ of $Q^m$ at $\bar{x}$ along direction *v* satisfies $(Q^m)'(\bar{x})(v) = \nabla(Q^m)(\bar{x})(v)$. To this end, let $J \in \mathcal{F}_{\sigma_\rho}$ be arbitrarily given. For any real number $t \neq 0$, similarly to (6.14), we calculate

$$\max\left\{\left\|\frac{Q^m(\bar{x}+tv)-Q^m(\bar{x})-t\nabla(Q^m)(\bar{x})(v)}{t}\right\|_{\rho,n(i)} : \|\cdot\|_{\rho,n(i)} \in J\right\}$$

$$= \max\left\{\left|\frac{\binom{m}{2}\bar{t}_{n(i)}^{m-2}t^2 s_{n(i)}^2 + \binom{m}{3}\bar{t}_{n(i)}^{m-3}t^3 s_{n(i)}^3 + \dots + \binom{m}{m}t^m s_{n(i)}^m}{t}\right|^\rho : \|\cdot\|_{\rho,n(i)} \in J\right\}$$

$$= \max\left\{\frac{\left|ts_{n(i)}\right|^{2\rho}\left|\binom{m}{2}\bar{t}_{n(i)}^{m-2} + \binom{m}{3}\bar{t}_{n(i)}^{m-3}s_{n(i)} + \dots + \binom{m}{m}s_{n(i)}^{m-2}\right|^\rho}{|t|^\rho} : \|\cdot\|_{\rho,n(i)} \in J\right\}$$

$$= \max\left\{|t|^\rho\left|s_{n(i)}\right|^{2\rho}\left|\binom{m}{2}\bar{t}_{n(i)}^{m-2} + \binom{m}{3}\bar{t}_{n(i)}^{m-3}s_{n(i)} + \dots + \binom{m}{m}s_{n(i)}^{m-2}\right|^\rho : \|\cdot\|_{\rho,n(i)} \in J\right\}$$

$$= |t|^\rho \max\left\{\left|s_{n(i)}\right|^{2\rho}\left|\binom{m}{2}\bar{t}_{n(i)}^{m-2} + \binom{m}{3}\bar{t}_{n(i)}^{m-3}s_{n(i)} + \dots + \binom{m}{m}s_{n(i)}^{m-2}\right|^\rho : \|\cdot\|_{\rho,n(i)} \in J\right\}$$

$\to 0$, as $t \to 0$.

Here, for the given and fixed $\bar{x} = (\bar{t}_1, \bar{t}_2, \dots) \in \zeta_\beta$, $v = \{s_n\}_{n=1}^{\infty} \in \sigma_\rho \backslash \{\theta_{\sigma_\rho}\}$ and with the finite set $J$, the factor before the above limit satisfies

$$\max\left\{\left|s_{n(i)}\right|^{2\rho}\left|\binom{m}{2}\bar{t}_{n(i)}^{m-2} + \binom{m}{3}\bar{t}_{n(i)}^{m-3}s_{n(i)} + \cdots + \binom{m}{m}s_{n(i)}^{m-2}\right|^{\rho} : \|\cdot\|_{\rho,n(i)} \in J\right\} < \infty. \qquad \square$$

Similarly to Theorem 6.6, we can prove the following theorem by Theorem 6.7.

**Corollary 6.8**. *Let m be a positive integer and let $Q_m$ be a polynomial type operator on $\sigma_\rho$ defined by* (6.17). *Then, $Q_m$ is Gâteaux differentiable on $\sigma_\rho$ such that, for any $\bar{x} = (\bar{t}_1, \bar{t}_2, \dots) \in \sigma_\rho$, the Gâteaux derivative $(Q_m)'(\bar{x})$ of $Q_m$ at $\bar{x}$ satisfies that*

(i) *If $m > 1$, we have*

$$(Q_m)'(\bar{x}) = \nabla Q_m(\bar{x}) = \begin{pmatrix} \sum_{i=1}^{m} i a_i \bar{t}_1^{i-1} & 0 & \dots \\ 0 & \sum_{i=1}^{m} i a_i \bar{t}_2^{i-1} & 0 \\ \vdots & 0 & \ddots \end{pmatrix}.$$

*Here, $(Q_m)'(\bar{x})$ is represented by an $\infty \times \infty$ matrix that defines a pointwise multiplication operator on $\sigma_\rho$ satisfying that, for every $v = \{t_n\}_{n=1}^{\infty} \in \sigma_\rho \backslash \{\theta_{\sigma_\rho}\}$,*

$$(Q_m)'(\bar{x})(v) = \nabla Q_m(\bar{x})(v) = \{t_n\}_{n=1}^{\infty} \begin{pmatrix} \sum_{i=1}^{m} i a_i \bar{t}_1^{i-1} & 0 & \dots \\ 0 & \sum_{i=1}^{m} i a_i \bar{t}_2^{i-1} & 0 \\ \vdots & 0 & \ddots \end{pmatrix} = \left\{\sum_{i=1}^{m} i a_i \bar{t}_n^{i-1} t_n\right\}_{n=1}^{\infty}.$$

(ii) *When $m = 1$, we have*

$$(Q_1)'(\bar{x}) = \nabla Q_1(\bar{x}) = \begin{pmatrix} a_1 & 0 & \dots \\ 0 & a_1 & 0 \\ \vdots & 0 & \ddots \end{pmatrix}.$$

*For every $v = \{t_n\}_{n=1}^{\infty} \in \sigma_\rho \backslash \{\theta_{\sigma_\rho}\}$, we have*

$$(Q_m)'(\bar{x})(v) = \nabla Q_1(\bar{x})v = a_1 v = \{a_1 t_n\}_{n=1}^{\infty} \in \sigma_\rho.$$

## 7. Another Example of Not Seminorm-Constructed Topological Vector Space $(\mathbb{S}, \tau_{\mathbb{S}})$

Let $\mathbb{S}$ denote the collection of all sequences of real numbers, which is a real vector space with its origin denoted by $\theta = (0, 0, \dots)$. We construct a countable family $\mathbb{F}_{\mathbb{S}}$ of *F*-seminorms $\|\cdot\|_k$ on $\mathbb{S}$ (the notation is different from the *F*-seminorms $\|\cdot\|_{\rho,k}$ on $\sigma_\rho$ used in the previous section) as follows. For $k \in \mathbb{N}$, define $\|\cdot\|_k \colon \mathbb{S} \to \mathbb{R}_+$ by

$$\|x\|_k = \frac{|t_k|}{1+|t_k|}, \text{ for any } x = \{t_n\}_{n=1}^{\infty} \in \mathbb{S}.$$

$\|\cdot\|_k$ has the following properties.

(i) $\|\cdot\|_k$ is a *F*-seminorm on $\mathbb{S}$, for every $k \in \mathbb{N}$;

(ii) Let *a* be a real number. Then,

(a) $|a| < 1 \implies \|ax\|_k \geq |a|\|x\|_k$;

(b) $|a| \geq 1 \implies \|ax\|_k \leq |a|\|x\|_k$.

Let $\tau_{\mathbb{S}}$ be the topology on $\mathbb{S}$ that is generated by this countable family $\mathbb{F}_{\mathbb{S}}$ of $F$-seminorms $\|\cdot\|_k$ on $\mathbb{S}$. Then, $(\mathbb{S}, \tau_{\mathbb{S}})$ is a topological vector space that is not have seminorm constructed. Let $\mathcal{F}_{\mathbb{S}}$ denote the collection of all nonempty finite subsets of $\mathbb{F}_{\mathbb{S}}$. Similarly to the previous two sections, we consider the Fréchet differentiability and Gâteaux differentiability of some polynomial type operator on the space $\mathbb{S}$. Let $m$ be a positive integer. The $m^{\text{th}}$ power operator $R^m: \mathbb{S} \to \mathbb{S}$ is defined by

$$R^m(x) = \{t_n^m\}_{n=1}^{\infty}, \text{ for every } x = \{t_n\}_{n=1}^{\infty} \in \mathbb{S}. \tag{7.1}$$

If $m = 1$, then $R^1$ is the identity mapping that is a continuous and linear mapping on $\mathbb{S}$. Hence, we only need to consider the cases $m > 1$. We first prove the continuity of the $m^{\text{th}}$ power mapping $R^m: \mathbb{S} \to \mathbb{S}$. Then, in the following proposition, we prove the Fréchet differentiability of $R^m$ on $\mathbb{S}$, which is stronger than the continuity. However, in order to practice the definition of the continuity (Definition 3.1), we will provide a direct proof for the continuity of $R^m$ on $\mathbb{S}$ first.

**Lemma 7.1**. *Let m be a positive integer with m > 1. Then the $m^{\text{th}}$ power operator $R^m: \mathbb{S} \to \mathbb{S}$ defined by* (7.1) *is a $\tau_{\mathbb{S}}$-continuous operator on* $\mathbb{S}$.

*Proof.* It is clear to see that the $m$th power mapping $R^m$ is well-defined on $\mathbb{S}$. For any $x = \{t_n\}_{n=1}^{\infty} \in \mathbb{S}$, by (7.1), we have

$$\|x\|_k = \frac{|t_k|}{1+|t_k|} \quad \text{and} \quad \|R^m(x)\|_k = \frac{|t_k|^m}{1+|t_k|^m}, \text{ for } k = 1, 2, \dots .$$

Next, for any given $\bar{x} = (\bar{t}_1, \bar{t}_2, \dots) \in \mathbb{S}$, we prove that $R^m$ is continuous at $\bar{x}$. Let $J \in \mathcal{F}_{\mathbb{S}}$ and $0 < \varepsilon < 1$ be arbitrarily given. Since $J \in \mathcal{F}_{\mathbb{S}}$, there is a positive integer $M$ such that

$$J \subseteq \{\|\cdot\|_k : k = 1, 2, \dots, M\}.$$

Let $I = \{\|\cdot\|_k : k = 1, \dots, M\}$. For this given fixed positive integer $M$, for any $x = (t_1, t_2, \dots) \in \mathbb{S}$, we write

$$\mathcal{P}_M(x) = \max\{|t_k| : k = 1, 2, \dots, M\}.$$

For the arbitrarily given $\bar{x} = (\bar{t}_1, \bar{t}_2, \dots) \in \mathbb{S}$, let $\delta = \frac{\varepsilon}{(1+2\mathcal{P}_M(\bar{x}))^m(\varepsilon+1)}$. It is clear that $\delta < \frac{1}{2}$. Then, we have a $\tau_{\mathbb{S}}$-open neighborhood $U_{I,\delta}(\bar{x})$ of $\bar{x}$ in $\mathbb{S}$. Then, for any $x = \{t_n\}_{n=1}^{\infty} \in U_{I,\delta}(\bar{x})$, it satisfies

$$\max\left\{\frac{|t_k - \bar{t}_k|}{1+|t_k - \bar{t}_k|} : \|\cdot\|_k \in I\right\} = \max\{\|x - \bar{x}\|_k : \|\cdot\|_k \in I\} < \delta. \tag{7.2}$$

By $\delta < \frac{1}{2}$, this implies that

$$\max\{|t_k - \bar{t}_k| : \|\cdot\|_k \in I\} < \frac{\delta}{1-\delta} < 1, \text{ for any } x = \{t_n\}_{n=1}^{\infty} \in U_{I,\delta}(\bar{x}).$$

This induces that, for any $k = 1, 2, \dots, M$, we have

$$|t_k| < 1 + |\bar{t}_k| \leq 1 + \mathcal{P}_M(\bar{x}), \text{ for any } x = \{t_n\}_{n=1}^{\infty} \in U_{I,\delta}(\bar{x}). \tag{7.3}$$

Then, by the property that $J \subseteq I$ and (7.2) and (7.3), we have

$$\max\{\|R^m(x) - R^m(\bar{x})\|_k : \|\cdot\|_k \in J\}$$

$$\leq \max \{\|R^m(x) - R^m(\bar{x})\|_k : \|\cdot\|_k \in I\}$$

$$= \max \left\{\frac{|t_k^m - \bar{t}_k^m|}{1+|t_k^m - \bar{t}_k^m|} : \|\cdot\|_k \in I\right\}$$

$$= \max \left\{\frac{|t_k^{m-1} + t_k^{m-2}\bar{t}_2 + \cdots + \bar{t}_k^{m-1}||t_k - \bar{t}_k|}{1+|t_k^{m-1} + t_k^{m-2}\bar{t}_2 + \cdots + \bar{t}_k^{m-1}||t_k - \bar{t}_k|} : \|\cdot\|_k \in I\right\}$$

$$\leq \max \left\{\frac{(|t_k^{m-1}| + |t_k^{m-2}||\bar{t}_k| + \cdots + |\bar{t}_k^{m-1}|)|t_k - \bar{t}_k|}{1+(|t_k^{m-1}| + |t_k^{m-2}||\bar{t}_k| + \cdots + |\bar{t}_k^{m-1}|)|t_k - \bar{t}_k|} : \|\cdot\|_k \in I\right\}$$

$$\leq \max \left\{\frac{((1+|\bar{t}_k|)^{m-1} + (1+|\bar{t}_k|)^{m-2}|\bar{t}_k| + \cdots + |\bar{t}_k^{m-1}|)|t_k - \bar{t}_k|}{1+((1+|\bar{t}_k|)^{m-1} + (1+|\bar{t}_k|)^{m-2}|\bar{t}_k| + \cdots + |\bar{t}_k^{m-1}|)|t_k - \bar{t}_k|} : \|\cdot\|_k \in I\right\}$$

$$\leq \max \left\{\frac{(1+2|\bar{t}_k|)^m |t_k - \bar{t}_k|}{1+(1+2|\bar{t}_k|)^m |t_k - \bar{t}_k|} : \|\cdot\|_k \in I\right\}$$

$$\leq (1 + 2\mathcal{P}_M(\bar{x}))^m \max \left\{\frac{|t_k - \bar{t}_k|}{1+(1+2q_M(\bar{x}))^m |t_k - \bar{t}_k|} : \|\cdot\|_k \in I\right\}$$

$$\leq (1 + 2\mathcal{P}_M(\bar{x}))^m \max \left\{\frac{|t_k - \bar{t}_k|}{1+|t_k - \bar{t}_k|} : \|\cdot\|_k \in I\right\}$$

$$< (1 + 2\mathcal{P}_M(\bar{x}))^m \, \delta$$

$$< \varepsilon.$$

By the extended $\varepsilon$-$\delta$ language of the continuity given in Definition 3.3, this lemma is proved. □

Next, we prove the Fréchet differentiability of the power operator $R^m$ on $\mathbb{S}$. We see that the formulations of the Fréchet derivatives of the power operator $R^m$ on $\mathbb{S}$ given by (7.4) and (7.5) below are formally same to (6.12) and (6.13).

**Theorem 7.2**. *Let m be a positive integer. Then* $R^m : \mathbb{S} \to \mathbb{S}$ *is Fréchet differentiable on* $\mathbb{S}$ *such that for any* $\bar{x} = (\bar{t}_1, \bar{t}_2, \ldots) \in \mathbb{S}$, $\nabla(R^m)(\bar{x})$ *has the following representation.*

(i) *For m* > 1, *we have*

$$\nabla(R^m)(\bar{x}) = \begin{pmatrix} m\bar{t}_1^{m-1} & 0 & \ldots \\ 0 & m\bar{t}_2^{m-1} & 0 \\ \vdots & 0 & \ddots \end{pmatrix}. \tag{7.4}$$

*Here,* $\nabla(R^m)(\bar{x})$ *is an* $\infty \times \infty$ *matrix that defines a pointwise multiplication operator on* $\mathbb{S}$ *satisfying that, for every* $x = \{t_n\}_{n=1}^{\infty} \in \mathbb{S}$, *we have*

$$\nabla(R^m)(\bar{x})(x) = \{t_n\}_{n=1}^{\infty} \begin{pmatrix} m\bar{t}_1^{m-1} & 0 & \ldots \\ 0 & m\bar{t}_2^{m-1} & 0 \\ \vdots & 0 & \ddots \end{pmatrix} = \{m\bar{t}_n^{m-1} t_n\}_{n=1}^{\infty} \in \mathbb{S}. \tag{7.5}$$

(ii) *In particular*, *when m* = 1, *we have*

$$\nabla(R)(\bar{x}) = \begin{pmatrix} 1 & 0 & \ldots \\ 0 & 1 & 0 \\ \vdots & 0 & \ddots \end{pmatrix}.$$

*It is defined as a constant pointwise multiplication operator on* $\mathbb{S}$ *such that*

$$\nabla(R)(\bar{x})(x) = x = \{t_n\}_{n=1}^{\infty}, \textit{for every } x = \{t_n\}_{n=1}^{\infty} \in \mathbb{S}.$$

*Proof*. If $m = 1$, then this proposition is clear. So, we only prove the case of $m > 1$. The ideas of the proof of this theorem are similar to the proof of Theorem 6.4. However, the techniques used in this proof are different from the proof of Theorem 6.4. Hence, we show the details of this proof.

For an arbitrarily given $\bar{x} = (\bar{t}_1, \bar{t}_2, \dots) \in \mathbb{S}$, by definition (7.4), $m\bar{x}^{m-1} = \{m\bar{t}_n^{m-1}\}_{n=1}^{\infty} \in \mathbb{S}$ and the pointwise multiplication $\nabla(R^m)(\bar{x})\colon \mathbb{S} \to \mathbb{S}$ is linear. Before we prove (7.4), we prove the continuity of $\nabla(R^m)(\bar{x})$ given in (7.4). To this end, by the linearity of $\nabla(R^m)(\bar{x})$, we only prove that $\nabla(R^m)(\bar{x})$ is continuous at the origin $\theta_{\mathbb{S}}$.

Let $J \in \mathcal{F}_{\mathbb{S}}$ and $\varepsilon > 0$ be arbitrarily given. Then, there is a positive integer $M$ such that

$$J \subseteq \{\|\cdot\|_k \colon k = 1, 2, \dots, M\}.$$

Let $I = \{\|\cdot\|_k \colon k = 1, \dots, M\}$. Similarly to the proof of Lemma 7.1, for this given $M$, we write

$$\mathcal{P}_M(x) = \max\{|t_k| \colon k = 1, 2, \dots, M\}, \text{ for any } x = (t_1, t_2, \dots) \in \mathbb{S}. \tag{7.6}$$

For the arbitrarily given $\bar{x} = (\bar{t}_1, \bar{t}_2, \dots) \in \mathbb{S}$, let $\delta = \frac{\varepsilon}{(1+2\mathcal{P}_M(\bar{x}))^m(\varepsilon+1)}$. It is clear that $\delta < \frac{1}{2}$. Then, we have a $\tau_{\mathbb{S}}$-open neighborhood $U_{I,\delta}$ of $\theta_{\mathbb{S}}$ in $\mathbb{S}$. Then, for any $x = \{t_n\}_{n=1}^{\infty} \in U_{I,\delta}$, it satisfies

$$\max\left\{\frac{|t_k|}{1+|t_k|} \colon \|\cdot\|_k \in I\right\} = \max\{\|x - \bar{x}\|_k \colon \|\cdot\|_k \in I\} < \delta.$$

Then, by (7.5), (7.6) and the above inequality, similar to the proof of Lemma 7.1, we have

$$\begin{aligned}
&\max\{\|\nabla(R^m)(\bar{x})(x)\|_k \colon \|\cdot\|_k \in J\} \\
&\leq \max\{\|\nabla(R^m)(\bar{x})(x)\|_k \colon \|\cdot\|_k \in I\} \\
&= \max\left\{\frac{|m\bar{t}_k^{m-1}t_k|}{1+|m\bar{t}_k^{m-1}t_k|} \colon \|\cdot\|_k \in I\right\} \\
&= m\max\left\{\frac{|\bar{t}_k^{m-1}||t_k|}{1+m|\bar{t}_k^{m-1}||t_k|} \colon \|\cdot\|_k \in I\right\} \\
&\leq m\max\left\{\frac{|\bar{t}_k^{m-1}||t_k|}{1+|\bar{t}_k^{m-1}||t_k|} \colon \|\cdot\|_k \in I\right\} \\
&\leq m\max\left\{\frac{(|\bar{t}_k|+1)^{m-1}|t_k|}{1+(|\bar{t}_k|+1)^{m-1}|t_k|} \colon \|\cdot\|_k \in I\right\} \\
&\leq m(\mathcal{P}_M(\bar{x})+1)^{m-1}\max\left\{\frac{|t_k|}{1+(|\bar{t}_k|+1)^{m-1}|t_k|} \colon \|\cdot\|_k \in I\right\} \\
&\leq m(\mathcal{P}_M(\bar{x})+1)^{m-1}\max\left\{\frac{|t_k|}{1+|t_k|} \colon \|\cdot\|_k \in I\right\} \\
&< m(\mathcal{P}_M(\bar{x})+1)^{m-1}\,\delta \\
&< \varepsilon.
\end{aligned}$$

By the extended $\varepsilon$-$\delta$ language of the continuity given in Definition 3.3, this proves the continuity of $\nabla(R^m)(\bar{x})$ at the origin $\theta_{\mathbb{S}}$ in $\mathbb{S}$. Next, we prove the Fréchet differentiability of $R^m: \mathbb{S} \to \mathbb{S}$. More precisely speaking, we prove (7.4) and (7.5). For any $u = \{s_n\}_{n=1}^{\infty} \in \mathbb{S}$, we calculate

$$R^m(\bar{x}+u) - R^m(\bar{x}) - \nabla(R^m)(\bar{x})(u)$$

$$= \{(\bar{t}_n + s_n)^m - \bar{t}_n^m - m\bar{t}_n^{m-1}s_n\}_{n=1}^{\infty}$$

$$= \left\{\binom{m}{2}\bar{t}_n^{m-2}s_n^2 + \binom{m}{3}\bar{t}_n^{m-3}s_n^3 + \cdots + \binom{m}{m}s_n^m\right\}_{n=1}^{\infty}. \tag{7.7}$$

This implies that, for any $k \in \mathbb{N}$, we have

$$\|R^m(\bar{x}+u) - R^m(\bar{x}) - \nabla(R^m)(\bar{x})(u)\|_k$$

$$= \frac{\left|\binom{m}{2}\bar{t}_k^{m-2}s_k^2 + \binom{m}{3}\bar{t}_k^{m-3}s_k^3 + \cdots + \binom{m}{m}s_k^m\right|}{1 + \left|\binom{m}{2}\bar{t}_k^{m-2}s_k^2 + \binom{m}{3}\bar{t}_k^{m-3}s_k^3 + \cdots + \binom{m}{m}s_k^m\right|}$$

$$= \frac{|s_k|^2\left|\binom{m}{2}\bar{t}_k^{m-2} + \binom{m}{3}\bar{t}_k^{m-3}s_k^1 + \cdots + \binom{m}{m}s_k^{m-2}\right|}{1 + |s_k|^2\left|\binom{m}{2}\bar{t}_k^{m-2} + \binom{m}{3}\bar{t}_k^{m-3}s_k^1 + \cdots + \binom{m}{m}s_k^{m-2}\right|}. \tag{7.8}$$

Let $K \in \mathcal{F}_{\mathbb{S}}$ and $0 < \varepsilon < 1$ be an arbitrarily given, which induces that $W_{K,\varepsilon}$ is an arbitrarily given $\tau_{\mathbb{S}}$-open neighborhood of $\theta_{\mathbb{S}}$ in $\mathbb{S}$. This implies that there is a positive integer $M$ such that

$$K \subseteq \{\|\cdot\|_k : k = 1, 2, \ldots, M\}.$$

For the fixed positive integer $M$, by (7.6), we have $\mathcal{P}_M(x) = \max\{|t_k| : k = 1, 2, \ldots, M\}$, for any $x = (t_1, t_2, \ldots) \in \mathbb{S}$. It is clear that, for $\bar{x} = (\bar{t}_1, \bar{t}_2, \ldots) \in \mathbb{S}$, we have $\mathcal{P}_M(\bar{x}) < \infty$. Let

$$I = \{\|\cdot\|_k : k = 1, 2, \ldots, M\} \in \mathcal{F}_{\mathbb{S}} \quad \text{and} \quad \delta = \frac{\varepsilon}{2\left(1 + \mathcal{P}_M(\bar{x})\right)^m}.$$

It is clear that $\delta < \frac{1}{2}$. We have an $\tau_{\mathbb{S}}$-open neighborhood $U_{I,\delta}$ of $\theta_{\mathbb{S}}$ in $\mathbb{S}$. Then, for the arbitrarily given $K \in \mathcal{F}_{\mathbb{S}}$ and $0 < \varepsilon < 1$, we first prove that $\nabla(R^m)(\bar{x})$ satisfies condition (DZ) in Definition 4.8 with respect to the constructed $\tau_{\mathbb{S}}$-open neighborhood $U_{I,\delta}$ of $\theta_{\mathbb{S}}$. Let $u = \{s_n\}_{n=1}^{\infty} \in U_{I,\delta}$. Suppose that

$$\max\left\{\frac{|s_k|}{1+|s_k|} : k = 1, 2, \ldots, M\right\} = \max\{\|u\|_k : \|\cdot\|_k \in I\} = 0.$$

This implies that $s_k = 0$, for $k = 1, 2, \ldots, M$. Then, we obtain that

$$\max\{\|R^m(\bar{x}+u) - R^m(\bar{x}) - \nabla(R^m)(\bar{x})(u)\|_k : \|\cdot\|_k \in K\}$$

$$\leq \max\{\|R^m(\bar{x}+u) - R^m(\bar{x}) - \nabla(R^m)(\bar{x})(u)\|_k : \|\cdot\|_k \in I\}$$

$$= \max\{\|R^m(\bar{x}) - R^m(\bar{x})\|_k : \|\cdot\|_k \in I\}$$

$$= 0.$$

Hence, we proved that for the arbitrarily given $K \in \mathcal{F}_{\mathbb{S}}$ and $0 < \varepsilon < 1$, $\nabla(R^m)(\bar{x})$ satisfies condition (DZ) in Definition 4.8 with respect to the constructed $\tau_{\mathbb{S}}$-open neighborhood $U_{I,\delta}$ of $\theta_{\mathbb{S}}$. Next, for the arbitrarily given $K \in \mathcal{F}_{\mathbb{S}}$ and $0 < \varepsilon < 1$, we prove that $\nabla(R^m)(\bar{x})$ satisfies condition (DR) in Definition 4.8 with respect to the constructed $\tau_{\mathbb{S}}$-open neighborhood $U_{I,\delta}$ of $\theta_{\mathbb{S}}$. Let $u = \{s_n\}_{n=1}^{\infty} \in U_{I,\delta}$. Suppose

$$0 < \max\left\{\frac{|s_k|}{1+|s_k|} : k = 1, 2, \dots, M\right\} = \max\{\|u\|_k : \|\cdot\|_k \in I\} < \delta. \tag{7.9}$$

Since the function $\frac{\lambda}{1+\lambda}$ of $\lambda$ is a strictly increasing function for $\lambda \in [0, \infty)$, then the condition $\delta < \frac{1}{2}$ and (7.9) imply that

$$\max\{|s_k| : \|\cdot\|_k \in I\} < 1.$$

Since $\|\cdot\|_k$ is a $F$-seminorm on $\mathbb{S}$, the above inequality implies that

$$\frac{|s_k|^2}{1+|s_k|^2} \leq \frac{|s_k|}{1+|s_k|}, \text{for } k = 1, 2, \dots, M \tag{7.10}$$

Then, for the given $u = \{s_n\}_{n=1}^{\infty} \in U_{I,\delta}$ satisfying (7.9) and (7.10), by (7.8), we estimate

$$\begin{aligned}
&\|R^m(\bar{x}+u) - R^m(\bar{x}) - \nabla(R^m)(\bar{x})(u)\|_k \\
&\leq \frac{|s_k|^2\left(\binom{m}{2}(\mathcal{P}_M(\bar{x}))^{m-2} + \binom{m}{3}(\mathcal{P}_M(\bar{x}))^{m-3}|s_k| + \cdots + \binom{m}{m}|s_k|^{m-2}\right)}{1+|s_k|^2\left(\binom{m}{2}(\mathcal{P}_M(\bar{x}))^{m-2} + \binom{m}{3}(\mathcal{P}_M(\bar{x}))^{m-3}|s_k| + \cdots + \binom{m}{m}|s_k|^{m-2}\right)} \\
&\leq \frac{|s_k|^2\left(\binom{m}{2}(\mathcal{P}_M(\bar{x}))^{m-2} + \binom{m}{3}(\mathcal{P}_M(\bar{x}))^{m-3} + \cdots + \binom{m}{m}\right)}{1+|s_k|^2\left(\binom{m}{2}(\mathcal{P}_M(\bar{x}))^{m-2} + \binom{m}{3}(\mathcal{P}_M(\bar{x}))^{m-3} + \cdots + \binom{m}{m}\right)} \\
&\leq \frac{|s_k|^2\left(1+\mathcal{P}_M(\bar{x})\right)^m}{1+|s_k|^2\left(1+\mathcal{P}_M(\bar{x})\right)^m} \\
&= \left(1+\mathcal{P}_{\mathrm{M}}(\bar{\mathrm{x}})\right)^{\mathrm{m}} \frac{|\mathrm{s_k}|^2}{1+|\mathrm{s_k}|^2\left(1+\mathrm{q_M}(\bar{\mathrm{x}})\right)^{\mathrm{m}}} \\
&\leq \left(1+\mathcal{P}_M(\bar{x})\right)^m \frac{|s_k|^2}{1+|s_k|^2} \\
&\leq \left(1+\mathcal{P}_M(\bar{x})\right)^m \frac{|s_k|}{1+|s_k|}.
\end{aligned} \tag{7.11}$$

We write $\Delta := \max\left\{\frac{|s_i|}{1+|s_i|} : \|\cdot\|_i \in I\right\}$. By $K \subseteq I$, it follows that

$$\frac{\frac{|s_k|}{1+|s_k|}}{\Delta} \leq 1, \text{ for each } k \text{ with } \|\cdot\|_k \in K. \tag{7.12}$$

Then by (7.7) and (7.11), for $u = \{s_n\}_{n=1}^{\infty} \in U_{I,\delta}$ satisfying (7.9), we have

$$\begin{aligned}
&\max\left\{\left\|\frac{R^m(\bar{x}+u) - R^m(\bar{x}) - \nabla(R^m)(\bar{x})(u)}{\max\left\{\frac{|s_i|}{1+|s_i|} : \|\cdot\|_i \in I\right\}}\right\|_k : \|\cdot\|_k \in K\right\} \\
&= \max\left\{\left\|\frac{R^m(\bar{x}+u) - R^m(\bar{x}) - \nabla(R^m)(\bar{x})(u)}{\max\left\{\frac{|s_i|}{1+|s_i|} : \|\cdot\|_i \in I\right\}}\right\|_k : \|\cdot\|_k \in K\right\}
\end{aligned}$$

$$= \max\left\{ \left\| \left\{ \frac{\binom{m}{2}\bar{t}_n^{m-2}s_n^2 + \binom{m}{3}\bar{t}_n^{m-3}s_n^3 + \cdots + \binom{m}{m}s_n^m}{\Delta} \right\}_{n=1}^{\infty} \right\|_k : \|\cdot\|_k \in K \right\}$$

$$= \max\left\{ \frac{\left| \frac{\binom{m}{2}\bar{t}_k^{m-2}s_k^2 + \binom{m}{3}\bar{t}_k^{m-3}s_k^3 + \cdots + \binom{m}{m}s_k^m}{\Delta} \right|}{1 + \left| \frac{\binom{m}{2}\bar{t}_k^{m-2}s_k^2 + \binom{m}{3}\bar{t}_k^{m-3}s_k^3 + \cdots + \binom{m}{m}s_k^m}{\Delta} \right|} : \|\cdot\|_k \in K \right\}$$

$$\leq \max\left\{ \frac{\frac{|s_k|^2\left(1+\mathcal{P}_M(\bar{x})\right)^m}{\Delta}}{1 + \frac{|s_k|^2\left(1+\mathcal{P}_M(\bar{x})\right)^m}{\Delta}} : \|\cdot\|_k \in K \right\}$$

$$= \left(1+\mathcal{P}_M(\bar{x})\right)^m \max\left\{ \frac{\frac{|s_k|^2}{\Delta}}{1 + \frac{|s_k|^2\left(1+\mathcal{P}_M(\bar{x})\right)^m}{\Delta}} : \|\cdot\|_k \in K \right\}$$

$$\leq \left(1+\mathcal{P}_M(\bar{x})\right)^m \max\left\{ \frac{\frac{|s_k|^2}{\Delta}}{1 + \frac{|s_k|^2}{\Delta}} : \|\cdot\|_k \in K \right\}$$

$$= \left(1+\mathcal{P}_M(\bar{x})\right)^m \max\left\{ \frac{\frac{\frac{|s_k|}{1+|s_k|}|s_k|(1+|s_k|)}{\Delta}}{1 + \frac{\frac{|s_k|}{1+|s_k|}|s_k|(1+|s_k|)}{\Delta}} : \|\cdot\|_k \in K \right\}$$

$$\leq \left(1+\mathcal{P}_M(\bar{x})\right)^m \max\left\{ \frac{|s_k|(1+|s_k|)}{1+|s_k|(1+|s_k|)} : \|\cdot\|_k \in K \right\} \qquad \text{(by (7.12))}$$

$$\leq 2\left(1+\mathcal{P}_M(\bar{x})\right)^m \max\left\{ \frac{|s_k|}{1+|s_k|(1+|s_k|)} : \|\cdot\|_k \in K \right\}$$

$$\leq 2\left(1+\mathcal{P}_M(\bar{x})\right)^m \max\left\{ \frac{|s_k|}{1+|s_k|} : \|\cdot\|_k \in K \right\}$$

$$\leq 2\left(1+\mathcal{P}_M(\bar{x})\right)^m \max\left\{ \frac{|s_i|}{1+|s_i|} : \|\cdot\|_i \in I \right\}$$

$$\leq 2\left(1+\mathcal{P}_M(\bar{x})\right)^m \delta$$

$$< \varepsilon.$$

Hence, we proved that for the arbitrarily given $\tau_{\mathbb{S}}$-open neighborhood $W_{K,\varepsilon}$ of $\theta_{\mathbb{S}}$, $\nabla(R^m)(\bar{x})$ satisfies condition (i) in Definition 4.8 with respect to the constructed $\tau_{\mathbb{S}}$-open neighborhood $U_{I,\delta}$ of $\theta_{\mathbb{S}}$. □

Similarly to Theorem 6.7, the space $(\mathbb{S}, \tau_{\mathbb{S}})$ is a Hausdorff topological vector space and it is not seminorm constructed. We need to prove directly the Gâteaux differentiability of polynomial type operators $R^m$ in the space $(\mathbb{S}, \tau_{\mathbb{S}})$ without using its Fréchet differentiability. Meanwhile, similarly to Theorem 6.7, we will find that the Gâteaux derivatives of polynomial type operators $R^m$ coincide with their Fréchet derivatives.

**Theorem 7.3**. *Let m be a positive integer. Then, the* $m^{\text{th}}$ *power operator* $R^m$ *in* $(\mathbb{S}, \tau_{\mathbb{S}})$ *is Gâteaux differentiable on* $\mathbb{S}$ *such that for any* $\bar{x} = (\bar{t}_1, \bar{t}_2, \ldots) \in \mathbb{S}$, *the Gâteaux derivative* $(R^m)'(\bar{x})$ *of* $R^m$ *at* $\bar{x}$

*satisfies that,*

(i) *For $m > 1$, we have*

$$(R^m)'(\bar{x}) = \nabla(R^m)(\bar{x}) = \begin{pmatrix} m\bar{t}_1^{m-1} & 0 & \cdots \\ 0 & m\bar{t}_2^{m-1} & 0 \\ \vdots & 0 & \ddots \end{pmatrix}.$$

*Now, for every $v = \{s_n\}_{n=1}^{\infty} \in \mathbb{S}\backslash\{\theta_{\mathbb{S}}\}$*

$$(R^m)'(\bar{x})(v) = \nabla(R^m)(\bar{x})(v) = \{s_n\}_{n=1}^{\infty} \begin{pmatrix} m\bar{t}_1^{m-1} & 0 & \cdots \\ 0 & m\bar{t}_2^{m-1} & 0 \\ \vdots & 0 & \ddots \end{pmatrix} = \{m\bar{t}_n^{m-1}s_n\}_{n=1}^{\infty} \in \mathbb{S}.$$

(ii) *When $m = 1$, we have*

$$R'(\bar{x}) = \nabla(R)(\bar{x}) = \begin{pmatrix} 1 & 0 & \cdots \\ 0 & 1 & 0 \\ \vdots & 0 & \ddots \end{pmatrix}.$$

*It is defined as a constant pointwise multiplication operator on $\mathbb{S}$ such that*

$$R'(\bar{x})(v) = v = \{s_n\}_{n=1}^{\infty}, \textit{for every } v = \{s_n\}_{n=1}^{\infty} \in \mathbb{S}\backslash\{\theta_{\mathbb{S}}\}.$$

*Proof*. Part (ii) is clear. So, we only prove (i). Let $\bar{x} = (\bar{t}_1, \bar{t}_2, \dots) \in \mathbb{S}\backslash\{\theta_{\mathbb{S}}\}$. Let $K \in \mathcal{F}_{\mathbb{S}}$ be arbitrarily given. Let $v = \{s_n\}_{n=1}^{\infty} \in \mathbb{S}\backslash\{\theta_{\mathbb{S}}\}$ be arbitrarily given. For any real number $t \neq 0$, similarly to (7.7) and (7.8), we calculate

$$\max\left\{\left\|\frac{R^m(\bar{x}+tv)-R^m(\bar{x})-\nabla(R^m)(\bar{x})(tv)}{t}\right\|_k : \|\cdot\|_k \in K\right\}$$

$$= \max\left\{\left\|\frac{R^m(\bar{x}+tv)-R^m(\bar{x})-t\nabla(R^m)(\bar{x})(v)}{t}\right\|_k : \|\cdot\|_k \in K\right\}$$

$$= \max\left\{\left\|\left\{\frac{\binom{m}{2}\bar{t}_n^{m-2}t^2s_n^2+\binom{m}{3}\bar{t}_n^{m-3}t^3s_n^3+\cdots+\binom{m}{m}t^m s_n^m}{t}\right\}_{n=1}^{\infty}\right\|_k : \|\cdot\|_k \in K\right\}$$

$$= \max\left\{\frac{\left|\frac{\binom{m}{2}\bar{t}_k^{m-2}t^2s_k^2+\binom{m}{3}\bar{t}_k^{m-3}t^3s_k^3+\cdots+\binom{m}{m}t^m s_k^m}{t}\right|}{1+\left|\frac{\binom{m}{2}\bar{t}_k^{m-2}t^2s_k^2+\binom{m}{3}\bar{t}_k^{m-3}t^3s_k^3+\cdots+\binom{m}{m}t^m s_k^m}{t}\right|} : \|\cdot\|_k \in K\right\}$$

$$= t\max\left\{\frac{|s_k|^2\left|\binom{m}{2}\bar{t}_k^{m-2}+\binom{m}{3}\bar{t}_k^{m-3}ts_n+\cdots+\binom{m}{m}t^{m-2}s_k^{m-2}\right|}{1+t|s_k|^2\left|\binom{m}{2}\bar{t}_k^{m-2}+\binom{m}{3}\bar{t}_k^{m-3}ts_n+\cdots+\binom{m}{m}t^{m-2}s_k^{m-2}\right|} : \|\cdot\|_k \in K\right\}$$

$$\to 0, \text{ as } t \to 0.$$ □

Here, for the given and fixed $\bar{x} = (\bar{t}_1, \bar{t}_2, \dots) \in \mathbb{S}$ and $v = \{s_n\}_{n=1}^{\infty} \in \mathbb{S}\backslash\{\theta_{\mathbb{S}}\}$, and with the finite set $K$, the factor before the above limit satisfies

$$\max\left\{\frac{|s_k|^2\left|\binom{m}{2}\bar{t}_k^{m-2}+\binom{m}{3}\bar{t}_k^{m-3}ts_n+\cdots+\binom{m}{m}t^{m-2}s_k^{m-2}\right|}{1+t|s_k|^2\left|\binom{m}{2}\bar{t}_k^{m-2}+\binom{m}{3}\bar{t}_k^{m-3}ts_n+\cdots+\binom{m}{m}t^{m-2}s_k^{m-2}\right|} : \|\cdot\|_k \in K\right\} < \infty, \text{ for } |t| < 1.$$ □

To finish this section, we reconsider the space $\mathbb{S}$ equipped with a different topology. Define $d_{\mathbb{S}}: \mathbb{S} \times \mathbb{S} \to \mathbb{R}_+$ by

$$d_{\mathbb{S}}(x, y) = \sum_{n=1}^{\infty} \frac{1}{2^n} \frac{|t_n - s_n|}{1+|t_n - s_n|}, \text{ for } x = \{t_n\}_{n=1}^{\infty},\ y = \{s_n\}_{n=1}^{\infty} \in \mathbb{S}. \tag{7.13}$$

By this definition, $d_{\mathbb{S}}$ is a translation invariant metric on $\mathbb{S}$. It is known that the metric space $(\mathbb{S}, d_{\mathbb{S}})$ is not locally convex. Based on this translation invariant metric $d_{\mathbb{S}}$ on $\mathbb{S}$, define $p_{\mathbb{S}}: \mathbb{S} \to \mathbb{R}_+$ by

$$p_{\mathbb{S}}(x) = d(x, \theta) = \sum_{n=1}^{\infty} \frac{1}{2^n} \frac{|t_n|}{1+|t_n|}, \text{ for } x = \{t_n\}_{n=1}^{\infty} \in \mathbb{S}.$$

Then $p_{\mathbb{S}}: \mathbb{S} \to \mathbb{R}_+$ is an $F$-norm. Let $\tau_{\mathbb{S}}$ be the topology on $\mathbb{S}$ induced by the $F$-norm $p_{\mathbb{S}}$. With respect to this topology $\tau_{\mathbb{S}}$, $(\mathbb{S}, \tau_{\mathbb{S}})$ is an $F$-space, in which the topology $p_{\mathbb{S}}$ on $\mathbb{S}$ is induced by one $F$-norm $p_{\mathbb{S}}$.

## 8. Single-valued Mappings from $(\sigma_\rho, \tau_\rho)$ to $(\mathbb{S}, \tau_{\mathbb{S}})$

In this paper, Gâteaux (directional) differentiability (Definition 4.1) and Fréchet differentiability (Definition 4.8) of mappings are defined, in which the considered mappings are from one general topological vector space to another one. In these definitions, the domain spaces may be different from the range spaces. However, in sections 5, 6 and 7, we only consider the differentiability of operators from a general topological vector space to itself. In this section, we study the differentiability of operators from a general topological vector space $(\sigma_\rho, \tau_\rho)$ to $(\mathbb{S}, \tau_{\mathbb{S}})$, which is a different topological vector space, both of them are not seminorm-constructed.

Let $(\sigma_\rho, \tau_\rho)$ be the not seminorm-constructed topological vector space with given $0 < \rho < 1$ studied in Section 6, in which the topology $\tau_\rho$ is induced by the countable family $\mathbb{F}_{\sigma_\rho}$ of $F$-seminorms $\|\cdot\|_{\rho,k}: \sigma_\rho \to \mathbb{R}_+$, where $k = 1, 2, \ldots$ defined by

$$\|x\|_{\rho,k} = |t_k|^\rho, \text{ for any } x = \{t_n\}_{n=1}^{\infty} \in \sigma_\rho.$$

Let $(\mathbb{S}, \tau_{\mathbb{S}})$ be the not seminorm-constructed topological vector space, in which the topology $\tau_{\mathbb{S}}$ is induced by the countable family $\mathbb{F}_{\mathbb{S}}$ of $F$-seminorms $\|\cdot\|_k: \mathbb{S} \to \mathbb{R}_+$, where $k = 1, 2, \ldots$ and defined by

$$\|x\|_k = \frac{|t_k|}{1+|t_k|}, \text{ for any } x = \{t_n\}_{n=1}^{\infty} \in \mathbb{S}.$$

Let $m$ be a positive integer. We define the $m^{\text{th}}$ power operator $\Lambda^m: \sigma_\rho \to \mathbb{S}$ as follows.

$$\Lambda^m(x) = \{t_n^m\}_{n=1}^{\infty} \in \mathbb{S}, \text{ for any } x = \{t_n\}_{n=1}^{\infty} \in \sigma_\rho. \tag{8.1}$$

**Theorem 8.1**. *Let m be a positive integer. Then* $\Lambda^m: \sigma_\rho \to \mathbb{S}$ *is Fréchet differentiable on* $\sigma_\rho$ *and for any* $\bar{x} = (\bar{t}_1, \bar{t}_2, \ldots) \in \sigma_\rho$, $\nabla(\Lambda^m)(\bar{x})$ *satisfies that,*

(i) *For m* > 1, *we have*

$$\nabla(\Lambda^m)(\bar{x}) = \begin{pmatrix} m\bar{t}_1^{m-1} & 0 & \cdots \\ 0 & m\bar{t}_2^{m-1} & 0 \\ \vdots & 0 & \ddots \end{pmatrix}.$$

*Here*, $\nabla(\Lambda^m)(\bar{x})$ *is an* $\infty \times \infty$ *matrix that defines a pointwise multiplication operator on* $\sigma_\rho$ *satisfying that, for every* $x = \{t_n\}_{n=1}^{\infty} \in \sigma_\rho$, *we have*

$$\nabla(\Lambda^m)(\bar{x})(x) = \{t_n\}_{n=1}^{\infty}\begin{pmatrix} m\bar{t}_1^{m-1} & 0 & \dots \\ 0 & m\bar{t}_2^{m-1} & 0 \\ \vdots & 0 & \ddots \end{pmatrix} = \{m\bar{t}_n^{m-1}t_n\}_{n=1}^{\infty} \in \mathbb{S}. \tag{8.2}$$

(ii) *In particular, when* $m = 1$, *we have*

$$\nabla(\Lambda^1)(\bar{x}) = \begin{pmatrix} 1 & 0 & \dots \\ 0 & 1 & 0 \\ \vdots & 0 & \ddots \end{pmatrix}.$$

*It is defined as a constant pointwise multiplication operator on* $\sigma_\rho$ *such that*

$$\nabla(\Lambda^1)(\bar{x})(x) = x = \{t_n\}_{n=1}^{\infty}, \textit{for every } x = \{t_n\}_{n=1}^{\infty} \in \mathbb{S}.$$

*Proof*. It is clear that this theorem holds for $m = 1$. Hence, we only prove part (i) of this theorem for $m > 1$ by proving (8.2). Let $K \in \mathcal{F}_{\mathbb{S}}$ and $0 < \varepsilon < 1$ be arbitrarily given, which induce an $\tau_{\mathbb{S}}$-open neighborhood $W_{K,\varepsilon}$ of $\theta_{\mathbb{S}}$ in $\mathbb{S}$. This implies that there is a positive integer $M$ such that

$$K \subseteq \{\|\cdot\|_k : k = 1, 2, \dots, M\}.$$

For the fixed positive integer $M$, similarly to the previous section, we define

$$\mathcal{P}_M(x) = \max\{|t_k| : k = 1, 2, \dots, M\}, \text{ for any } x = (t_1, t_2, \dots) \in \mathbb{S}.$$

It is clear that, for any $\bar{x} = (\bar{t}_1, \bar{t}_2, \dots) \in \mathbb{S}$, we have $\mathcal{P}_M(\bar{x}) < \infty$. Let

$$I = \left\{\|\cdot\|_{\rho,k} : k = 1, 2, \dots, M\right\} \in \mathcal{F}_{\sigma_\rho} \quad \text{and} \quad \delta = \frac{\varepsilon}{\left(1+\mathcal{P}_M(\bar{x})\right)^m}.$$

With the above $I \in \mathcal{F}_{\sigma_\rho}$ and $0 < \delta < 1$, we have an $\tau_\rho$-open neighborhood $U_{I,\delta}$ of $\theta_{\sigma_\rho}$ in $\sigma_\rho$. Then, for the arbitrarily given $\tau_{\mathbb{S}}$-open neighborhood $W_{K,\varepsilon}$ of $\theta_{\mathbb{S}}$, we first prove that $\nabla(\Lambda^m)(\bar{x})$ satisfies condition (DZ) in Definition 4.8 with respect to the constructed $\tau_\rho$-open neighborhood $U_{I,\delta}$ of $\theta_{\sigma_\rho}$ in $\sigma_\rho$.

Let $u = \{s_n\}_{n=1}^{\infty} \in U_{I,\delta}$. Suppose that

$$\max\{|s_k|^\rho : k = 1, 2, \dots, M\} = \max\left\{\|u\|_{\rho,k} : \|\cdot\|_{\rho,k} \in I\right\} = 0.$$

This implies that $s_k = 0$, for $k = 1, 2, \dots, M$. Then, we obtain that

$$\begin{aligned}
&\max\{\|\Lambda^m(\bar{x}+u) - \Lambda^m(\bar{x}) - \nabla(\Lambda^m)(\bar{x})(u)\|_k : \|\cdot\|_k \in K\} \\
&\leq \max\{\|\Lambda^m(\bar{x}+u) - \Lambda^m(\bar{x}) - \nabla(\Lambda^m)(\bar{x})(u)\|_k : k = 1, 2, \dots, M\} \\
&= \max\{\|\Lambda^m(\bar{x}) - \Lambda^m(\bar{x})\|_k : k = 1, 2, \dots, M\} \\
&= 0.
\end{aligned}$$

Hence, we proved that for the arbitrarily given $\tau_{\mathbb{S}}$-open neighborhood $W_{K,\varepsilon}$ of $\theta_{\mathbb{S}}$, $\nabla(\Lambda^m)(\bar{x})$ satisfies condition (DZ) in Definition 4.8 with respect to the constructed $\tau_\rho$-open neighborhood $U_{I,\delta}$ of $\theta_{\sigma_\rho}$ in $\sigma_\rho$. Next, for the arbitrarily given $\tau_{\mathbb{S}}$-open neighborhood $W_{K,\varepsilon}$ of $\theta_{\mathbb{S}}$ in $\mathbb{S}$, we prove that $\nabla(\Lambda^m)(\bar{x})$ satisfies condition (DR) in Definition 4.8 with respect to the constructed $\tau_\rho$-open neighborhood $U_{I,\delta}$ of $\theta_{\sigma_\rho}$ in $\sigma_\rho$.

Let $u = \{s_n\}_{n=1}^{\infty} \in U_{I,\delta}$. Suppose that

$$0 < \max\{|s_i|^\rho : i = 1, 2, \dots, M\} = \max\{\|u\|_{\rho,i} : \|\cdot\|_{\rho,i} \in I\} < \delta. \tag{8.3}$$

Since for the given fixed $0 < \rho < 1$, the function $\lambda^\rho$ of $\lambda$ is a strictly increasing function for $\lambda \in [0, \infty)$, then the conditions $\delta < \frac{1}{2}, \frac{1}{\rho} > 1$ and (8.3) imply that

$$\max\{|s_i|^\rho : \|\cdot\|_{\rho,i} \in I\} < \delta^{\frac{1}{\rho}} < \delta < 1. \tag{8.4}$$

Since $\|\cdot\|_k$ is a $F$-seminorm on $\mathbb{S}$, the above inequality implies that $\|s_k u\|_k \leq \|u\|_k$, that is

$$\frac{|s_k|^2}{1+|s_k|^2} \leq \frac{|s_k|}{1+|s_k|}, \text{for } k = 1, 2, \dots, M \tag{8.5}$$

Then, for this given $u = \{s_n\}_{n=1}^{\infty} \in U_{I,\delta}$ satisfying (8.3)–(8.5), by the definition of $\mathcal{P}_M$, we estimate

$$\max\left\{\left\|\frac{\Lambda^m(\bar{x}+u)-\Lambda^m(\bar{x})-\nabla(\Lambda^m)(\bar{x})(u)}{\max\{\|u\|_{\rho,i}:\|\cdot\|_{\rho,i}\in I\}}\right\|_k : \|\cdot\|_k \in K\right\}$$

$$= \max\left\{\left\|\frac{\Lambda^m(\bar{x}+u)-\Lambda^m(\bar{x})-\nabla(\Lambda^m)(\bar{x})(u)}{\max\{|s_i|^\rho:\|\cdot\|_{\rho,i}\in I\}}\right\|_k : \|\cdot\|_k \in K\right\}$$

$$= \max\left\{\left\|\frac{\{\binom{m}{2}\bar{t}_n^{m-2}s_n^2+\binom{m}{3}\bar{t}_n^{m-3}s_n^3+\cdots+\binom{m}{m}s_n^m\}_{n=1}^{\infty}}{\max\{|s_i|^\rho:\|\cdot\|_{\rho,i}\in I\}}\right\|_k : \|\cdot\|_k \in K\right\}$$

$$= \max\left\{\frac{\left|\frac{\binom{m}{2}\bar{t}_k^{m-2}s_k^2+\binom{m}{3}\bar{t}_k^{m-3}s_k^3+\cdots+\binom{m}{m}s_k^m}{\max\{|s_i|^\rho:\|\cdot\|_{\rho,i}\in I\}}\right|}{1+\left|\frac{\binom{m}{2}\bar{t}_k^{m-2}s_k^2+\binom{m}{3}\bar{t}_k^{m-3}s_k^3+\cdots+\binom{m}{m}s_k^m}{\max\{|s_i|^\rho:\|\cdot\|_{\rho,i}\in I\}}\right|} : \|\cdot\|_k \in K\right\}$$

$$= \max\left\{\frac{|s_k|^2\left|\frac{\left|\binom{m}{2}\bar{t}_k^{m-2}+\binom{m}{3}\bar{t}_k^{m-3}s_k^1+\cdots+\binom{m}{m}s_k^{m-2}\right|}{\max\{|s_i|^\rho:\|\cdot\|_{\rho,i}\in I\}}\right|}{1+|s_k|^2\left|\frac{\left|\binom{m}{2}\bar{t}_k^{m-2}+\binom{m}{3}\bar{t}_k^{m-3}s_k^1+\cdots+\binom{m}{m}s_k^{m-2}\right|}{\max\{|s_i|^\rho:\|\cdot\|_{\rho,i}\in I\}}\right|} : \|\cdot\|_k \in K\right\}$$

$$= \max\left\{\frac{|s_k|^2\left|\frac{\left|\binom{m}{2}\bar{t}_k^{m-2}+\binom{m}{3}\bar{t}_k^{m-3}s_k^1+\cdots+\binom{m}{m}s_k^{m-2}\right|}{\max\{|s_i|^\rho:i=1,\dots,M\}}\right|}{1+|s_k|^2\left|\frac{\left|\binom{m}{2}\bar{t}_k^{m-2}+\binom{m}{3}\bar{t}_k^{m-3}s_k^1+\cdots+\binom{m}{m}s_k^{m-2}\right|}{\max\{|s_i|^\rho:i=1,\dots,M\}}\right|} : \|\cdot\|_k \in K\right\}$$

$$\leq \max\left\{\frac{\frac{|s_k|^2(1+\mathcal{P}_M(\bar{x}))^m}{\max\{|s_i|^\rho:i=1,\dots,M\}}}{1+\frac{|s_k|^2(1+\mathcal{P}_M(\bar{x}))^m}{\max\{|s_i|^\rho:i=1,\dots,M\}}} : \|\cdot\|_k \in K\right\}$$

$$= \left(1 + \mathcal{P}_M(\bar{x})\right)^m \max \left\{ \frac{\frac{|s_k|^2}{\max\{|s_i|^\rho : i=1,\ldots,M\}}}{1 + \frac{|s_k|^2 \left(1+\mathcal{P}_M(\bar{x})\right)^m}{\max\{|s_i|^\rho : i=1,\ldots,M\}}} : \|\cdot\|_k \in K \right\}$$

$$\leq \left(1 + \mathcal{P}_M(\bar{x})\right)^m \max \left\{ \frac{|s_k|^{2-\rho}}{1 + \frac{|s_k|^2}{\max\{|s_i|^\rho : i=1,\ldots,M\}}} : \|\cdot\|_k \in K \right\}$$

$$= \left(1 + \mathcal{P}_M(\bar{x})\right)^m \max \{|s_k|^{2-\rho} : \|\cdot\|_k \in K\}$$

$$= \left(1 + \mathcal{P}_M(\bar{x})\right)^m \left(\max \{|s_k|^{\rho} : \|\cdot\|_k \in K\}\right)^{\frac{2-\rho}{\rho}}$$

$$\leq \left(1 + \mathcal{P}_M(\bar{x})\right)^m \left(\max \{|s_k|^{\rho} : k = 1, \ldots, M\}\right)^{\frac{2-\rho}{\rho}}$$

$$\leq \left(1 + \mathcal{P}_M(\bar{x})\right)^m \delta^{\frac{2-\rho}{\rho}}$$

$$< \left(1 + \mathcal{P}_M(\bar{x})\right)^m \delta$$

$$< \varepsilon.$$

Hence, we proved that for the arbitrary $\tau_{\mathbb{S}}$-open neighborhood $W_{K,\varepsilon}$ of $\theta_{\mathbb{S}}$, $\nabla(\Lambda^m)(\bar{x})$ satisfies condition (DR) in Definition 4.8 with respect to the constructed $\tau_\rho$-open neighborhood $U_{I,\delta}$ of $\theta_{\sigma_\rho}$ in $\sigma_\rho$. □

## 9. Some Applications of Differentiation to Vector Optimizations in Partially Ordered Topological Vector Spaces

### 9.1. Review Partially Ordered Topological Vector Spaces

Throughout this section, if not otherwise stated, we let $(X, \tau_X)$ and $(Y, \tau_Y)$ be Hausdorff topological vector spaces, in which the topologies $\tau_X$ and $\tau_Y$ are induced by families $\mathbb{F}_X$ and $\mathbb{F}_Y$ of positive $F$-seminorms, respectively, which are as studied in sections 2, 3 and 4. Recall that, by (2.2) to (2.5), for any $x_0 \in X$, $I \in \mathcal{F}_X$ and $\delta > 0$, we write

$$U_{I,\delta}(x_0) = \{x \in X : \max\{p(x - x_0) : p \in I\} < \delta\}$$

and
$$U(\mathbb{F}_X)(x_0) = \{U_{I,\lambda}(x_0) : I \in \mathcal{F}_X, \lambda > 0\}.$$

$U(\mathbb{F}_X)(x_0)$ forms a basis of $\tau_X$ around point $x_0$. Let $y_0 \in Y$, $J \in \mathcal{F}_Y$ and $\lambda > 0$. Then we write

$$V_{J,\lambda}(y_0) = \{y \in Y : \max\{q(y - y_0) : q \in J\} < \lambda\}$$

and
$$V(\mathbb{F}_Y)(Y_0) = \{V_{J,\lambda}(y_0) : J \in \mathcal{F}_Y, \lambda > 0\}. \tag{9.1}$$

$V(\mathbb{F}_Y)(Y_0)$ forms a basis of $\tau_Y$ around point $y_0$. In this section, we first review the concepts of partial orders on topological vector spaces induced by closed, convex and pointed cones. By these partial orders, we define the order maximum and order minimum values along given directions of single-valued mappings. Then, in contrast to optimization properties in calculus, we investigate the connection between vector optimization values and Gâteaux and Fréchet derivatives with respect to some single-valued mappings.

Let $C$ be a nonempty closed convex and pointed cone in $X$ with $C \neq \{\theta_X\}$. Let $\preccurlyeq_C$ be the partial order on $X$ generated by $C$ such that, for any $x_1, x_2 \in X$,

$$x_1 \preccurlyeq_C x_2 \qquad \Leftrightarrow \qquad x_2 - x_1 \in C.$$

More strictly, for any $x_1, x_2 \in X$,

$$x_1 \prec_C x_2 \qquad \Leftrightarrow \qquad x_2 - x_1 \in C\backslash\{\theta_X\}.$$

With this partial order $\preccurlyeq_C$, $(X, \tau_X, \preccurlyeq_C)$ is called a partially ordered Hausdorff topological vector space, in which the topology $\tau_X$ is induced by a family $\mathbb{F}_X$ of positive $F$-seminorms on $X$ and the partial order $\preccurlyeq_C$ is defined by a nonempty closed convex and pointed cone $C$ in $X$.

Similarly to the partial order $\preccurlyeq_C$ on $X$, we define a partial order on $Y$. Let $K$ be a nonempty closed convex and pointed cone in $Y$ with $K \neq \{\theta_Y\}$. Similarly to the definition of the partial order $\preccurlyeq_C$ on $X$, we define the partial order $\preccurlyeq_K$ on $Y$ generated by $K$. Then, with this partial order $\preccurlyeq_K$ on $Y$, $(Y, \tau_Y, \preccurlyeq_K)$ is also a partially ordered Hausdorff topological vector space, in which the topology $\tau_Y$ is induced by a family $\mathbb{F}_Y$ of positive $F$-seminorms on $Y$ and the partial order $\preccurlyeq_K$ is defined by a nonempty closed convex and pointed cone $K$ in $Y$.

### 9.2. Vector Optimization in Partially Ordered Topological Vector Spaces

In calculus, for a given differentiable single valued function, a critical point of this function is a necessary condition for this considered function to take a (local or absolute) extreme value at it. In this subsection, we extend this property to single-valued mappings with respect to Gâteaux differentiability and Fréchet differentiability and ordered extrema in partially ordered topological vector spaces. We start to define generalized critical points and ordered extrema of single-valued mappings.

Throughout this subsection, if not otherwise stated, we always let $(X, \tau_X)$ be a Hausdorff topological vector space. Let $(Y, \tau_Y, \preccurlyeq_K)$ be a partially ordered Hausdorff topological vector space. Let $A$ be a nonempty convex $\tau_X$-open subset in $X$. Let $T: A \to Y$ be a single-valued mapping and let $\bar{x} \in A$.

**Definition 9.1.** Suppose that $T$ is Gâteaux differentiable at $\bar{x}$. If

$$T'(\bar{x})(v) = \theta_Y, \text{ for every } v \in X\backslash\{\theta_X\},$$

then $\bar{x}$ is called a generalized critical point of $T$.

**Definition 9.2.** Let $\bar{x} \in A$ and let $v \in X\backslash\{\theta_X\}$. If

$$T(\bar{x} + tv) \preccurlyeq_K T(\bar{x}), \text{ for all real } t \text{ with } \bar{x} + tv \in A, \tag{9.2}$$

then, $T$ is said to take a directional $\preccurlyeq_K$-maximum value in $A$ at $\bar{x}$ along direction $v$. $T(\bar{x})$ is called the directional $\preccurlyeq_K$-maximum value of $T$ in $A$ along direction $v$. The directional $\preccurlyeq_K$-minimum value of $T$ at $\bar{x}$ along direction $v \in X\backslash\{\theta_X\}$ can be similarly defined.

More strictly speaking, if $T$ satisfies the following order-inequality

$$T(x) \preccurlyeq_K T(\bar{x}), \text{ for all } x \in A, \tag{9.3}$$

then, $T$ is said to take (absolute) $\preccurlyeq_K$-maximum value in $A$ at $\bar{x}$ and $T(\bar{x})$ is called the (absolute) $\preccurlyeq_K$-maximum value of $T$ in $A$. The (absolute) $\preccurlyeq_K$-minimum value of $T$ in $A$ at $\bar{x}$ can be similarly defined. In either of (absolute) $\preccurlyeq_K$-maximum value or (absolute) $\preccurlyeq_K$-minimum value, they are called $\preccurlyeq_K$-extreme values. In this case, we say that $\bar{x}$ is an $\preccurlyeq_K$-extreme point of the mapping $T$.

Suppose that $T$ is Gâteaux differentiable at $\bar{x}$. More precisely, in Corollary 9.4 below, we prove

$$\bar{x} \text{ is an } \preccurlyeq_K\text{-extreme point of } T \quad \Longrightarrow \quad \bar{x} \text{ is an generalized critical point of } T.$$

In the counter examples 9.10 and 9.11 below, we show that

$$\bar{x} \text{ is an generalized critical point of } T \quad \not\Rightarrow \quad \bar{x} \text{ is an } \preccurlyeq_K\text{-extreme point of } T.$$

In next theorem, we investigate the connection between generalized critical points and ordered extreme points of single-valued mappings in partially ordered Hausdorff topological vector spaces, which naturally generalizes the connection between critical points and extreme points of real valued functions in calculus.

**Theorem 9.3**. *Let $\bar{x} \in A$ and $v \in X\backslash\{\theta_X\}$. If $T$ takes a directional $\preccurlyeq_K$-maximum ($\preccurlyeq_K$-minimum) value at $\bar{x}$ in direction $v$ and $T$ is Gâteaux directionally differentiable at $\bar{x}$ in direction $v$, then $T'(\bar{x}, v) = \theta_Y$.*

*Proof*. Suppose that $T$ takes its $\preccurlyeq_K$-minimum value at point $\bar{x}$ along direction $v$ (the $\preccurlyeq_K$-maximum case can be similarly proved). At first, we prove $T'(\bar{x}, v) \succcurlyeq_K \theta_Y$, this is, $T'(\bar{x}, v) \in K$. Assume, by the way of contradiction, that $T'(\bar{x}, v) \notin K$, this is equivalent to $T'(\bar{x}, v) \in Y\backslash K$. Notice that $V(\mathbb{F}_Y)(T'(\bar{x}, v))$ forms an $\tau_Y$-open neighborhood basis of $Y$ around $T'(\bar{x}, v)$, here

$$V(\mathbb{F}_Y)(T'(\bar{x}, v)) = \{V_{J,\lambda}(T'(\bar{x}, v)) : J \in \mathcal{F}_Y, \lambda > 0\}.$$

Since $Y\backslash K$ is an $\tau_Y$-open subset of $Y$, by $T'(\bar{x}, v) \in Y\backslash K$, there is an $\tau_Y$-open neighborhood (It is defined in (9.1)) $V_{J,\varepsilon}(T'(\bar{x}, v))$ of $T'(\bar{x}, v)$ in $Y$, for some $J \in \mathcal{F}_Y$ and $\varepsilon > 0$ such that $V_{J,\varepsilon}(T'(\bar{x}, v)) \subseteq Y\backslash K$. This induces

$$\left(V_{J,\varepsilon}(T'(\bar{x}, v))\right) \cap K = \emptyset. \tag{9.4}$$

Notice that by definition (9.1), we have

$$V_{J,\varepsilon}(T'(\bar{x}, v)) = \{y \in Y : \max\{q(y - T'(\bar{x}, v)) : q \in J\} < \varepsilon\}. \tag{9.5}$$

(9.4) and (9.5) together imply that

$$\min\{q(y - T'(\bar{x}, v)) : q \in J\} \geq \varepsilon, \text{ for any } y \in K. \tag{9.6}$$

Since $T$ is Gâteaux differentiable at point $\bar{x}$ along direction $v$, by (4.1) in Definition 4.1, for the above given $\varepsilon$-open neighborhood $V_{J,\varepsilon}(T'(\bar{x}, v))$ of $T'(\bar{x}, v)$ in $Y$, there is $\delta_1 > 0$ such that, for real number $t$, we have

$$0 < |t| < \delta_1 \quad \Longrightarrow \quad \text{Max}\left\{q\left(\frac{T(\bar{x}+tv)-T(\bar{x})}{t} - T'(\bar{x}, v)\right) : q \in J\right\} < \varepsilon. \tag{9.7}$$

By the conditions in this theorem that $A$ is a convex $\tau_X$-open subset in $X$ and $\bar{x} \in A$, and by the property that $U(\mathbb{F}_X)(\bar{x})$ forms an $\tau_X$-open neighborhood basis of $X$ around $\bar{x}$, there is an $\tau_X$-open neighborhood $U_{I,\delta_2}(\bar{x})$ of $\bar{x}$ in $X$, for some $I \in \mathcal{F}_X$ and $\delta_2 > 0$ such that $U_{I,\delta_2}(\bar{x}) \subseteq A$ satisfying

$$U_{I,\delta_2}(\bar{x}) = \{x \in X: \max\{p(x-\bar{x}): p \in I\} < \delta_2\}. \tag{9.8}$$

By the properties (ii) and (iii) in Lemma 2.5 of $F$-seminorms on $X$, for this $v \in X\backslash\{\theta_X\}$, we have

$$\max\{p(tv): p \in I\} \to 0 \text{ as } t \to 0. \tag{9.9}$$

Then, for the above fixed $\delta_2 > 0$, by (9.9), there is $\delta_3 > 0$ with $0 < \delta_3 < \delta_2$, such that

$$|t| < \delta_3 \quad \Longrightarrow \max\{p(tv): p \in I\} < \delta_2 \Longrightarrow \bar{x} + tv \in U_{I,\delta_2}(\bar{x}) \subseteq A. \tag{9.10}$$

(9.10) implies that

$$|t| < \delta_3 \quad \Longrightarrow \quad T(\bar{x} + tv) \text{ is well-defined.} \tag{9.11}$$

Let $\delta = \min\{\delta_1, \delta_3\}$. By (9.11) and (9.7) we obtain that

$$0 < |t| < \delta \quad \Longrightarrow \quad \mathrm{Max}\left\{q\left(\frac{T(\bar{x}+tv)-T(\bar{x})}{t} - T'(\bar{x}, v)\right): q \in J\right\} < \varepsilon. \tag{9.12}$$

In particular, let $t > 0$ in (9.12), we get

$$0 < t < \delta \quad \Longrightarrow \quad \mathrm{Max}\left\{q\left(\frac{T(\bar{x}+tv)-T(\bar{x})}{t} - T'(\bar{x}, v)\right): q \in J\right\} < \varepsilon. \tag{9.13}$$

By the assumption that $T(\bar{x} + tv) - T(\bar{x}) \succcurlyeq_K \theta_Y$, it yields that

$$0 < t < \delta \quad \Longrightarrow \quad \frac{T(\bar{x}+tv)-T(\bar{x})}{t} \in K.$$

By $t > 0$ in (9.13), this creates a contradiction between (9.13) and (9.6). Hence $T'(\bar{x}, v) \in K$. This proves

$$T'(\bar{x}, v) \succcurlyeq_K \theta_Y. \tag{9.14}$$

Notice that $-K$ is also a closed convex pointed cone in $Y$, and

$$T'(\bar{x}, v) \preccurlyeq_K \theta_Y \text{ if and only if } T'(\bar{x}, v) \in -K.$$

Similar to the proof of (9.14), when we prove $T'(\bar{x}, v) \preccurlyeq_K \theta_Y$, in particular in (9.12), let $t < 0$, we get

$$0 < -t < \delta \quad \Longrightarrow \quad \mathrm{Max}\left\{q\left(\frac{T(\bar{x}+tv)-T(\bar{x})}{t} - T'(\bar{x}, v)\right): q \in J\right\} < \varepsilon. \tag{9.15}$$

Since $T(\bar{x} + tv) - T(\bar{x}) \succcurlyeq_K \theta_Y$, by $t < 0$ in (9.15), it yields that

$$0 < -t < \delta \quad \Longrightarrow \quad \frac{T(\bar{x}+tv)-T(\bar{x})}{t} \in -K.$$

Then, similar to the proof of (9.14), we can prove

$$T'(\bar{x}, v) \preccurlyeq_K \theta_Y. \tag{9.16}$$

By (9.14) and (9.16), and by the properties of partial orders in topological vector spaces, we obtain that

$T'(\bar{x}, v) = \theta_Y$. □

We have the following consequences of Theorem 9.3 immediately.

**Corollary 9.4**. *Let $\bar{x} \in A$. Suppose that T is Gâteaux differentiable at $\bar{x}$. Then*

(i) *If T takes (absolute) $\preccurlyeq_K$-maximum or (absolute) $\preccurlyeq_K$-minimum value in A at $\bar{x}$, then $\bar{x}$ is an generalized critical point of T with*

$$T'(\bar{x})(v) = \theta_Y, \textit{for every } v \in X\backslash\{\theta_X\}.$$

(ii) *Inversely, if*

$$T'(\bar{x})(v) \neq \theta_Y, \textit{for every } v \in X\backslash\{\theta_X\},$$

*then T takes neither $\preccurlyeq_K$-maximum, nor $\preccurlyeq_K$-minimum value at $\bar{x}$.*

*Proof.* This corollary follows from Theorem 9.3 and the proof is omitted here. □

**Corollary 9.5**. *Let $(X, \|\cdot\|_X)$ be a normed vector space. Let $(Y, \tau_Y, \preccurlyeq_K)$ be a partially ordered Hausdorff topological vector space. Let A be a nonempty convex $\|\cdot\|_X$-open subset in X. Let $T: A \to Y$ be a single-valued mapping. Let $\bar{x} \in A$. Suppose that T is Fréchet differentiable at $\bar{x}$. Then we have*

(i) *If T takes (absolute) $\preccurlyeq_K$-maximum or (absolute) $\preccurlyeq_K$-minimum value in A at $\bar{x}$, then*

$$\nabla T(\bar{x})(v) = T'(\bar{x})(v) = \theta_Y, \textit{for every } v \in X\backslash\{\theta_X\}.$$

(ii) *Inversely, if*

$$\nabla T(\bar{x})(v) \neq \theta_Y, \textit{for every } v \in X\backslash\{\theta_X\},$$

*then T takes neither $\preccurlyeq_K$-maximum, nor $\preccurlyeq_K$-minimum value at point $\bar{x}$.*

*Proof*. By Corollary 4.20, the condition that $(X, \|\cdot\|_X)$ is a normed vector space implies that the Fréchet differentiability of $T$ at $\bar{x}$ induces the Gâteaux differentiability of $T$ at $\bar{x}$, which satisfies $\nabla T(\bar{x}) = T'(\bar{x})$. Then, this corollary follows from Theorem 9.3 and Corollary 9.4 immediately. □

In particular, let $Y = \mathbb{R}$. $(\mathbb{R}, \leq)$ is the special totally ordered Hilbert space, in which the ordinary order $\leq$ is induced by $[0, \infty)$. Let $(Y, \tau_Y, \preccurlyeq_K) = (\mathbb{R}, \leq)$ in Theorem 9.3, Corollary 9.4 and Corollary 9.5, we have the following ordinary maximum and minimum properties.

**Corollary 9.6**. *Let $(X, \tau_X)$ be a Hausdorff topological vector space. Let A be a nonempty convex $\tau_X$-open subset in X and let $T: A \to \mathbb{R}$ be a real valued functional. Let $\bar{x} \in A$. We have that*

(i) *Let $v \in X\backslash\{\theta_X\}$. If T takes a directional maximum (minimum) value at $\bar{x}$ along direction v and T is Gâteaux differentiable at $\bar{x}$ along direction v, then*

$$T'(\bar{x}, v) = 0.$$

(ii) *If T takes a directional maximum (minimum) value in A at $\bar{x}$ and T is Gâteaux differentiable at $\bar{x}$, then*

$$T'(\bar{x})(v) = 0, \textit{for every } v \in X\backslash\{\theta_X\}.$$

(iii) *Suppose that T is Gâteaux differentiable at $\bar{x}$. If T satisfies*

$$T'(\bar{x})(v) \neq 0, \textit{for every } v \in X\backslash\{\theta_X\},$$

*then T takes neither maximum, nor minimum value at point $\bar{x}$.*

*Proof.* This corollary follows from Theorem 9.3 and the proof is omitted here. □

**Corollary 9.7**. *Let $(X, \tau_X)$ be a Hausdorff locally convex topological vector space. Let A be a nonempty convex $\tau_X$-open subset in X. Let $\bar{x} \in A$. Let $T: A \to \mathbb{R}$ be a real valued functional. Suppose that T is Fréchet differentiable at $\bar{x}$. Then*

(i) *If T takes its maximum (minimum) value in A at point $\bar{x}$, then*

$$\nabla T(\bar{x})(v) = T'(\bar{x})(v) = 0, \textit{for every } v \in X\backslash\{\theta_X\}.$$

(ii) *Inversely, if*

$$\nabla T(\bar{x})(v) \neq 0, \textit{for every } v \in X\backslash\{\theta_X\},$$

*then T takes neither maximum, nor minimum value at point $\bar{x}$.*

*Proof.* This corollary follows from Theorem 9.3 and the proof is omitted here. □

We provide two examples below to verify the results of Corollary 9.4:

*T* is Gâteaux differentiable at $\bar{x}$ and *T* takes ordered extreme at $\bar{x}$ and $\bar{x}$ is an generalized critical point of *T*.

**Example 9.8**. In this example, we consider the topological vector space $\mathcal{S}(\mathbb{R}^n, \mathbb{C})$ with origin $\theta_{\mathcal{S}}$ studied in section 5. At first, we introduce a partial order on $\mathcal{S}(\mathbb{R}^n, \mathbb{C})$. Let *K* be the positive cone in $\mathcal{S}(\mathbb{R}^n, \mathbb{C})$ that is defined by

$$K = \{f \in \mathcal{S}(\mathbb{R}^n, \mathbb{C}): f(x) \geq 0, \text{for every } x \in \mathbb{R}^n\}.$$

It is clear that *K* is a nonempty closed convex and pointed cone in $\mathcal{S}(\mathbb{R}^n, \mathbb{C})$ with $K \neq \{\theta_{\mathcal{S}}\}$. Let $\preccurlyeq_K$ be the partial order on $\mathcal{S}(\mathbb{R}^n, \mathbb{C})$ induced by *K*. For an arbitrarily given positive integer *n*, let $P^{2n}$ be the $2n^{\text{th}}$ power operator on $\mathcal{S}(\mathbb{R}^n, \mathbb{C})$ studied in section 5. For any $f \in \mathcal{S}(\mathbb{R}^n, \mathbb{C})$, we have

$$P^{2n}(f)(x) = f(x)^{2n} \geq 0, \text{ for any } x \in \mathbb{R}^n.$$

By the definition of *K*, this implies that

$$P^{2n}(f) \succcurlyeq_K P^{2n}(\theta_{\mathcal{S}}) = \theta_{\mathcal{S}} =, \text{ for any } f \in \mathcal{S}(\mathbb{R}^n, \mathbb{C}).$$

Hence, $P^{2n}$ takes absolute $\preccurlyeq_K$-minimum value at $\theta_{\mathcal{S}}$ on $\mathcal{S}(\mathbb{R}^n, \mathbb{C})$. On the other hand, by Proposition 5.10, $P^{2n}$ is Gâteaux differentiable at $\theta_{\mathcal{S}}$ and the Gâteaux derivative of $P^{2n}$ at $\theta_{\mathcal{S}}$ satisfies that,

$$(P^{2n})'(\theta_{\mathcal{S}})(u) = 2n \cdot 0 \cdot u = \theta_{\mathcal{S}}, \text{ for every } u \in \mathcal{S}(\mathbb{R}^n, \mathbb{C})\backslash\{\theta_{\mathcal{S}}\}.$$

This implies that the point $\theta_{\mathcal{S}}$ is a generalized critical point of $P^{2n}$ that verified Corollary 9.4 in this case.

**Example 9.9**. In this example, we consider the topological vector space $(\mathbb{S}, \tau_{\mathbb{S}})$ with origin $\theta_{\mathbb{S}}$ studied in section 7. Let $K$ be the positive cone in $(\mathbb{S}, \tau_{\mathbb{S}})$ that is defined by

$$K = \{u = \{s_n\}_{n=1}^{\infty} \in (\mathbb{S}, \tau_{\mathbb{S}}): s_n \geq 0, \text{for each } n = 1, 2, \ldots\}.$$

It is clear that $K$ is a nonempty closed convex and pointed cone in $(\mathbb{S}, \tau_{\mathbb{S}})$ with $K \neq \{\theta_{\mathbb{S}}\}$. Let $\preccurlyeq_K$ be the partial order on $(\mathbb{S}, \tau_{\mathbb{S}})$ induced by $K$. For an arbitrarily given positive integer $m$, let $R^{2m}$ be the $2m^{\text{th}}$ power operator on $(\mathbb{S}, \tau_{\mathbb{S}})$ studied in section 7. $R^{2m}$ is defined by

$$R^{2m}(u) = \{s_n^{2m}\}_{n=1}^{\infty} \in (\mathbb{S}, \tau_{\mathbb{S}}), \text{ for any } u = \{s_n\}_{n=1}^{\infty} \in (\mathbb{S}, \tau_{\mathbb{S}}).$$

Similarly to the proof of Example 9.8, we can show that $R^{2m}$ takes absolute $\preccurlyeq_K$-minimum value at $\theta_{\mathbb{S}}$ on $(\mathbb{S}, \tau_{\mathbb{S}})$. And, the point $\theta_{\mathbb{S}}$ is a generalized critical point of $R^{2m}$ that verified Corollary 9.4 in this case.

In contrast with the property that critical points are only necessary conditions to be extreme points of single-valued functions in calculus, we provide some counter examples below to show that in partially ordered topological vector spaces, generalized critical points are only necessary conditions to be ordered extreme points of single-valued mappings. In calculus, it is well-known that

$$\bar{x} \text{ is a critical point of } T \quad \nRightarrow \quad \bar{x} \text{ is an extreme point of } T.$$

By this special case, it immediately implies that, in general

$$\bar{x} \text{ is a generalized critical point of } T \quad \nRightarrow \quad \bar{x} \text{ is an ordered extreme point of } T. \qquad (9.17)$$

However, in addition to the standard case, we provide a counter example in general topological vector spaces to demonstrate (9.17).

**Example 9.10**. Let $\mathcal{S}(\mathbb{R}^n, \mathbb{C})$ be equipped with the partially order $\preccurlyeq_K$ as given in Example 9.8. Let $P^3$ be the 3rd power operator on $\mathcal{S}(\mathbb{R}^n, \mathbb{C})$ studied in section 5. By Proposition 5.10, $P^3$ is Gâteaux differentiable at $\theta_{\mathcal{S}}$ and the Gâteaux derivative of $P^3$ at $\theta_{\mathcal{S}}$ satisfies that,

$$(P^3)'(\theta_{\mathcal{S}})(u) = 3 \cdot 0 \cdot u = \theta_{\mathcal{S}}, \text{ for every } u \in \mathcal{S}(\mathbb{R}^n, \mathbb{C}) \backslash \{\theta_{\mathcal{S}}\}.$$

This implies that the origin $\theta_{\mathcal{S}}$ is a generalized critical point of $P^3$. Define

$$h(x) = e^{-\frac{1}{x_1^2 + \cdots + x_n^2}}, \text{ for any } x = (x_1, \ldots, x_n) \in \mathbb{R}^n.$$

It is well-known that

$$h \in \mathcal{S}(\mathbb{R}^n, \mathbb{C}) \quad \text{and} \quad h \in K \backslash \{\theta_{\mathcal{S}}\}, \text{ that is } h \succ_K \theta_{\mathcal{S}}.$$

Then, for real number $t$, we have

$$t > 0 \Longrightarrow th \in K \backslash \{\theta_{\mathcal{S}}\}, \text{ that is } th \succ_K \theta_{\mathcal{S}};$$

and

$$t < 0 \Longrightarrow th \in -K \backslash \{\theta_{\mathcal{S}}\}, \text{ that is } th \prec_K \theta_{\mathcal{S}}.$$

By definition, we have

$$P^3(h)(x) = h^3(x) = e^{-\frac{3}{x_1^2+\cdots+x_n^2}}, \text{ for any } x = (x_1, \dots, x_n) \in \mathbb{R}^n.$$

This implies that, for any real $t$, we have

$$P^3(th)(x) = t^3h^3(x) = t^3e^{-\frac{3}{x_1^2+\cdots+x_n^2}}, \text{ for any } x = (x_1, \dots, x_n) \in \mathbb{R}^n.$$

This immediately implies that

$$t > 0 \implies P^3(th) \in K\backslash\{\theta_S\}, \text{ that is } P^3(th) \succ_K \theta_S;$$

and
$$t < 0 \implies P^3(th) \in -K\backslash\{\theta_S\}, \text{ that is } P^3(th) \prec_K \theta_S.$$

Hence, the origin $\theta_S$ is not an $\preccurlyeq_K$-extreme of the 3rd power operator $P^3$ on $\mathcal{S}(\mathbb{R}^n, \mathbb{C})$. □

**Example 9.11**. We continue Example 9.9. Let $(\mathbb{S}, \tau_\mathbb{S})$ be equipped with the partial order $\preccurlyeq_K$ given in Example 9.9. By Theorem 7.3, the 3rd power operator $R^3$ on $(\mathbb{S}, \tau_\mathbb{S})$ is Gâteaux differentiable at $\theta_\mathbb{S}$ such that, and the Gâteaux derivative $(R^3)'(\theta_\mathbb{S})$ is given by

$$(R^3)'(\theta_\mathbb{S}) = \nabla(R^3)(\theta_\mathbb{S}) = \begin{pmatrix} 0 & 0 & \cdots \\ 0 & 0 & 0 \\ \vdots & 0 & \ddots \end{pmatrix}.$$

This satisfies that

$$(R^3)'(\theta_\mathbb{S})(v) = \{s_n\}_{n=1}^{\infty}\begin{pmatrix} 0 & 0 & \cdots \\ 0 & 0 & 0 \\ \vdots & 0 & \ddots \end{pmatrix} = \theta_\mathbb{S}, \text{ for every } v = \{s_n\}_{n=1}^{\infty} \in \mathbb{S}\backslash\{\theta_\mathbb{S}\}.$$

Hence, $\theta_\mathbb{S}$ is a generalized critical point of the 3rd power operator $R^3$ in $(\mathbb{S}, \tau_\mathbb{S})$. Let $w = \{1\}_{n=1}^{\infty} \in \mathbb{S}\backslash\{\theta_\mathbb{S}\}$. It is clear that $w \in K\backslash\{\theta_\mathbb{S}\}$, that is, $w \succ_K \theta_\mathbb{S}$. Similarly to the proof of Example 9.8, we have that

$$t > 0 \implies R^3(tw) = \{t^3\}_{n=1}^{\infty} \in K\backslash\{\theta_\mathbb{S}\}, \text{ that is } R^3(tw) \succ_K \theta_\mathbb{S};$$

and
$$t < 0 \implies R^3(tw) = \{t^3\}_{n=1}^{\infty} \in -K\backslash\{\theta_\mathbb{S}\}, \text{ that is } R^3(tw) \prec_K \theta_\mathbb{S}.$$

Hence, the origin $\theta_\mathbb{S}$ is not a $\preccurlyeq_K$-extreme of the 3rd power operator $R^3$ on $(\mathbb{S}, \tau_\mathbb{S})$. □

### 9.3. Order Monotonic Mappings in Partially Ordered Topological Vector Spaces

In calculus, the monotone property of a differentiable real function on an open interval is described by positive or negative derivatives of the considered function. In this subsection, we extend this property to single-valued mappings in general topological vector spaces.

In this subsection, if not otherwise stated, we always let $(X, \tau_X, \preccurlyeq_C)$ and $(Y, \tau_Y, \preccurlyeq_K)$ be partially ordered Hausdorff topological vector spaces. Let $A$ be a nonempty convex $\tau_X$-open subset in $X$ and let $T: A \to Y$ be a single-valued mapping.

**Definition 9.12.** If for any $x_1, x_2 \in A$, we have

$$x_1 \preccurlyeq_C x_2 \implies T(x_1) \preccurlyeq_K T(x_2),$$

then, $T$ is said to be $\preccurlyeq_C$-$\preccurlyeq_K$ increasing on $A$, or $T$ is said to be order increasing on $A$, if the orders are known. Order decreasing can be similarly defined. In particular, when $(Y, \tau_Y, \preccurlyeq_K) = (\mathbb{R}, \leq)$, then $T$ is said to be increasing on $A$.

**Theorem 9.13**. *Suppose that $T$ is Gâteaux differentiable on $A$. If $T$ is order increasing on $A$, then for any $\bar{x} \in A$, we have*

$$T'(\bar{x})(v) \succcurlyeq_K \theta_Y, \text{ for any } v \succ_C \theta_X \text{ (That is } v \in C\backslash\{\theta_X\}).$$

*Proof*. Suppose that $T$ is order increasing on $A$. Since $A$ is a nonempty convex $\tau_X$-open subset in $X$, for the given $\bar{x} \in A$ and $v \in C\backslash\{\theta_X\}$, by the continuity of addition operator: $X \times X \to X$ and scalar multiplication operator: $\mathbb{R} \times X \to X$, there is $\delta_0 > 0$ such that

$$|t| < \delta_0 \quad \Longrightarrow \quad \bar{x} + tv \in A.$$

Notice that, by $v \succ_C \theta_X$, we have that

$$t > 0 \;\; \Longrightarrow \;\; \bar{x} + tv \succcurlyeq_C \bar{x} \quad \text{and} \quad t < 0 \;\; \Longrightarrow \;\; \bar{x} + tv \preccurlyeq_C \bar{x}.$$

By the condition that $T$ is order increasing on $A$, this implies that

$$\frac{T(\bar{x}+tv)-T(\bar{x})}{t} \succcurlyeq_K \theta_Y, \text{ for any } t \text{ with } 0 < |t| < \delta_0. \tag{9.18}$$

By (9.18), similar to the proofs of Theorem 9.3 and Theorem 4.2. we will show that $T'(\bar{x})(v) \succcurlyeq_K \theta_Y$, this is, $T'(\bar{x})(v) \in K$. Assume, by contradiction, that $T'(\bar{x})(v) \notin K$, this is equivalent to $T'(\bar{x})(v) \in Y\backslash K$. Notice that $V(\mathbb{F}_Y)(T'(\bar{x})(v))$ forms an $\tau_Y$-open neighborhood basis around $T'(\bar{x})(v)$ in $Y$ with

$$V(\mathbb{F}_Y)(T'(\bar{x})(v)) = \{V_{J,\lambda}(T'(\bar{x})(v)): J \in \mathcal{F}_Y, \lambda > 0\}.$$

Since $Y\backslash K$ is an $\tau_Y$-open subset of $Y$, by $T'(\bar{x})(v) \in Y\backslash K$, there is an $\tau_Y$-open neighborhood (It is defined in (9.1)) $V_{J,\varepsilon}(T'(\bar{x})(v))$ of $T'(\bar{x})(v)$ in $Y$ for some $J \in \mathcal{F}_Y$ and $\varepsilon > 0$ such that $V_{J,\varepsilon}(T'(\bar{x})(v)) \subseteq Y\backslash K$. This induces

$$\left(V_{J,\varepsilon}(T'(\bar{x})(v))\right) \cap K = \emptyset. \tag{9.19}$$

Notice that by definition (9.1), we have

$$V_{J,\varepsilon}\big(T'(\bar{x})(v)\big) = \{y \in Y: \max\{q\big(y - T'(\bar{x})(v)\big): q \in J\} < \varepsilon\}. \tag{9.20}$$

(9.19) and (9.20) together imply that

$$\min\{q\big(y - T'(\bar{x})(v)\big): q \in J\} \geq \varepsilon, \text{ for any } y \in K. \tag{9.21}$$

Since $T$ is Gâteaux differentiable at point $\bar{x}$, then $T$ is Gâteaux directional differentiable at point $\bar{x}$ along the given (fixed) direction $v$. By definition (4.1) in Definition 4.1, for the above given $\varepsilon$-open neighborhood $V_{J,\varepsilon}(T'(\bar{x})(v))$ of $T'(\bar{x})(v)$ in $Y$, there is a positive number $\delta_1 < 1$ such that, for real number $t$, we have

$$0 < |t| < \delta_1 \quad \Longrightarrow \quad \mathrm{Max}\left\{q\left(\frac{T(\bar{x}+tv)-T(\bar{x})}{t} - T'(\bar{x})(v)\right): q \in J\right\} < \varepsilon. \tag{9.22}$$

By the conditions in this theorem that $A$ is a convex $\tau_X$-open subset in $X$ and $\bar{x} \in A$, and by the property that $U(\mathbb{F}_X)(\bar{x})$ forms an $\tau_X$-open neighborhood basis of $X$ around $\bar{x}$, there is an $\tau_X$-open neighborhood $U_{I,\delta_2}(\bar{x})$ of $\bar{x}$ in $X$, for some $I \in \mathcal{F}_X$ and $0 < \delta_2 < \delta_1$ such that $U_{I,\delta_2}(\bar{x}) \subseteq A$ satisfying

$$U_{I,\delta_2}(\bar{x}) = \{x \in X : \max\{p(x - \bar{x}) : p \in I\} < \delta_2\}.$$

By the properties (ii) and (iii) in Lemma 2.5 of $F$-seminorms on $X$, for this $v \in C\backslash\{\theta_X\}$, we have

$$\max\{p(tv) : p \in I\} \to 0 \text{ as } t \to 0. \tag{9.23}$$

Then, for the above fixed $\delta_2 > 0$, by (9.23), there is $\delta_3 > 0$ with $0 < \delta_3 < \delta_2$, such that

$$|t| < \delta_3 \implies \max\{p(tv) : p \in I\} < \delta_2 \implies \bar{x} + tv \in U_{I,\delta_2}(\bar{x}) \subseteq A. \tag{9.24}$$

(9.24) implies that

$$|t| < \delta_3 \implies T(\bar{x} + tv) \text{ is well-defined.} \tag{9.25}$$

Let $\delta = \min\{\delta_1, \delta_3\}$. By (9.25) and (9.22) we obtain that

$$0 < |t| < \delta \implies \operatorname{Max}\left\{q\left(\frac{T(\bar{x}+tv)-T(\bar{x})}{t} - T'(\bar{x})(v)\right) : q \in J\right\} < \varepsilon. \tag{9.26}$$

By (9.18), it yields that

$$0 < |t| < \delta \implies \frac{T(\bar{x}+tv)-T(\bar{x})}{t} \in K. \tag{9.27}$$

This creates a contradiction between (9.27) and (9.22). Hence $T'(\bar{x}, v) \in K$. This proves

$$T'(\bar{x}, v) \succcurlyeq_K \theta_Y. \qquad \square$$

**Corollary 9.14**. *Let $(X, \|\cdot\|_X, \preccurlyeq_C)$ be a partially ordered normed vector space. Let $(Y, \tau_Y, \preccurlyeq_K)$ be a partially ordered Hausdorff topological vector space. Let A be a nonempty convex $\|\cdot\|_X$-open subset in X. Let $T: A \to Y$ be a single-valued mapping. Suppose that T is Fréchet differentiable on A. If T is order increasing on A, then for any $\bar{x} \in A$, we have*

$$\nabla T(\bar{x})(v) = T'(\bar{x})(v) \succcurlyeq_K \theta_Y, \textit{ for any } v \succ_C \theta_X \textit{ (That is } v \in C\backslash\{\theta_X\}).$$

*Proof.* By using Theorem 9.13, the proof of this corollary is similar to the proof of Corollary 9.6. □

**Corollary 9.15**. *Let $(X, \tau_X, \preccurlyeq_C)$ be a partially ordered Hausdorff topological vector space. Let A be a nonempty convex $\tau_X$-open subset in X and let $T: A \to \mathbb{R}$ be a real valued function. Suppose that T is Gâteaux differentiable on A. If T is increasing on A, then for any $\bar{x} \in A$, we have*

$$T'(\bar{x})(v) \geq 0, \textit{ for any } v \succ_C \theta_X \textit{ (That is } v \in C\backslash\{\theta_X\}).$$

*Proof*. In Theorem 9.13, let $(Y, \tau_Y, \preccurlyeq_K) = (\mathbb{R}, \leq)$, which proves this corollary immediately. □

**Corollary 9.16**. *Let $(X, \|\cdot\|_X, \preccurlyeq_C)$ be a partially ordered normed vector space. Let A be a nonempty convex $\|\cdot\|_X$-open subset in X. Let $T: A \to \mathbb{R}$ be a real valued function. Let $\bar{x} \in A$. Suppose that T is Fréchet differentiable on A. If T is increasing on A, then for any $\bar{x} \in A$, we have*

$$\nabla T(\bar{x})(v) = T'(\bar{x})(v) \geq 0, \textit{ for any } v \succ_C \theta_X \textit{ (That is } v \in C\backslash\{\theta_X\}).$$

*Proof*. By Corollaries 9.14 and 9.15, this corollary is proved immediately. □

## 10. Conclusion and Remarks

The cornerstone of this paper is to use Theorem 1.1 that states that for any topological vector space, its topology can be induced by a family of $F$-seminorms and to use the $F$-seminorm structures to develop an extended $\varepsilon$-$\delta$ language to define the continuity and differentiability of mappings between general topological vector spaces. Meanwhile, some related techniques are developed, which are used to prove and to find the explicit Gâteaux derivatives and Fréchet derivatives of some operators, which include continuous and linear operators and polynomial type operators.

We showed that the definitions of Gâteaux directional differentiability and Gâteaux differentiability (Definition 4.1) and Fréchet differentiability (Definition 4.8) in general topological vector spaces are natural extensions of the differentiability in normed vector (Banach) spaces.

It is worth to note again that, under our definitions, in the considered topological vector spaces $\mathcal{S}(\mathbb{R}^n, \mathbb{C})$ in section 5, $(\sigma_\rho, \tau_\rho)$ in section 6 and $(\mathbb{S}, \tau_{\mathbb{S}})$ in section 7, for the power type operators, their Fréchet derivatives and Gâteaux derivatives are the same. It is known that by (1.3), in normed vector spaces, if an operator is Fréchet differentiable at a given point, then it is Gâteaux differentiable at this point and the Fréchet derivatives and Gâteaux derivatives are exactly same.

There are some questions for interested readers.

1. We study the connection between Gâteaux (directional) differentiability and Fréchet differentiability of mappings between general topological vector spaces. Let $(X, \tau_X)$ and $(Y, \tau_Y)$ be Hausdorff topological vector spaces. Let $T\colon X \to Y$ be a single-valued mapping. Let $\bar{x} \in X$. Does, in general,

$$T \text{ is Fréchet differentiable at } \bar{x} \quad \Longrightarrow \quad T \text{ is Gâteaux differentiable at } \bar{x}?$$

2. In section 4, contrast with ordinary differentiation in calculus, several properties of Fréchet derivatives were proved. Is it possible to prove a certain type of chain rule for Fréchet derivatives in general topological vector spaces?
3. Let $(X, \tau_X)$, $(Y, \tau_Y)$ and $(Z, \tau_Z)$ be Hausdorff topological vector spaces, in which the topologies $\tau_X$, $\tau_Y$ and $\tau_Z$ are induced by families $\mathbb{F}_X$, $\mathbb{F}_Y$ and $\mathbb{F}_Z$ of positive $F$-seminorms, respectively. The families $\mathbb{F}_X$, $\mathbb{F}_Y$ and $\mathbb{F}_Z$ induce the corresponding $F$-seminorm bases $\{U_{I,\lambda}\colon \lambda > 0, I \in \mathcal{F}_X\}$, $\{V_{J,\lambda}\colon \lambda > 0, J \in \mathcal{F}_Y\}$ and $\{W_{K,\lambda}\colon \lambda > 0, K \in \mathcal{F}_Z\}$ in $X$, $Y$ and $Z$, respectively. Let $\theta_X$, $\theta_Y$ and $\theta_Z$ respectively denote the origins of $X$, $Y$ and $Z$. Let $T\colon X \to Y$ and $S\colon Y \to Z$ be single-valued mappings. Let $\bar{x} \in X$. Can we prove the following chain-rules?

   *Suppose $T$ is Fréchet differentiable at $\bar{x}$ and $S$ is Fréchet differentiable at $T(\bar{x})$, which means that both $\nabla T(\bar{x})\colon X \to Y$ and $\nabla S(T(\bar{x}))\colon Y \to Z$ exist. Then $S \circ T\colon X \to Z$ is Fréchet differentiable at $\bar{x}$ such that*

$$\nabla(S \circ T)(\bar{x}) = \nabla S(T(\bar{x})) \circ \nabla T(\bar{x}).$$

4. In theorem 4.2 and Corollary 4.3, we proved the uniqueness of Gâteaux (directional) derivatives for general topologica vector spaces. However, regarding to the uniqueness of Fréchet derivatives. in Theorem 4.10, we proved this property under the following conditions for the space $X$.

(SC) $X$ is seminorm constructed;
(SB) $\mathbb{F}_X$ is boundedness under the following sense

$$\sup\{p(u): p \in \mathbb{F}_X\} < \infty, \text{ for every } u \in X,$$

We consider the conditions (SC) and (SB) to be very strong. So, can we prove the uniqueness of Fréchet derivatives without conditions (SC) and (SB)?

5. Let $(Z, \tau_X)$ and $(Y, \tau_Y)$ be Hausdorff topological vector spaces, in which there are two families $\mathbb{F}_X^1$ and $\mathbb{F}_X^2$ of positive $F$-seminorms on $X$ such that $\tau_X$ is induced by each one of $\mathbb{F}_X^1$ and $\mathbb{F}_X^2$, independently. Meanwhile, there are two families $\mathbb{F}_Y^1$ and $\mathbb{F}_Y^2$ of positive $F$-seminorms on $Y$ such that $\tau_Y$ is induced by each one of $\mathbb{F}_Y^1$ and $\mathbb{F}_Y^2$, independently. In Theorem (4.23), we proved that, if space $X$ satisfies condition (SX) and $Y$ satisfies condition (SY), then

   $T$ is Fréchet differentiable at $\bar{x}$ with $\tau_X$ and $\tau_Y$ being induced by $\mathbb{F}_X^1$ and $\mathbb{F}_Y^1$, respectively,
   $\Longrightarrow$ $T$ is Fréchet differentiable at $\bar{x}$ with $\tau_X$ and $\tau_Y$ being induced by $\mathbb{F}_X^2$ and $\mathbb{F}_Y^2$, respectively, that have the same Fréchet derivative.

   Can we prove the above results (Theorem 4.23) without the conditions (SX) and (SY)?

6. Let $\rho \in (0,1)$ and let $\mathcal{L}_\rho[0, 1]$ denote the vector space of Lebesgue measurable real valued functions defined on [0, 1] satisfying that

$$\|f\|_\rho := \int_0^1 |f(t)|^\rho \, dt < \infty, \text{ for any } f \in \mathcal{L}_\rho[0, 1].$$

   It is known that $\|\cdot\|_\rho$ is a $F$-norm on $\mathcal{L}_\rho[0, 1]$, and therefore, $(\mathcal{L}_\rho[0, 1], \|\cdot\|_\rho)$ is a $F$-space, which is a topological vector space that the topology is induced by a single $F$-norm $\|\cdot\|_\rho$ on $\mathcal{L}_\rho[0, 1]$. $(\mathcal{L}_\rho[0, 1], \|\cdot\|_\rho)$ is not locally convex.

   Define some single-valued mappings on $\mathcal{L}_\rho[0, 1]$ and investigate their Gâteaux (directional) differentiability and Fréchet differentiability.

7. Let $\gamma \in (0,1]$ and let $\mathcal{M}_\gamma[0, 1]$ denote the vector space of Lebesgue measurable real valued functions defined on [0, 1] satisfying that $q_\gamma: \mathcal{M}_\gamma[0, 1] \to \mathbb{R}_+$ and

$$q_\gamma(f) := \int_0^1 \frac{|f(t)|^\gamma}{1+|f(t)|^\gamma} \, dt, \text{ for any } f \in \mathcal{M}_\gamma[0, 1].$$

   Then, one could check that $q_\gamma$ is a $F$-norm on $\mathcal{M}_\gamma[0, 1]$, and therefore, $(\mathcal{M}_\gamma[0, 1], q_\gamma)$ is a $F$-space, which is a topological vector space that its topology is induced by a single $F$-norm $q_\gamma$ on $\mathcal{M}_\gamma[0, 1]$. $(\mathcal{M}_\gamma[0, 1], q_\gamma)$ is not locally convex.

   Define some single-valued mappings on $\mathcal{M}_\gamma[0, 1]$ and investigate their Gâteaux (directional) differentiability and Fréchet differentiability.